\newif{\ifarxiv}
\arxivtrue

\documentclass{article}
\pdfoutput=1 
\usepackage[T1]{fontenc}
\usepackage{amsmath,amsfonts,amssymb,enumitem,url,graphicx,orcidlink}
\usepackage{bm}
\usepackage{prooftree}

\DeclareRobustCommand{\VAN}[3]{#2} 

\newtheorem{theorem}{Theorem}[section]
\newtheorem{lemma}[theorem]{Lemma}
\newtheorem{proposition}[theorem]{Proposition}
\newtheorem{corollary}[theorem]{Corollary}
\newtheorem{fact}[theorem]{Fact}
\newtheorem{definition}[theorem]{Definition}

\newcommand\eqdef{{\buildrel\text{def}\over=}}

\newtheorem{remark}[theorem]{Remark}
\newtheorem{example}[theorem]{Example}
\numberwithin{equation}{section}

\newcommand\proof{\emph{Proof.}\quad}
\newcommand\qed{\hfill $\Box$\vskip0.2em}

\newcommand\card{\mathrm{card}\;}
\newcommand\sob{\mathrm{sob}\;}
\newcommand\rsob{\mathrm{rsob}\;}

\newcommand\Sober{\mathcal S}
\newcommand\Hoare{\mathcal H}
\newcommand\V{{\mathsf V}}
\newcommand\HV{\Hoare_\V}
\newcommand\HVz{\Hoare_{0 \V}}
\newcommand\fin{{\mathrm{fin}}}
\newcommand\Open{\mathcal O}
\newcommand\real{\mathbb{R}}

\newcommand\nat{\mathbb{N}}
\newcommand\pow{\mathbb{P}}
\newcommand\dc{\mathop{\downarrow}}
\newcommand\upc{\mathop{\uparrow}}
\newcommand\limp{\mathrel{\Rightarrow}}
\newcommand\rk{\mathrm{rk}}
\newcommand\diff{\smallsetminus}
\newcommand\lex{{\mathrm{lex}}}
\newcommand\Step{\mathrm{Step}}
\newcommand\Inc{\mathrm{Inc}}

\newcommand\Top{\mathbf{Top}}

\newcommand\tA{\mathtt{A}}

\newcommand\tP{\mathtt{P}}
\newcommand\tQ{\mathtt{Q}}
\newcommand\tR{\mathtt{R}}
\newcommand\bA{\bm{A}}
\newcommand\bB{\bm{B}}
\newcommand\bC{\bm{C}}

\newcommand\bP{\bm{P}}
\newcommand\bQ{\bm{Q}}
\newcommand\bR{\bm{R}}
\newcommand\bU{\bm{U}}
\newcommand\bV{\bm{V}}

\newcommand\adjsw[1]{{#1}^\circ}

\newcommand\adjmul[1]{\widehat{#1}}

\newcommand{\bigtriangleright}{\mathop{\mbox{\LARGE$\triangleright$}}}

\newcommand\mopen{\{\mkern-\thinmuskip|}
\newcommand\mclose{|\mkern-\thinmuskip\}}
\newcommand\mempty{\pmb{\emptyset}}
\newcommand\mpl{\cdot} 

\newcommand\red\relax
\newcommand\blue\relax
\newcommand\purple\relax

\begin{document}

\renewcommand{\bf}{\bfseries}
\renewcommand{\sc}{\scshape}

\title{Statures and Sobrification Ranks of Noetherian Spaces}

\author{Jean Goubault-Larrecq\thanks{The first author was supported by grant ANR-17-CE40-0028 of the
    French National Research Agency ANR (project BRAVAS)} $^1$
  and Bastien Laboureix$^{1,2}$\\[0.8em]
  $^1$ Universit\'e Paris-Saclay, CNRS, ENS Paris-Saclay,\\
  Laboratoire M\'ethodes Formelles, 91190, Gif-sur-Yvette, France\\
  \url{goubault@lsv.fr}
  \orcidlink{0000-0001-5879-3304}\\[0.5em]
  $^2$ Universit\'e de Lorraine, LORIA, 54000 Nancy, France\\
  \url{bastien.laboureix@hotmail.fr}
}

\maketitle

\begin{abstract}
  There is a rich theory of maximal order types of well-partial-orders
  (wpos), pioneered by de Jongh and Parikh (1977) and Schmidt (1981).
  Every wpo is Noetherian in its Alexandroff topology, and there are
  more; this prompts us to investigate an analogue of that theory in
  the wider context of Noetherian spaces.

  The notion of maximal order type does not seem to have a direct
  analogue in Noetherian spaces per se, but the equivalent notion of
  stature, investigated by Blass and Gurevich (2008) does: we define
  the stature $||X||$ of a Noetherian space $X$ as the ordinal rank of
  its poset of proper closed subsets.  We obtain formulas for statures
  of sums, of products, of the space of words on a space $X$, of the
  space of finite multisets on $X$, in particular.  They confirm
  previously known formulas on wpos, and extend them to Noetherian
  spaces.

  The proofs are, by necessity, rather different from their wpo
  counterparts, and rely on explicit characterizations of the
  sobrifications of the corresponding spaces, as obtained by Finkel
  and the first author (2020).

  We also give formulas for the statures of some natural Noetherian
  spaces that do not arise from wpos: spaces with the cofinite
  topology, Hoare powerspaces, powersets, and spaces of words on $X$
  with the so-called prefix topology.

  Finally, because our proofs require it, and also because of its
  independent interest, we give formulas for the ordinal ranks of the
  sobrifications of each of those spaces, which we call their
  sobrification ranks.


\end{abstract}
\ifarxiv
\relax
\else
\tableofcontents
\fi

\section{Introduction}
\label{sec:intro}

A \emph{well-quasi-ordering} $\leq$ (or \emph{wqo}) on a set $P$ is a
preordering such that every sequence ${(x_n)}_{n \in \nat}$ is good,
namely: there are indices $m<n$ such that $x_m \leq x_n$.  Then
$(P, \leq)$, or just $P$ for short, is a \emph{well-quasi-order}; we
abbreviate this as wqo as well.  A \emph{well-partial-order}
(\emph{wpo}) is an antisymmetric wqo.  A wpo is the same thing as a
partial ordering whose linear extensions are all well-founded
\cite{Wolk:wqo}.  In that case, a well-known result by de Jongh and
Parikh states that the supremum of the ordinal ranks of those linear
extensions is attained \cite[Theorem~2.13]{dJP:maxot}; that supremum
is called the \emph{maximal order type} $o (P)$ of the wpo $P$.

\blue{The study of maximal order types of wpos grew out of questions
  in infinite combinatorics, pioneered by de Jongh and Parikh
  \cite{dJP:maxot} and Schmidt \cite{Schmidt:maxot}.  Schmidt found
  the maximal order types of wpos of finite words, and of wpos of
  finite trees (a.k.a., terms of first-order logic).  This has
  profound consequences in proof theory where, roughly speaking, the
  maximal order type of $P$ coincides with the largest
  proof-theoretical ordinal of natural logical theories that cannot
  prove that $P$ is well-partial-ordered, see
  \cite{Simpson:nonprov,RW:kruskal} for example.  This also has
  implications on upper bounds of the lengths of rewrite sequences of
  those rewrite systems that can be shown to terminate using some
  familiar well-founded orderings \cite{Weiermann:fin:kruskal}, and in
  studying the exact complexity of certain decidable, but highly
  complex, problems in computer science, and notably in verification
  \cite{FFSS:ack,SS:higman}.  }

The purpose of this paper is to propose an extension of the theory of
maximal order types of wpos to Noetherian spaces.  A \emph{Noetherian
  space} is a topological space in which every open subset is compact,
meaning that each of its open covers has a finite subcover; we do not
require separation. (See Section~9.7 of \cite{JGL:topo} for more
information on Noetherian spaces.) The basic premise of this work is
that a preordered set $(P, \leq)$ is wqo if and only if $P$ is
Noetherian in the Alexandroff topology of $\leq$
\cite[Proposition~9.7.17]{JGL:topo}.  \blue{This was already used in
  \cite{JGL-icalp10} in order to generalize some classic wqo-based
  verification algorithms \cite{ABJ:SRE,FS:wsts}.  With the eventual
  aim of understanding the complexity of those algorithms, one may
  think of following a similar route as in the wpo case
  \cite{FFSS:ack,SS:higman}.  A first step in this endeavor---and our
  only concern in this paper---is to find a notion that would
  generalize the notion of maximal order type to all Noetherian spaces
  $X$, and to study it.}
This will be the \emph{stature} of $X$
(Definition~\ref{defn:stature}).

Explicitly, we will define the stature of $X$ as the ordinal rank of
its poset of proper closed subsets, imitating a characterization of
the maximal order type due to---and taking the name from---Blass and
Gurevich \cite{BG:stature}. A related notion is the ordinal rank of
the poset of \emph{irreducible} closed subsets of $X$, which we will
call the \emph{sobrification rank} of $X$. We will see that the
stature and sobrification rank of a Noetherian space are related, and
this will help us in proofs of formulas for the stature of certain
Noetherian spaces. Additionally, this notion of sobrification rank has
independent interest, even on wpos: the ideal Karp-Miller procedure of
\cite{BFG-lmcs20} requires, among other conditions, that the
sobrification rank of the state space of the well-structured
transition system under study is strictly less than $\omega^2$ (see
Section~5 of that paper).


\subsection{Outline}
\label{sec:outline}

Section~\ref{sec:prelim} provides some preliminary notions, mostly on
ordinals and on topology.
\blue{We explore a few ways of defining an adequate generalization of
  statures to Noetherian spaces in Section~\ref{sec:stature}.  We
  settle on the definition we gave above, and we show that it
  coincides with the maximal order type of $X$ (Blass and Gurevich's
  own stature) in the case where $X$ is a
  wpo.  
  On wpos, that ordinal rank coincides with the length of a longest
  chain in the poset of closed sets (minus $1$), but one should be
  warned that this coincidence fails in the realm of Noetherian
  spaces, except for countable Noetherian spaces, as we see in the
  same section.}

We introduce our notion of sobrification rank in
Section~\ref{sec:dimension}, and we show that, up to some $+1$ or $1+$
terms, the sobrification rank is smaller than the stature, and the
stature is smaller than $\omega$ to the power of the sobrification
rank.  Section~\ref{sec:bf-few-tools} is devoted to a few technical
tools, which will allow us to compare the statures and sobrification
ranks of spaces $X$ and $Y$ once we have certain maps from $X$ to $Y$,
and to compute the stature of a space $X$ from the statures of a
cofinal family of proper closed subsets of $X$.

We compute the statures and sobrification ranks of finite $T_0$
spaces, and of well-founded chains in various topologies, in
Section~\ref{sec:finite-spaces}.  We do the same for spaces with a
cofinite topology in Section~\ref{sec:spac-witha-cofin}, for
topological sums in Section~\ref{sec:sums}, for lexicographic sums in
Section~\ref{sec:bf-lexic-sums}, for topological products in
Section~\ref{sec:products}, for Hoare powerspaces and powersets in
Section~\ref{sec:hoare-powerspaces}, for spaces of finite words with
the so-called word topology in Section~\ref{sec:bf-finite-words}
(generalizing the case of wpos of words explored by de Jongh and
Parikh \cite{dJP:maxot} and Schmidt \cite{Schmidt:maxot}), for spaces
of so-called heterogeneous words in the prefix topology in
Section~\ref{sec:bf-heter-words}, and for spaces of finite multisets
in Section~\ref{sec:bf-finite-multisets} (generalizing the case of
wpos of multisets explored by Aschenbrenner and Pong
\cite{AP:monomial}, Weiermann \cite{Weiermann:multiset} and van der
Meeren, Rathjen and Weiermann \cite{vdMRW:bigVeblen}).  \blue{We have
  chosen this particular list of constructions since it matches
  roughly those introduced \cite[Section~9.7]{JGL:topo} and further
  studied in \cite{FGL:partI}.  Finite words are particularly useful
  in several applications in verification, where they are typically
  used to model the contents of communication queues; see
  \cite[Section~6]{JGL-icalp10} for example.}

In all cases, we obtain exact formulae for stature and sobrification rank,
except for the sobrification rank of spaces of multisets and for the stature of
Hoare powerspaces, for which we obtain non-matching lower
and upper bounds; but those bounds are optimal, as we will demonstrate.

\blue{The formulae that we obtain for statures of constructions that
  already exist in the world of wqos, namely for finite sets, sums,
  products, finite words, and finite multisets, are exactly the same
  as in wqo case, and no surprise should be expected there.  This is a
  good thing: the notions and results we obtain smoothly extend those
  already known on wqos---though with pretty different proofs, as in
  the cases of finite words and finite multisets.}

\blue{We obtain new formulae for constructions that make sense on
  Noetherian spaces but not on wqos.}  Spaces with a cofinite
topology, Hoare powerspaces, powersets, and spaces of words with the
prefix topology are examples of Noetherian spaces that do not arise
from wpos, and therefore form proper generalizations of the theory of
wpos.  One should also note that our results on statures on spaces
obtained from constructions that are classical on wqos, such as $X^*$,
\blue{despite yielding the same formulae,} are not consequences of the
result from the wqo literature that they generalize. Finally, our
study of sobrification ranks seems new as well.

We conclude in Section~\ref{sec:bf-conclusion}.

\section{Preliminaries}
\label{sec:prelim}

Given any function $f$ from a product set $X \times Y$ to $Z$,
the application of $f$ to a pair $(x, y)$ will be written as $f (x,
y)$, not $f ((x, y))$.  The symbols $\subseteq$, $\subsetneq$, and
$\not\subseteq$ stand for inclusion (or equality), strict inclusion,
and the negation of inclusion respectively.

We usually write $\leq$ for the preordering of any preordered space,
and $<$ for its strict part.  A function $f$ is \emph{monotonic} if
and only if $x \leq y$ implies $f (x) \leq f (y)$, and an
\emph{order embedding} if and only if $x \leq y$ and
$f (x) \leq f (y)$ are equivalent, for all points $x$ and $y$.

If all the considered preorderings are orderings, then every order
embedding is injective and monotonic, and every injective monotonic
map is strictly monotonic; we say that $f$ is a \emph{strictly
  monotonic} between posets if and only if $x < y$ implies
$f (x) < f (y)$.

\subsection{Ordinals}

We assume some basic familiarity with ordinals, ordinal sum
$\alpha+\beta$, ordinal multiplication $\alpha\beta$, and ordinal
exponentiation $\alpha^\beta$.  \blue{They are all defined by
  induction on $\beta$, e.g., $\alpha^0 \eqdef 1$,
  $\alpha^{\beta+1} \eqdef \alpha^\beta \alpha$, and
  $\alpha^\beta \eqdef \sup \{\alpha^\gamma \mid \gamma < \beta\}$ for
  every limit ordinal $\beta$; in particular, one should not confuse
  ordinal exponentiation with cardinal exponentiation: if $\alpha$ and
  $\beta$ are countable cardinals, say, then $\alpha^\beta$ is
  countable.  Also, one should beware that ordinal addition and
  ordinal multiplication are not commutative.  We refer to
  \cite{Jech:sets} for more information.}
  
Addition $\alpha+\beta$ is monotonic in $\alpha$, namely
$\alpha\leq \alpha'$ implies $\alpha+\beta \leq \alpha'+\beta$; it is
strictly monotonic and continuous in $\beta$, namely $\beta < \beta'$
implies $\alpha+\beta < \alpha+\beta'$, and
$\alpha + \sup_{i \in I} \beta_i = \sup_{i \in I} (\alpha+\beta_i)$
for every non-empty family ${(\beta_i)}_{i \in I}$ of ordinals.
Similarly, multiplication $\alpha\beta$ is monotonic in $\alpha$ and
$\beta$, continuous in $\beta$, and strictly monotonic in $\beta$ if
$\alpha \neq 0$, and exponentiation $\alpha^\beta$ is monotonic in
$\alpha$ and $\beta$, continuous in $\beta$, and strictly monotonic in
$\beta$ if $\alpha \geq 2$.  Additionally, addition is
left-cancellative: $\alpha+\beta = \alpha+\gamma$ implies
$\beta = \gamma$.

Every ordinal $\alpha$ can be written in a unique way as a finite sum
$\omega^{\alpha_1} + \cdots + \omega^{\alpha_m}$ with
$\alpha \geq \alpha_1 \geq \cdots \geq \alpha_m$, $m \in \nat$.  This
is the so-called \emph{Cantor normal form} of $\alpha$.

The ordinals of the form $\omega^\beta$ are exactly the \emph{additively
indecomposable} ordinals, namely the ordinals $\alpha$ such that any
finite sum of ordinals strictly smaller than $\alpha$ is still
strictly smaller than $\alpha$.  The additively decomposable ordinals
are those whose Cantor normal form is such that $m \neq 1$.

Similarly, the \emph{multiplicatively indecomposable} ordinals, namely
the ordinals $\alpha$ such that any product of ordinals strictly
smaller than $\alpha$ is still strictly smaller than $\alpha$, are
$0$, $1$, $2$, and those of the form $\omega^{\omega^\beta}$.

One can compare ordinals
$\alpha \eqdef \omega^{\alpha_1} + \cdots + \omega^{\alpha_m}$ and
$\beta \eqdef \omega^{\beta_1} + \cdots + \omega^{\beta_n}$ in Cantor
normal form by: $\alpha \leq \beta$ if and only if the list
$\alpha_1, \cdots, \alpha_m$ is lexicographically smaller than or
equal to the list $\beta_1, \cdots, \beta_n$, namely, either the two
lists are equal or there is an index $i$ with $1\leq i\leq \min (m,n)$
such that $\alpha_1=\beta_1$, \ldots, $\alpha_{i-1}=\beta_{i-1}$, and
$\alpha_i < \beta_i$.

The \emph{natural} (or Hessenberg) \emph{sum} $\alpha \oplus \beta$
and the \emph{natural product} $\alpha \otimes \beta$ are defined as
follows, where $\alpha$ and $\beta$ are written in Cantor normal form,
respectively $\omega^{\alpha_1} + \cdots + \omega^{\alpha_m}$ and
$\omega^{\beta_1} + \cdots + \omega^{\beta_n}$: $\alpha \oplus \beta$
is equal to $\omega^{\gamma_1} + \cdots + \omega^{\gamma_{m+n}}$,
where $\gamma_1 \geq \cdots \geq \gamma_{m+n}$ is the list obtained by
sorting the list
$\alpha_1, \cdots, \alpha_m, \beta_1, \cdots, \beta_n$ in decreasing
order; $\alpha \otimes \beta$ is equal to
$\bigoplus_{\substack{1\leq i\leq m\\1\leq j\leq n}} \omega^{\alpha_i
  \oplus \beta_j}$.  Those operations are associative and commutative.
Natural sum is strictly monotonic in both arguments; natural product
is monotonic in each argument, and strictly monotonic in each argument
provided that the other one is non-zero.

The following is an equivalent inductive definition of
$\alpha \oplus \beta$ by well-founded recursion on $(\alpha, \beta)$:
$\alpha \oplus \beta$ is the smallest ordinal strictly larger than all
the ordinals $\alpha' \oplus \beta$ with $\alpha' < \alpha$ and
$\alpha \oplus \beta'$ with $\beta' < \beta$.

One should note that the ordinals of the form $\omega^\alpha$ are
\emph{$\oplus$-indecomposable} as well, namely that any finite natural
sum of ordinals strictly smaller than $\omega^\alpha$ is still
strictly smaller than $\omega^\alpha$. Similarly, the ordinals of the
form $\omega^{\omega^\alpha}$ are \emph{$\otimes$-indecomposable}, in
the sense that any finite natural product of ordinals strictly smaller
than $\omega^{\omega^\alpha}$ is still strictly smaller than
$\omega^{\omega^\alpha}$.

The \emph{ordinal rank} (or \emph{height}) $|P|$ of a well-founded
poset $P$ is the least ordinal $\alpha$ such that there is a strictly
monotonic map from $P$ into $\alpha$, namely into the set of ordinals
strictly smaller than $\alpha$.  This can be defined equivalently as
follows \blue{\cite[Exercise~2G.7]{Moschovakis:descset}}.  By
well-founded induction on $x \in P$, we let $\rk_P (x)$ be the least
ordinal strictly larger than $\rk_P (x')$, for every $x' < x$, namely
$\rk_P (x) \eqdef \sup_{x' < x} (\rk_P (x') + 1)$.  (We take the
supremum of an empty family to be $0$, so the rank of a minimal
element is always $0$.)  Then $|P|$ is defined as
$\sup_{x \in P} (\rk_P (x) + 1)$, or equivalently as
$\rk_{P^\top} (\top)$, where $P^\top $ is $P$ with a fresh element
$\top$ added above all others.

One can show that for every ordinal $\alpha < |P|$, there is an
element $x \in P$ such that $\rk_P (x) = \alpha$, and that for every
$x \in P$ and for every $\alpha < \rk_P (x)$, there is an element
$y \in P$ such that $y < x$ and $\rk_P (y) = \alpha$
\cite[Lemma~13]{BG:stature}.

We allow ourselves to write $\alpha-1$ for the unique ordinal of which
$\alpha$ is the successor, if $\alpha$ is a successor ordinal.  For
example, if $P$ is a well-founded poset with a largest element, then
$|P|-1$ makes sense.  We will also use the notation $\alpha-1$ for
zero and limit ordinals, in which case we agree that $\alpha-1$ simply
denotes $\alpha$ itself.

A \emph{chain} $D$ in a poset $P$ is a non-empty, totally ordered
subset.  When $P$ is well-founded, $D$ is isomorphic to a unique
non-zero ordinal, which happens to be $|D|$, and which is called the
\emph{length} of the chain $D$.  The \emph{chain length} $\ell (P)$ of
a well-founded poset $P$ is the supremum of the lengths of its chains.
It is fairly easy to see that $\ell (P) \leq |P|$ for every
well-founded poset $P$, but the reverse inequality fails in general.
We will give a counterexample as a special case of
Proposition~\ref{prop:l!=rk} below, see also the paragraph before
Theorem~1 in \cite{Schmidt:height}, or the end of Section~3 in
\cite{Kriz:orderperfect}.

\subsection{Topology}

We refer the reader to \cite{JGL:topo} for basic notions of topology.
Some information on Noetherian spaces can be found in Section~9.7 of
that book.  Additional information can be found in \cite{FGL:partI}.

We write $\Open X$ for the lattice of open subsets of a topological
space $X$.

Every topological space $X$ has a \emph{specialization preordering},
which we will always write as $\leq$, and defined by $x \leq y$ if and
only if every open neighborhood of $x$ contains $y$. The closure of a
point $x$ is always equal to its downward closure
$\dc x \eqdef \{x' \in X \mid x' \leq x\}$ in its specialization
preordering, and therefore $x \leq y$ if and only if $x$ belongs to
the closure of $\{y\}$. $X$ is $T_0$ if and only if $\leq$ is a
partial ordering. The specialization preordering of any subspace $F$
of $X$ is the restriction $\leq_{|F}$ of $\leq$ to $F$.

There are several topologies on a preordered set $(P, \leq)$ whose
specialization preordering is $\leq$. The finest one is the
\emph{Alexandroff topology}, whose open subsets are exactly the
\emph{upwards-closed} subsets, namely the subsets $U$ such that every
element larger than an element already in $U$ is itself in $U$. The
coarsest one is the \emph{upper topology}, which is the coarsest one
such that the downward closure $\dc x$ of any point, with respect to
$\leq$, is closed. Its closed subsets are the intersections of \emph{finitary
  closed} subsets, namely of subsets of the form
$\dc \{x_1, \ldots, x_n\} \eqdef \bigcup_{i=1}^n \dc x_i$. An
intermediate topology is the \emph{Scott topology}, which is of
fundamental importance in computer science and elsewhere; its open
sets are the upwards-closed subsets $U$ such that every directed
family $D$ that has a supremum in $U$ already intersects $U$.

The \emph{extended Hoare powerdomain} $\Hoare_0 X$ of a topological
space $X$ is the poset of closed subsets of $X$, ordered by inclusion.
It turns out that $X$ is Noetherian if and only if $\Hoare_0 X$ is
well-founded, namely if and only if there is no strictly descending
chain of non-empty closed subsets of $X$.  This is a rephrasing of
Proposition~9.7.6 of \cite{JGL:topo}, which states it in terms of open
sets.

The \emph{Hoare powerdomain} $\Hoare X$ is $\Hoare_0 X$ minus the
empty set, and is a familiar object in domain theory
\cite[Section~IV.8]{GHKLMS:contlatt}.  $\Hoare X$ is also well-founded
if and only if $X$ is Noetherian.

\blue{One should note that we take a rather unusual view at $\Hoare_0$
  and $\Hoare$, and also at the upcoming construction $\Sober$: we
  consider them as mapping topological spaces to \emph{posets}, and
  Noetherian spaces to well-founded posets.  It is certainly possible
  to consider $\Hoare_0 X$ and $\Hoare X$ as topological spaces of
  their own, with suitable topologies, and we will do so in
  Section~\ref{sec:hoare-powerspaces}; but we will want to consider
  ordinal ranks of $\Hoare_0 X$ and $\Hoare X$, and this requires
  seeing the latter as well-founded posets, not spaces.
}

A closed subset $C$ of $X$ is \emph{irreducible} if and only it is
non-empty and, given any two closed subsets $F_1$ and $F_2$ of $X$
whose union contains $C$, one of them already contains $C$.  The
closures of points $\dc x$ are always irreducible closed, and the
space $X$ is \emph{sober} if and only if every irreducible closed
subset of $X$ is the closure of a unique point.  Given any space $X$,
one can form its \emph{sobrification} $\Sober X$: its points are the
irreducible closed subsets of $X$, and its open subsets are the sets
$\diamond U \eqdef \{C \in \Sober X \mid C \cap U \neq \emptyset\}$,
$U \in \Open X$.  The specialization ordering of $\Sober X$ is inclusion.

$\Sober X$ is always a sober space, and the map
$\eta_X \colon x \mapsto \dc x$ is continuous, and a topological
embedding if $X$ is $T_0$.  More is true: the inverse map
$\eta_X^{-1}$ is an order isomorphism between $\Open X$ and
$\Open (\Sober X)$, for every topological space $X$; its inverse is
the map $U \mapsto \diamond U$.  Since a space is Noetherian if and
only if its lattice of open sets has no infinite ascending chain, it
follows that $X$ is Noetherian if and only if $\Sober X$ is.

In a Noetherian space $X$, every closed subset $F$ is the union of
finitely many irreducible closed subsets $C_1$, \ldots, $C_m$.  See
Lemma~3.6 of \cite{FGL:partI}, for example.  This is a fundamental
result, which we will use often.  More generally, this holds if and
only if $X$ is a \emph{FAC space} \cite{JGL-tp19}, namely a space with
no infinite discrete subspace; every Noetherian space is FAC.

By keeping only the maximal elements (with respect to inclusion) in
the list $C_1$, \ldots, $C_m$, we obtain that $C_1$, \ldots, $C_m$ are
exactly the maximal irreducible closed subsets of $X$ included in $F$,
and this list is unique up to permutation.  We call the finitely many
maximal irreducible closed subsets of $X$ included in $F$ the
\emph{components} of $F$.


\section{Stature}
\label{sec:stature}

We will define a notion of stature for Noetherian spaces which, when
applied to wpos equipped with their Alexandroff topology, will
coincide with the usual notion of maximal order type.

A na\" ive attempt to do so would be the following.  The topological
analogue of a poset is a $T_0$ space.  The topological analogue of an
extension of an ordering is a coarser $T_0$ topology, and therefore
the topological analogue of a linear extension (namely, a
\emph{maximal} ordering extending the original ordering) should be a
minimal $T_0$ topology coarser than the original $T_0$ topology.
However, Larson \cite[Example~6]{Larson:mintop} showed that such
minimal $T_0$ topologies may fail to exist.  This is true even in the
realm of Noetherian spaces: Larson's example is $\real$ with its
cofinite topology, and every set is Noetherian in its cofinite
topology.

Instead, we turn to the following notion.

The \emph{stature} $||P||$ of a well-partial order $P$ is defined by
Blass and Gurevich \cite[Definition~28]{BG:stature} as the ordinal
rank of the poset of non-empty bad sequences ordered by reverse
extension, and coincides with the maximal order type of $P$
\cite[Theorem~10]{BG:stature}.  It also coincides with the ordinal
rank of the poset $\mathcal I (P)$ of proper downwards-closed subsets
of $P$, ordered by inclusion \cite[Proposition~31]{BG:stature}.
(Blass and Gurevich call the downwards-closed subsets \emph{ideals},
but we reserve this term to downwards-closed and directed subsets.  A
proper subset of $P$ is any subset of $P$ different from $P$, and that
includes the empty set.)  Those results had already appeared as
Proposition~2.2 in \cite{Kriz:orderperfect}\purple{, but were known even
  before, as a combination of \cite[Proposition~2]{BP:strat}
  and of
  \cite[Theorem~5.3]{PZ:krulldim} 
  (equivalently, of Theorems~III-2.4
  and~III-2.6 in Chapter~I of N. Zaguia's thesis
  \cite{Zaguia:idealchains}).}

The latter definition is the one that extends nicely to Noetherian
spaces.
\begin{definition}[Stature]
  \label{defn:stature}
  The \emph{stature} $||X||$ of a Noetherian space $X$ is
  $|\Hoare_0 X|-1$.
\end{definition}
Note that $|\Hoare_0 X|-1$ is well-defined, because $\Hoare_0 X$ is
well-founded, and has a largest element, namely $X$ itself.

Equivalently, $||X||$ is equal to the rank $\rk_{\Hoare_0 X} (X)$ of
its largest element $X$, or to the ordinal rank
$|\Hoare_0 X \diff \{X\}|$.

\blue{%
  \begin{remark}
    \label{rem:stat}
    Despite the fact that $\Hoare_0 X$ has just one more
    element than $\Hoare X$, it would be wrong to believe that
    $||X|| = |\Hoare X|$.  In general, $1+|\Hoare X| = ||X||+1$, but
    ordinal addition is not commutative.  For example, if $X$ is
    $\nat$ with the Alexandroff topology of its natural ordering, we
    have $|\Hoare X| = |\Hoare_0 X| = \omega+1$, and $||X|| = \omega$.
  \end{remark}
}

\blue{%
  Let us verify that Definition~\ref{defn:stature} generalizes the
  notion of stature (a.k.a.\ maximal order type) of well-partial
  orders to Noetherian spaces.  Given any poset $P$, considered as a
  topological space with its Alexandroff topology, the closed subsets
  of $P$ are exactly its downwards-closed subsets.  Therefore
  $\Hoare_0 P \diff \{P\}=\mathcal I (P)$, and hence, when $P$ is a
  well-partial order, our $||P||$ coincides with Blass and Gurevich's
  $|\mathcal I (P)|$.  Let us state this explicitly.
  \begin{fact}
    \label{fact:stature}
    For every well-partial order $P$, $||P||$ is the maximal order
    type of $P$.
  \end{fact}
}

For a well-partial order $P$, the stature also coincides with the
chain-length of $\mathcal I (P)$, and there is a chain of that precise
length.  This is the meaning of the equality
``$\mu (\dc A)=\lambda (\dc A)$'' in
\cite[Proposition~2.2]{Kriz:orderperfect}: in K\v r\'\i\v z's
notation, and with $A\eqdef P$, $\mu (\dc A)$ denotes the ordinal rank
of $\mathcal I (P)^\top$, and $\lambda (\dc A)$ denotes the chain
length of $\mathcal I (P)^\top$, provided that the supremum that
defines it is attained.  \purple{(See also
  \cite[Theorem~5.3]{PZ:krulldim}.)}  One may wonder whether a similar
coincidence would hold in the realm of Noetherian spaces, namely
whether $\ell (\Hoare_0 X) = ||X||+1$ for every Noetherian space $X$,
and the answer is negative.
\begin{proposition}
  \label{prop:l!=rk}
  For every ordinal $\alpha$, there is a non-empty, sober Noetherian
  space $X_\alpha$ whose stature is at least $\alpha$, and such that
  all the chains in $\Hoare_0 {X_\alpha}$ are countable.  In
  particular, when $\alpha \geq \aleph_1$,
  $\ell (\Hoare_0 {X_\alpha}) < ||X_\alpha||+1$.
\end{proposition}
\proof We build $X_\alpha$ by induction.  This is a slight variant of
a construction mentioned by D. Schmidt \cite[paragraph before
Theorem~1]{Schmidt:height}.

We build $X_\alpha$ as a poset and equip it with the upper topology.
It turns out that the sober Noetherian spaces are exactly the posets
that are well-founded and have properties~T and~W, in the upper
topology of their ordering \cite[Theorem~9.7.12]{JGL:topo}.
Property~T states that the whole set itself is finitary closed, and
Property~W states that any binary intersection $\dc x \cap \dc y$ is
finitary closed.  In a sober Noetherian space, every closed subset is
finitary closed.

We let $X_0$ be a one element set.  For every ordinal $\alpha$, we let
$X_{\alpha+1}$ be $X_\alpha^\top$.  Finally, for every limit ordinal
$\alpha$, we let $X_\alpha$ be
$(\coprod_{\beta < \alpha} X_\beta)^\top$.  In the disjoint sum
$\coprod_{\beta < \alpha} X_\beta$, any pair of elements coming from
$X_\beta$ and $X_\gamma$ with $\beta\neq \gamma$ is incomparable; the
elements coming from the same summand $X_\beta$ are compared as in
$X_\beta$.

It is easy to see that $X_\alpha$ is a (reversed) tree for every
ordinal $\alpha$, namely: there is a largest element (the
\emph{root}), and the upward closure
$\upc x \eqdef \{y \in X_\alpha \mid x \leq y\}$ of any element $x$ of
$X_\alpha$ is a finite chain.  This is proved by induction on
$\alpha$.  Then properties~T and~W are obvious.  $X_\alpha$ is also
clearly well-founded, by induction on $\alpha$ once again.  Therefore
$X_\alpha$ is sober and Noetherian in its upper topology.

Every element $x$ of $X_\alpha$ is the largest element of a poset
isomorphic to $X_\beta$ for some unique ordinal $\beta \leq \alpha$.
By induction on $\beta$, $\rk_{X_\alpha} (x) = \beta$.  Also, since
$\dc x$ is a proper closed subset of $\dc y$ if $x < y$, a similar
induction on $\beta$ yields
$\rk_{\Hoare_0 X_\alpha} (\dc x) \geq \beta$.  Taking
$\beta \eqdef \alpha$, we obtain that
$||X_\alpha|| = \rk_{\Hoare_0 X_\alpha} (X_\alpha)$ is at least
$\alpha$.

We now claim that all the chains in $\Hoare_0 {X_\alpha}$ are
countable.  This is again by induction on $\alpha$.  The only
interesting case is when $\alpha$ is a limit ordinal.  Let $D$ be a
chain in $\Hoare_0 {X_\alpha}$.  $D$ can contain at most one closed
set containing the top element of $X_\alpha$, since there is only one
such set, which is $X_\alpha$ itself.  Without loss of generality, we
will therefore assume that $D$ only contains closed subsets of
$X_\alpha$ that do not contain its top element $\top$.  All those
closed subsets $F$ are finitary, hence must be of the form
$\bigcup_{i=1}^n \dc x_i$, where each $x_i$ is taken from some summand
$X_{\beta_i}$ of the sum $\coprod_{\beta < \alpha} X_\beta$; let us
call the set $\{\beta_1, \cdots, \beta_n\}$ the \emph{support} $s(F)$
of $F$.  For any two closed subsets $F$, $F'$ not containing $\top$,
$F \subseteq F'$ implies $s(F) \subseteq s(F')$.  Therefore the family
$D' \eqdef \{s (F) \mid F \in D\}$ forms a chain of finite subsets of
$\alpha$.  It is easy to see that there can be at most one element of
$D'$ of each given cardinality, so that $D'$ is countable.  It follows
that there is a countable subset $I$ of $\alpha$ such that
$s (F) \subseteq I$ for every $F \in D$.

For every $\beta \in I$, let $D_{|\beta}$ denote the chain consisting
of those sets of the form $F \cap X_\beta$, where $F$ ranges over $D$.
This is a chain of closed subsets of $X_\beta$, which is therefore
countable by induction hypothesis.  We should mention a subtle point
here: those are closed in the subspace topology, but we require to
work with closed subsets in the \emph{upper} topology of $X_\beta$ to
be able to use the induction hypothesis.  They are indeed closed in
the upper topology of $X_\beta$, because they are downward closures of
finitely many points in the downward closed subset $X_\beta$.


Any element $F$ of $D$ can be written as the finite union of the sets
$F \cap X_\beta$, where $\beta$ ranges over $s (F)$.  It follows that
there are at most as many elements of $D$ as there are finite subsets
of the countable disjoint union $\biguplus_{\beta \in I} D_{|\beta}$,
hence that $D$ is countable.

If $\alpha \geq \aleph_1$, finally,
$\ell (\Hoare_0 {X_\alpha}) \leq \aleph_1 < \alpha+1 \leq
||X_\alpha||+1$.  \qed

To the contrary, Proposition~\ref{prop:l=rk} below states that every
Noetherian space $X$ with a \emph{countable} Hoare powerdomain is such
that $\ell (\Hoare_0 X) = ||X||+1$. Since having a countable Hoare
powerdomain may seem like a strange condition, we first note that this
is a familiar property.

\begin{lemma}
  \label{lemma:noeth:count}
  The following properties are equivalent for a Noetherian space $X$:
  \begin{enumerate}[label=(\roman*)]
  \item $X$ is second-countable;
  \item $\Hoare_0 X$ is countable;
  \item $\Sober X$ is countable.
  \end{enumerate}
\end{lemma}
\proof $(i)\limp(ii)$.  if $X$ is second-countable, then every open
subset $U$ can be written as a union of open sets from a countable base
$B$.  Since $X$ is Noetherian, $U$ is compact, so $U$ is already equal
to a finite subunion.  Hence $X$ only has countably many open sets,
and therefore only countably many closed sets.

$(ii)\limp(iii)$.  If $\Hoare_0 X$ is countable, then its subset
$\Sober X$ is countable.

$(iii)\limp(i)$.  Let us assume that $\Sober X$ is countable.  Since
$X$ is Noetherian, every closed set is a finite union of elements of
$\Sober X$, so $X$ has only countably many closed subsets.  \qed

\begin{proposition}
  \label{prop:l=rk}
  For every non-empty second-countable Noetherian space $X$,
  $\ell (\Hoare_0 X) = ||X||+1$, and there is a chain in $\Hoare_0 X$
  of maximal length.
\end{proposition}
\proof \purple{By Theorem~III-1.1 of
  \cite[Chapter~I]{Zaguia:idealchains}, every distributive,
  well-founded lattice of countable rank $\alpha$ contains a chain of
  length $\alpha$.  (A even more general statement appears in
  \cite[Remark~3.4]{Kriz:orderperfect}, but we will not need it.)  We
  apply this to the lattice $\Hoare_0 X$, using
  Lemma~\ref{lemma:noeth:count} in order to guarantee that its ordinal
  rank $||X||+1$ is countable.  \qed}
%

\begin{example}
  \label{exa:Xa}
  Let us look back at the counterexample of
  Proposition~\ref{prop:l!=rk}\blue{: there, we built a space
    $X_\alpha$, and we obtained the strict inequality
    $\ell (\Hoare_0 {X_\alpha}) < ||X_\alpha||+1$ when
    $\alpha \geq \aleph_1$.}  When $\alpha < \aleph_1$, namely when
  $\alpha$ is countable, we do have 
  \blue{equality:}
  $\ell (\Hoare_0 {X_\alpha}) = ||X_\alpha||+1$.  We use
  Proposition~\ref{prop:l=rk} in order to show this.  It suffices to
  observe that $X_\alpha$ is sober, hence isomorphic to
  $\Sober X_\alpha$, and countable; the latter is proved by an easy
  induction of the countable ordinal $\alpha$.
\end{example}

\begin{remark}
  \label{rem:catalog}
  Figure~2 of \cite{FGL:partI} displays a catalogue of Noetherian
  spaces.  (We will deal with most spaces in that list, to the
  exception of spectra of Noetherian rings and spaces of trees.)  This
  starts with a few basic Noetherian spaces, such as any finite poset
  or $\nat$ in their Alexandroff topology, or spectra of Noetherian
  rings.  The latter may fail to be second-countable, but the other
  basic Noetherian spaces are.  The catalogue then proceeds by
  building new Noetherian spaces from old, using Noetherian-preserving
  operators $F$---namely, if $X_1$, \ldots, $X_m$ are Noetherian, then
  so is $F (X_1, \cdots, X_m)$.  For example, any finite product, any
  finite sum of Noetherian spaces is Noetherian, the space $X^*$ of
  finite words over $X$ with a suitable topology (see
  Section~\ref{sec:bf-finite-words}) is Noetherian, and so on.  The
  main results of \cite{FGL:partI} consist in giving explicit
  descriptions of $\Sober (F (X_1, \cdots, X_m))$ from $\Sober {X_1}$,
  \ldots, $\Sober {X_m}$.  For example, $\Sober (X^*)$ consists of
  certain regular expressions called word products over $\Sober (X)$,
  which we will define and use in Section~\ref{sec:bf-finite-words}.
  It is then easy to see that if
  $\Sober (X)$ is countable, then so is $\Sober (X^*)$.  That extends
  to all the operators $F$ of \cite[Figure~2]{FGL:partI}.  As a
  consequence of Lemma~\ref{lemma:noeth:count}, and taking the proviso
  that we only consider spectra of Noetherian rings with countably
  many radical ideals, all the Noetherian spaces displayed there are
  second-countable.  This yields an ample supply of Noetherian spaces
  on which Proposition~\ref{prop:l=rk} applies.
\end{remark}

\section{Sobrification Rank}
\label{sec:dimension}

Since the closed subsets of a Noetherian space $X$ are the finite
unions of elements of $\Sober X$, it will be useful to also consider
the ordinal rank of $\Sober X$.  Here we see $\Sober X$ not as a
topological space, rather as a poset under its specialization
ordering, which is inclusion.  As such, it is a well-founded poset,
since it is included in $\Hoare_0 X$, which is itself well-founded
under inclusion.
\begin{definition}[Sobrification Rank]
  \label{defn:dimension}
  The \emph{sobrification rank} $\sob X$ of a Noetherian space is $|\Sober X|$.
\end{definition}

\begin{lemma}
  \label{lemma:dimension}
  For every non-empty Noetherian space $X$,
  \[
    \sob X = \max_{i=1}^m \rk_{\Sober X} (C_i) + 1,
  \]
  where $C_1$, \ldots, $C_m$ are the components of $X$ ($m \geq 1$).
\end{lemma}
\proof The components of $X$ exist because $X$ is closed in $X$, and
$m \geq 1$ because $X$ is non-empty.  For every $C \in \Sober X$, the
inclusion $C \subseteq X = C_1 \cup \cdots \cup C_m$ together with the
fact that $C$ is irreducible imply that $C$ is included in some $C_i$.
It follows that
$|\Sober X| = \max (\rk_{\Sober X} (C_1)+1, \cdots, \rk_{\Sober X}
(C_m)+1)$, whence the claim.  \qed

In particular, $|\Sober X|$ is a successor ordinal for every non-empty
Noetherian space $X$.  Hence the following definition makes sense.
\begin{definition}[Reduced sobrification rank]
  \label{defn:rdim}
  The \emph{reduced sobrification rank} of a non-empty Noetherian space $X$ is
  $\rsob X \eqdef \sob X - 1$.  This is equal to
  $\max_{i=1}^m \rk_{\Sober X} (C_i)$, where $C_1$, \ldots, $C_m$ are
  the components of $X$ ($m \geq 1$).
\end{definition}

There is a notion of \emph{Krull dimension} of Noetherian spaces,
defined as the supremum of the lengths, minus 1, of (finite non-empty)
chains of elements of $\Sober X$.  The result is a natural number or
infinity.  One can extend this notion and consider the ordinal-valued
chain length $\ell (\Sober X)$.
\begin{remark}
  \label{rem:dimension:dim}
  In contrast to Proposition~\ref{prop:l=rk}, $\ell (\Sober X)$ is in
  general not equal to $|\Sober X| = \sob X$, even when $X$ is
  countable.  The space $X_\alpha$ of Proposition~\ref{prop:l!=rk}
  with $\alpha \eqdef \omega$ is sober, so $\Sober {X_\omega}$ is
  isomorphic to $X_\omega$; $X_\omega$ only has finite chains (of
  arbitrary length), so
  $\ell (X_\omega) = \ell (\Sober {X_\omega}) = \omega$, while
  $|X_\omega| = |\Sober {X_\omega}| = \omega+1$.
\end{remark}

\begin{proposition}
  \label{prop:dim:stature}
  For every Noetherian space $X$,
  \begin{enumerate}
  \item $1+\sob X \leq ||X||+1$, or equivalently $1+\rsob X \leq ||X||$;
  \item $||X||+1 \leq \omega^{\sob X}$;
  \item if $X$ is non-empty, then $||X|| \leq \omega^{\rsob X}
    \otimes n$, where $n$ is the number of components of $X$.
  \end{enumerate}
\end{proposition}
\proof (1) The map $\rk_{\Hoare_0 X}$, once restricted to
$\Sober X \cup \{\emptyset\}$, is a strictly increasing map into the
set of ordinals strictly smaller than $|\Hoare_0 X| = ||X||+1$.  Hence
$|\Sober X \cup \{\emptyset\}| \leq ||X||+1$.  We now observe that
$|\Sober X \cup \{\emptyset\}| = 1+\sob X$, because the ordinal rank
of a disjoint union of two posets $A$ and $B$, where every element of
$A$ is below every element of $B$, is the sum of the ordinal ranks of
$A$ and $B$.

(2) For every $F \in \Hoare_0 X$, let $C_1$, \ldots, $C_m$ be the
components of $F$.  We recall that those are the maximal irreducible
closed subsets of $X$ included in $F$; in particular, they are
pairwise incomparable.  We define $\varphi (F)$ as
$\bigoplus_{i=1}^m \omega^{\rk_{\Sober X} (C_i)}$, and we claim that
$\varphi$ is a strictly increasing map.

Given any two elements $F$, $F'$ of $\Hoare_0 X$, let $C_1$, \ldots,
$C_m$ be the components of $F$, and $C'_1$, \ldots, $C'_n$ be those of
$F'$.  If $F \subseteq F'$, then each $C_i$ is included in $F'$, hence
in some $C'_j$, since $C_i$ is irreducible.  We pick one such $j$ and
call it $f (i)$.  Then we split the natural sum
$\varphi (F) = \bigoplus_{i=1}^n \omega^{\rk_{\Sober X} (C_i)}$ as the
natural sum of the quantities
$\alpha_j \eqdef \bigoplus_{i \in f^{-1} (j)} \omega^{\rk_{\Sober X}
  (C_i)}$, $1\leq j\leq n$.

We note that $\alpha_j \leq \omega^{\rk_{\Sober X} (C'_j)}$, with
equality if and only if there is a unique element $i$ in $f^{-1} (j)$
and $C_i=C'_j$.  Indeed, if $f^{-1} (j) = \{i\}$ and $C_i=C'_j$, then
$\alpha_j = \omega^{\rk_{\Sober X} (C_i)} = \omega^{\rk_{\Sober X}
  (C'_j)}$.  Otherwise, for every $i \in f^{-1} (j)$, $C_i$ is a
proper subset of $C'_j$: if $C_i$ were equal to $C'_j$ for some
$i \in f^{-1} (j)$, then $C_i$ would contain $C_{i'}$ for every
$i' \in f^{-1} (i)$, and since the sets $C_i$ are pairwise
incomparable, that would force $f^{-1} (j)$ to consist of exactly one
element.  Since $C_i \subsetneq C'_j$ for every $i \in f^{-1} (j)$,
$\rk_{\Sober X} (C_i) < \rk_{\Sober X} (C'_j)$, so
$\omega^{\rk_{\Sober X} (C_i)} < \omega^{\rk_{\Sober X} (C'_j)}$.
Since $\omega^{\rk_{\Sober X} (C'_j)}$ is $\oplus$-indecomposable,
$\alpha_j < \omega^{\rk_{\Sober X} (C'_j)}$.

Using the fact that natural sum is strictly monotonic in all of its
arguments, $\varphi (F) = \bigoplus_{j=1}^n \alpha_j$ is then smaller
than or equal to
$\bigoplus_{j=1}^n \omega^{\rk_{\Sober X} (C'_j)} = \varphi (F')$,
with equality if and only if $f^{-1} (j)$ consists of a single element
$i$ and $C_i=C'_j$, for every $j$ such that $1\leq j\leq n$.  The
latter case is equivalent to the fact that $f$ is a bijection, and
that the lists $C_1, \cdots, C_m$ and $C'_1, \cdots, C'_n$ are equal
up to permutation, namely to the equality $F=F'$.

It follows that $\varphi$ is a strictly monotonic map from
$\Hoare_0 X$ to a set of ordinals.  The largest value it takes is
$\varphi (X) = \bigoplus_{i=1}^n \omega^{\rk_{\Sober X} (C_i)}$, where
$C_1$, \ldots, $C_n$ are the components of $X$.  Using the additive
indecomposability of $\omega^{|\Sober X|}$, together with the fact
that $\rk_{\Sober X} (C_i) < |\Sober X|$ (see
Lemma~\ref{lemma:dimension}), we obtain that
$\varphi (X) < \omega^{|\Sober X|}$.

The existence of a strictly monotonic map
$\varphi \colon \Hoare_0 X \to \omega^{|\Sober X|}$ in turn implies
that
$||X||+1 = |\Hoare X| \leq \omega^{|\Sober X|} = \omega^{\sob X}$.

(3) We use the same map $\varphi$ as in item~(2).  Since $\varphi$ is
strictly monotonic, $\varphi (F) \geq \rk_{\Hoare_0 X} (F)$ for every
closed subset $F$ of $X$, by an easy well-founded induction on $F$.
Taking $F \eqdef X$, we obtain that
$||X|| = \rk_{\Hoare_0 X} (X) \leq \varphi (X) = \bigoplus_{i=1}^n
\omega^{\rk_{\Sober X} (C_i)}$, where $C_1$, \ldots, $C_n$ are the
components of $X$.  For each $C_i$,
$\rk_{\Sober X} (C_i) \leq \sob X-1 = \rsob X$, by
Lemma~\ref{lemma:dimension}, whence the claim.  \qed

\blue{While we are on the subject of sobrification, we should mention
  the following.
  \begin{fact}
    \label{fact:stat:sob}
    For every Noetherian space $X$, 
    $||\Sober X|| = ||X||$ and $\sob (\Sober X) = \sob X$.
  \end{fact}
  Indeed, the lattices of open subsets of $X$ and of $\Sober X$ are
  isomorphic through $U \mapsto \diamond U$.  Hence $X$ and $\Sober X$
  have isomorphic lattices of closed subsets ($\Hoare_0 X$ and
  $\Hoare_0 (\Sober X)$); the isomorphism also restricts to one
  between $\Sober X$ and $\Sober (\Sober X)$.}

\section{Direct and inverse images}
\label{sec:bf-few-tools}

We collect a few tools that we will use in order to evaluate the
stature of several kinds of spaces.  We will discover a few others
along the way.  We start with a technical lemma.
\begin{lemma}
  \label{lemma:surj}
  Let $f \colon X \to Y$ be a continuous map.
  \begin{enumerate}
  \item The map $\Sober f \colon \Sober X \to \Sober Y$, where for
    every $C \in \Sober X$, $\Sober f (C)$ is the closure $cl (f (C))$
    of the image $f (C)$ of $C$ by $f$, is a well-defined, monotonic
    map.
  \item If $f$ is surjective, then
    $f^{-1} \colon \Hoare_0 (Y) \to \Hoare_0 (X)$ is injective.
  \item If $f^{-1} \colon \Hoare_0 (Y) \to \Hoare_0 (X)$ is injective,
    then it is an order embedding.
  \end{enumerate}
\end{lemma}
\proof (1) For every $C \in \Sober X$, $\Sober f (C)$ is irreducible
closed, see \cite[Lemma~8.2.42]{JGL:topo} for example.  For a short
argument, if $\Sober f (C)$ is included in the union $F_1 \cup F_2$ of
two closed sets, then $C \subseteq f^{-1} (F_1 \cup F_2) = f^{-1}
(F_1) \cup f^{-1} (F_2)$, and the claim follows from the
irreducibility of $C$.  The monotonicity of $\Sober f$ is clear.

(2) Let $F$ and $F'$ be two closed subsets of $Y$ such that
$f^{-1} (F) = f^{-1} (F')$.  For every $y \in F$, we can write $y$ as
$f (x)$ for some $x \in X$, and then $x$ is in $f^{-1} (F)$.  Since
$f^{-1} (F) = f^{-1} (F')$, $f (x) = y$ is also in $F'$.  We show the
converse inclusion $F' \subseteq F$ similarly.

(3)
We show that $f^{-1} (F) \subseteq f^{-1} (F')$ is
equivalent to $F \subseteq F'$.  Indeed,
$f^{-1} (F) \subseteq f^{-1} (F')$ if and only if
$f^{-1} (F) \cup f^{-1} (F') = f^{-1} (F')$, if and only if
$f^{-1} (F \cup F') = f^{-1} (F')$, if and only if $F \cup F' = F'$,
if and only if $F \subseteq F'$.
\qed


We will say that a continuous map $f \colon X \to Y$ is \emph{Skula
  dense} if and only if $f^{-1} \colon \Hoare_0 (Y) \to \Hoare_0 (X)$
is injective.  By Lemma~\ref{lemma:surj}, every surjective map is
Skula dense.  The name ``Skula dense'' stems from the following
observation, which we present for completeness only, and is a slight
relaxation of \cite[Exercise~V-5.32]{GHKLMS:contlatt}.  The
\emph{Skula topology} \cite{Skula:noeth} on $Y$ (called the
\emph{$b$-topology} there, and sometimes also called the \emph{strong}
topology) has all crescents as basic open sets, where a
\emph{crescent} is a difference $U \diff V$ of two open sets.  This is
a remarkable topology.  For example, a topological space $Y$ is sober
Noetherian if and only if it is compact Hausdorff in its Skula
topology \cite[Theorem~3.1]{Hoffmann:sobr:rem}, and for any space $X$
embedded in a sober space $Y$, the Skula closure of $X$ in $Y$ is
homeomorphic to $\Sober X$ \cite[Proposition~3.4]{keimel08}.
\begin{lemma}
  \label{lemma:Sdense}
  A continuous map $f \colon X \to Y$ between topological spaces is
  Skula dense, in the sense that
  $f^{-1} \colon \Hoare_0 (Y) \to \Hoare_0 (X)$ is injective, if and
  only if the image of $f$ is dense in $Y$ with the Skula topology.
\end{lemma}
\proof If $f^{-1}$ is injective, then we claim that any non-empty
Skula open subset of $Y$ intersects the image of $f$.  Such a Skula
open set must contain a basic non-empty Skula open set $U \diff V$,
where $U$ and $V$ are open in $Y$.  Then
$f^{-1} (U \diff V) \neq f^{-1} (\emptyset)$, since $f^{-1}$ is
injective.  Hence $f^{-1} (U \diff V)$ is non-empty.  We pick
$x \in f^{-1} (U \diff V)$, then $f (x)$ is both in $U \diff V$ and in
the image of $f$.

Conversely, if $f$ has Skula dense image, let $U$ and $V$ be two open
subsets of $Y$ such that $f^{-1} (U) = f^{-1} (V)$.  Then both
$f^{-1} (U \diff V)$ and $f^{-1} (V \diff U)$ are empty, so neither
$U \diff V$ nor $V \diff U$ intersects the image of $f$.  Since the
latter is Skula dense, $U \diff V$ and $V \diff U$ must be empty,
whence $U=V$.  \qed

We remark that, as a consequence, the Skula dense maps between $T_0$
spaces are exactly the epimorphisms in the category of $T_0$ spaces
\cite[Exercise~V-5.33]{GHKLMS:contlatt}.

\begin{lemma}
  \label{lemma:stature:mono}
  Let $f \colon X \to Y$ be a continuous map from a Noetherian space
  $X$ to a topological space $Y$.  If $f$ is Skula dense,
  then $Y$ is Noetherian, and $||Y|| \leq ||X||$.
\end{lemma}
\proof Let $F$, $F'$ be two closed subsets of $Y$.  If
$F \subsetneq F'$, then $f^{-1} (F) \subseteq f^{-1} (F')$, and that
inclusion is strict because $f^{-1}$ is injective.  In particular, any
infinite decreasing sequence in $\Hoare_0 Y$ would be mapped through
$f^{-1}$ to an infinite descending sequence in $\Hoare_0 X$, which is
impossible since $X$ is Noetherian.  Therefore $Y$ is Noetherian.

The map
$\rk_{\Hoare_0 X} X \circ f^{-1} \colon \Hoare_0 Y \to ||X||+1$ is a
strictly increasing map, and this immediately entails that
$|\Hoare_0 Y| = ||Y||+1$ is less than or equal to $||X||+1$, hence
that $||Y|| \leq ||X||$.  \qed


\begin{remark}
  \label{rem:dim:mono}
  It is not the case that, under the assumptions of
  Lemma~\ref{lemma:stature:mono}, $\sob Y \leq \sob X$.
  Consider the poset $X \eqdef \{1, 2, 3\}$ with $1 < 2$ and $3$
  incomparable with both $1$ and $2$, and $Y \eqdef \{1, 2, 3\}$ with
  $1<2<3$.  Equip both with their Alexandroff topologies, and let $f$
  be the identity map.  Then $\sob Y = 3 \not\leq \sob X = 2$.
\end{remark}

A dual statement is as follows.  We will call \emph{full} any map
$f \colon X \to Y$ such that every open subset $U$ of $X$ can be
written as $f^{-1} (V)$ for some open subset $V$ of $Y$\footnote{A
  word on naming.  This notion was called ``almost open'' in
  \cite{JGL:topo}.  However, ``almost open'' is improper, for several
  reasons.  One is that it runs into a conflict with another notion,
  due to Klaus Keimel.  Another one is that an equivalent statement is
  that for every closed subset $C$ of $X$, there is a closed subset
  $D$ of $Y$ such that $C = f^{-1} (D)$, showing that the concept has
  little to do with openness.  In order to fix this, the concept was
  called ``initial'' in \cite[Section~4]{GL:infwords}, since an
  initial morphism in the concrete category of topological spaces is
  precisely a full, continuous map.  However, initiality also implies
  continuity, which fullness does not.  Fullness also evokes some form
  of surjectivity, and indeed $f$ is full if and only if
  $f^{-1} \colon \Open Y \to \Open X$ is surjective.}.  Every full map
$f \colon X \to Y$ is order-reflecting with respect to the
specialization preorderings, namely $f (x) \leq f (y)$ implies
$x \leq y$\footnote{Regarding naming considerations again, we may look
  at preordered sets as certain categories, and then an
  order-reflection is a special case of a full functor, namely one
  that is surjective on homsets.}.  A typical example is given by
topological embeddings, which are full, continuous injective maps.
\begin{lemma}
  \label{lemma:stature:mono:inj}
  Let $f \colon X \to Y$ be a full map from a topological space $X$ to
  a Noetherian space $Y$.
  \begin{enumerate}
  \item $X$ is Noetherian;
  \item $||X|| \leq ||Y||$;
  \item if $f$ is also continuous, then $\sob X \leq \sob Y$;
  \item if $f$ is not just full and continuous, but also Skula dense,
    then $||X||=||Y||$ and $\sob X = \sob Y$.
  \end{enumerate}
\end{lemma}
\proof Let us define $f_* (F)$ as $cl (f (F))$, for every closed
subset $F$ of $X$.  This is the same definition as $\Sober f$, except
that $f$ is no longer assumed to be continuous.

The map $f_*$ is monotonic, and we claim that it is injective.  Let us
assume that $f_* (F) = f_* (F')$, where $F$ and $F'$ are closed in
$X$.  For every open subset $U$ of $X$, we write $U$ as $f^{-1} (V)$
for some open subset $V$ of $Y$.  Then $U$ intersects $F$ if and only
if $V$ intersects $f (F)$.  An open set intersects a set $A$ if and
only if it intersects its closure $cl (A)$, so $U$ intersects $F$ if
and only if $V$ intersects $cl (f (F)) = f_* (F)$.  Similarly, $U$
intersects $F'$ if and only if $V$ intersects $f_* (F')$.  Since
$f_* (F) = f_* (F')$, $F$ and $F'$ intersect the same open subsets $U$
of $X$.  Therefore, they are equal.

Since $f_*$ is monotonic and injective, it is strictly monotonic.
Every infinite decreasing sequence in $\Hoare_0 X$ would be mapped by
$f_*$ to an infinite decreasing sequence in $\Hoare_0 Y$, showing (1).
Additionally,
$\rk_{\Hoare_0 Y} \circ f_* \colon \Hoare_0 X \to ||Y||+1$ is strictly
monotonic, showing that $|\Hoare_0 X|=||X||+1$ is less than or equal
to $||Y||+1$.  This shows (2).

If additionally $f$ is continuous, then $\Sober f$ is well-defined and
monotonic by Lemma~\ref{lemma:surj}~(1), and coincides with the
injective map $f_*$ on $\Sober X$, hence is strictly monotonic.  Then
$\rk_{\Sober Y} \circ f_* \colon \Sober X \to \sob Y$ is strictly
monotonic, showing (3).

Finally, if $f$ is also Skula dense, then $f^{-1}$ is an order
embedding by Lemma~\ref{lemma:surj}~(3).  The fact that $f$ is full
means that $f^{-1}$ is surjective.  Hence $f^{-1}$ defines an order
isomorphism between $\Hoare_0 Y$ and $\Hoare_0 X$, which proves (4).
\qed

\begin{remark}
  \label{rem:qhom}
  A full, continuous, Skula dense map $f \colon X \to Y$ is the same
  thing as a map $f$ such that $f^{-1}$ is an order-isomorphism of
  $\Hoare_0 Y$ onto $\Hoare_0 X$, or equivalently of $\Open Y$ onto
  $\Open X$.  Such maps are called \emph{quasihomeomorphisms} in
  \cite[Definition~V-5.8]{GHKLMS:contlatt}.
\end{remark}

In the sequel, $F$ and other subsets of $X$ are given the subspace
topology induced from $X$.  The inclusion map is a topological
embedding, by definition.

\begin{corollary}
  \label{corl:stature:subspace}
  For every subset $F$ of a Noetherian space $X$, $||F|| \leq ||X||$
  and $\sob F \leq \sob X$.
\end{corollary}

The following will allow us to trade ranks of closed subsets for
statures of closed \emph{subspaces}.
\begin{lemma}
  \label{lemma:rk:C}
  For every closed subset $F$ of a Noetherian space $X$, $||F|| =
  \rk_{\Hoare_0 X} (F)$.
\end{lemma}
\proof
By well-founded induction on $F$, observing that the closed subsets of
$F$, seen as a topological subspace of $X$, are exactly the closed
subsets of $X$ that are included in $F$.
\qed

It will often be the case that we are able to determine the stature of
larger and larger proper closed subsets $F$ of a space $X$.  The
following will allow us to determine the stature of $X$ as a
consequence.
\begin{proposition}
  \label{prop:limit}
  Let $X$ be a Noetherian space, and ${(F_i)}_{i \in I}$ be a family
  of closed subsets of $X$ that is cofinal in the sense that every
  proper closed subset of $X$ is included in some $F_i$.  Then
  $||X|| \leq \sup_{i \in I} (||F_i||+1)$, and equality holds if all
  the subsets $F_i$ are proper.
\end{proposition}
\proof For every proper closed subset $F$ of $X$, we have
$F \subseteq F_i$ for some $i \in I$ by cofinality, hence
$||F|| \leq ||F_i||$ by Corollary~\ref{corl:stature:subspace}.  Since
$||X|| = \rk_{\Hoare_0 X} (X)$ is the supremum of
$\rk_{\Hoare_0 X} (F)+1$ when $F$ ranges over the proper closed
subsets of $X$, and $\rk_{\Hoare_0 X} (F) = ||F||$ by
Lemma~\ref{lemma:rk:C}, $||X|| \leq \sup_{i \in I} (||F_i||+1)$.  If
every $F_i$ is proper, then $||F_i|| < ||X||$ for every $i \in I$, by
Corollary~\ref{corl:stature:subspace}, so $||F_i||+1 \leq ||X||$ for
every $i$.  \qed

\begin{remark}
  \label{rem:inductive:limit}
  Let $X$ be non-empty and Noetherian, and ${(F_i)}_{i \in I}$ be a
  cofinal family of proper closed subsets of $X$, as in
  Proposition~\ref{prop:limit}.
  One may form the colimit of the diagram formed by the subspaces
  $F_i$ and the corresponding inclusion maps.  The result is $X$ with
  the topology \emph{determined} by the subspace topologies on each
  $F_i$; writing $\tau$ for that topology, a subset $F$ of $X$ is
  $\tau$-closed if and only if $F \cap F_i$ is closed in $F_i$ for
  every $i \in I$.  The topology $\tau$ contains the topology of $X$,
  but is in general much finer, as it contains closed sets that are
  included in no single $F_i$.  In general, a colimit of Noetherian
  spaces will fail to be Noetherian.  A typical counterexample is the
  collection of subsets $\{0, \cdots, n\}$ of $X \eqdef \nat$, each
  with the discrete topology, which arises from this construction by
  giving $X$ the cofinite topology.  In that case, the topology $\tau$
  is the discrete topology on $\nat$.  The cofinite topology on $\nat$
  is Noetherian, the discrete topology is not.
\end{remark}

\section{Finite spaces, well-founded chains}
\label{sec:finite-spaces}

Let $X$ be a finite $T_0$ space, of cardinality $n$.  The topology of
$X$ is necessarily the Alexandroff topology of its specialization
ordering $\leq$, and $X$ is automatically sober.  Clearly,
$0 \leq \sob X \leq n$.  The sobrification rank of $X$ can be as low
as $0$, if $X$ is empty.  For non-empty spaces, the sobrification rank
of $X$ can be as low as $1$, if $\leq$ is the equality ordering
(namely, if $X$ is $T_1$), and as high as $n$, if $X$ is a chain.

The stature of $X$, however, must be $n$ in all cases.
\begin{lemma}
  \label{lemma:finite}
  The stature $||X||$ of any finite $T_0$ space $X$ of cardinality $n
  \in \nat$ is $n$.
\end{lemma}
\proof This can be obtained from the fact that the maximal order type
of any finite poset with cardinality $n$ is $n$.  A direct proof is
equally easy.  By well-founded induction on $F \in \Hoare_0 X$,
$\rk_{\Hoare_0 X} (F)$ is less than or equal to the cardinality of
$F$.  Therefore $||X|| = \rk_{\Hoare_0 X} (X) \leq n$.  In order to
prove the reverse inequality, we exhibit a chain of length $n+1$ in
$\Hoare_0 X$.  We build points $x_i$, $1\leq i\leq n$, and closed
subsets $F_i$ of $X$, $0 \leq i \leq n$, by induction on $i$, such
that $F_i = \{x_1, \cdots, x_i\}$ for every $i$ and
$F_0 \subsetneq F_1 \subsetneq \cdots \subsetneq F_n$.  We do this by
letting $F_0 \eqdef \emptyset$, and, at stage $i$, by picking a
minimal point $x_i$ in $X \diff F_{i-1}$, so that
$F_i \eqdef F_{i-1} \cup \dc x_i = F_{i-1} \cup \{x_i\}$.  \qed

One can see that $\sob \nat = |\nat_\omega| = \omega+1$, where $\nat$
is given the Alexandroff topology of its usual ordering.  Here
$\nat_\omega = \Sober {\nat}$ is $\nat$ plus a fresh element $\omega$
on top of all others, with the upper topology
\cite[Theorem~5.4]{FGL:partI}.  Every non-empty closed subset of
$\nat$ is irreducible, so
$||\nat|| = |\Hoare_0 \nat|-1 = (1+\sob \nat)-1 = \omega$.

We generalize this to ordinals $\alpha$.  The upper topology on
$\alpha$ coincides with its Scott topology, and its open subsets are
the sets $\upc \beta$, where $\beta$ ranges over the non-limit
ordinals strictly smaller than $\alpha$, plus the empty set.  It is
indeed easy to see that such sets are Scott-open, while $\upc \beta$
is not when $\beta$ is a limit ordinal.

We recall that the notation $\alpha-1$ makes sense even when $\alpha$
is a limit ordinal, in which case it denotes $\alpha$ itself.
\begin{lemma}
  \label{lemma:ordinal}
  Let $\alpha$ be any ordinal.
  \begin{enumerate}
  \item With its Alexandroff topology, $\alpha$ has stature $\alpha$;
    its sobrification rank is $\alpha$ if $\alpha$ is finite, $\alpha+1$
    otherwise.
  \item With its upper (=Scott) topology, $\alpha$ has stature
    $\alpha$ if $\alpha$ is finite, $\alpha-1$ otherwise; its
    sobrification rank is $\alpha+1$ if $\alpha$ is a limit ordinal, $\alpha$
    otherwise.
  \end{enumerate}
\end{lemma}
\proof (1) Let us agree to write $\alpha$ for both the poset $\alpha$
and the space $\alpha$ with its Alexandroff topology.  The stature of
$\alpha$ is its maximal order type, since $\alpha$ is a wpo.  That
maximal order type is just the order type of $\alpha$, namely $\alpha$
itself, see \cite[Lemma~46]{BG:stature}, for instance.

The (downwards-)closed subsets of $\alpha$ are totally ordered by
inclusion.  Hence if $C$, $C_1$, $C_2$ are three closed subsets of
$\alpha$ such that $C \subseteq C_1 \cup C_2$, then $C$ is included in
the larger of $C_1$ and $C_2$; hence every non-empty closed subset $C$
of $\alpha$ is irreducible.  It follows that
$\Hoare_0 \alpha = \Sober \alpha \cup \{\emptyset\}$, so
$|\Sober \alpha \cup \{\emptyset\}| = 1+\sob \alpha$ is equal to
$|\Hoare_0 \alpha| = ||\alpha||+1 = \alpha+1$.  If $\alpha$ is finite,
then so is $1+\sob \alpha$, and $1+\sob\alpha = \sob\alpha+1$ is equal
to $\alpha+1$, so $\sob \alpha = \alpha$.  If $\alpha$ is infinite,
then so is $1+\sob \alpha$, and then $1+\sob\alpha = \sob\alpha$;
therefore $\sob\alpha = \alpha+1$.

(2) We now agree to write $\alpha$ for the space $\alpha$ with its
upper topology.  Its closed subsets are exactly the empty set, the
sets $\dc \beta$, with $\beta < \alpha$, plus the whole set $\alpha$
if $\alpha$ is a limit ordinal.  Hence, if $\alpha$ is a limit
ordinal, then $|\Hoare_0 X| = 1+\alpha+1$, so $||X|| = 1+\alpha$;
since $\alpha$ is infinite, $1+\alpha=\alpha$, which is equal to
$\alpha-1$ by convention, since $\alpha$ is infinite.  If $\alpha$ is not a
limit ordinal, then $|\Hoare_0 X| = 1+\alpha$ instead; then, if
$\alpha$ is finite, then $||X|| = \alpha$, otherwise $\alpha$ is a
successor ordinal and $||X|| = \alpha-1$.

All the non-empty closed subsets are irreducible, so a similar
analysis applies to $\Sober \alpha$.  If $\alpha$ is a limit ordinal,
then $|\Sober X| = \alpha+1$.  Otherwise, $|\Sober X| = \alpha$.  \qed


\section{Spaces with a cofinite topology}
\label{sec:spac-witha-cofin}

Any space obtained by equipping a set $X$ with the cofinite topology is
Noetherian, since $\Hoare_0 X$ is obviously well-founded in that case.
The specialization ordering of such a space is equality, which is not
a wqo unless $X$ is finite.  Note that every such space is $T_1$.

\begin{theorem}
  \label{thm:dim:cofinite}
  The sobrification rank of a space $X$ with the cofinite topology is $0$ if
  $X$ is empty, $1$ if $X$ is finite and non-empty, and $2$ if $X$ is
  infinite.
\end{theorem}
\proof The case of the empty space is obvious.  If $X$ is finite, then
the cofinite topology is the discrete topology.  In that case, $X$ is
sober, so $\Sober X$ and $X$ are homeomorphic, and no point is
strictly above any other, so the rank of every point is $0$.  It
follows that, if $X$ is non-empty, then $\sob X = 1$.

If $X$ is infinite, then the irreducible closed subsets of $X$ are the
sets of the form $\{x\}$ with $x \in X$, plus the whole space itself.
In order to see this, it suffices to show that $X$ is irreducible, and
that no proper closed subset with at least two points is.
As far as the first point is concerned, if $X$ is included in the
union of two closed subsets $F_1$ and $F_2$, then at least one of them
is infinite (since $X$ is) and therefore equal to $X$, since the only
infinite closed subset of $X$ is $X$ itself.
As far as the second point is concerned, let $F$ be a proper closed
subset of $X$.  Hence $F$ is finite.  If $F$ contained at least two
points $x$ and $y$, then it would be included in the union of the two
closed sets $F \diff \{x\}$ and $\{x\}$ without being included in
either.

Now the rank of each irreducible closed subset of the form $\{x\}$ is
$0$.  Then the rank of $X$ in $\Sober X$ is $1$, and therefore
$\sob X = 2$.  \qed

\begin{theorem}
  \label{thm:stature:cofinite}
  The stature $||X||$ of a space $X$ with the cofinite topology is
  the cardinality of $X$ if $X$ is finite, and $\omega$ otherwise.
\end{theorem}
\proof When $X$ is finite, this is Lemma~\ref{lemma:finite}.  Let us
assume $X$ infinite.  For every proper closed subset $F$ of $X$,
$\rk_{\Hoare_0 X} (F) = ||F||$ by Lemma~\ref{lemma:rk:C}, and this is
equal to the cardinality of $F$, since $F$ is finite, by
Lemma~\ref{lemma:finite}.  Those values span the whole of $\nat$ as
$F$ varies, since $X$ is infinite.  It follows that
$||X|| = \rk_{\Hoare_0 X} (X)$, which is the least ordinal strictly
larger than those, is equal to $\omega$.  \qed

As a corollary, spaces with a cofinite topology yield examples of
spaces with very low sobrification rank and stature, and with arbitrarily high
cardinality.

\section{Sums}
\label{sec:sums}

Let $P+Q$ denote the coproduct of the two posets $P$ and $Q$ in the
category of posets and monotonic maps.  This is the disjoint union
of $P$ and $Q$, where all elements of $P$ are incomparable with all
elements of $Q$, and the ordering relations inside $P$ and inside $Q$
are preserved.

The maximal order type $o (P+Q)$ of the sum of two wpos $P$ and $Q$ is
equal to $o (P) \oplus o (Q)$.  We have a similar result for statures,
which we will prove using the following lemma, which one can find in
\cite[Section~4.8.3]{Fraisse:relations}, for example.
\begin{lemma}
  \label{lemma:rk:x}
  Let $P$ and $Q$ be two well-founded posets.  For all $p \in P$ and
  $q \in Q$, $\rk_{P \times Q} (p, q) = \rk_P (p) \oplus \rk_Q (q)$.
\end{lemma}


Let $X+Y$ denote the topological sum of two topological spaces $X$
and $Y$, namely their coproduct in the category $\Top$.  This is
Noetherian as soon as both $X$ and $Y$ are.
\begin{proposition}
  \label{prop:stature:sum}
  For all Noetherian spaces $X$ and $Y$,
  $||X+Y|| = ||X|| \oplus ||Y||$.
\end{proposition}
\proof Every closed subset $F$ of $X+Y$ can be written in a unique way
as the disjoint union of $F \cap X$ and of $F \cap Y$, which are
closed subsets of $X$ and of $Y$, respectively.  It follows that
$\Hoare_0 (X+Y)$ and $\Hoare_0 X \times \Hoare_0 Y$ are
order-isomorphic.  Then
$||X+Y|| = \rk_{\Hoare_0 (X+Y)} (X+Y) = \rk_{\Hoare_0 X \times
  \Hoare_0 Y} (X, Y) = \rk_{\Hoare_0 X} (X) \oplus \rk_{\Hoare_0 Y}
(Y) = ||X|| \oplus ||Y||$, where the next-to-last equality is by
Lemma~\ref{lemma:rk:x}.  \qed

\begin{lemma}
  \label{lemma:rk:+}
  Given any two well-founded posets $P$ and $Q$,
  $|P+Q| = \max (|P|, |Q|)$.
\end{lemma}
\proof $|P+Q|$ is the least ordinal strictly larger than
$\rk_{P+Q} (p) = \rk_P (p)$ for every $p \in P$ and than
$\rk_{P+Q} (q) = \rk_Q (q)$ for every $q \in Q$.  \qed

\begin{proposition}
  \label{prop:dim:sum}
  For all Noetherian spaces $X$ and $Y$,
  $\sob (X+Y) = \max (\sob X, \allowbreak \sob Y)$
  and $\rsob (X+Y) = \max (\rsob X, \allowbreak \rsob Y)$.
\end{proposition}
\proof By Lemma~\ref{lemma:rk:+}, since the sobrification of $X+Y$ is
the disjoint sum of the posets $\Sober X$ and $\Sober Y$.  The latter
can be seen by realizing that any irreducible closed subset of $X+Y$
must be included in $X$ or in $Y$.  \qed

\section{Lexicographic sums}
\label{sec:bf-lexic-sums}

There is another notion of sum, which we call the \emph{lexicographic
  sum} of two spaces.
\begin{definition}[Lexicographic sum]
  \label{defn:+lex}
  The lexicographic sum $X +_\lex Y$ of two topological spaces $X$ and
  $Y$ is the disjoint sum of $X$ and $Y$, and its open subsets are the
  open subsets of $Y$, plus the sets of the form $U+Y$, where $U$ is
  an open subset of $X$.
\end{definition}
If $X$ and $Y$ are Noetherian, then so is $X +_\lex Y$, since $X+Y$ is
and every topology coarser than a Noetherian topology is Noetherian.
The specialization preordering of $X +_\lex Y$ is the lexicographic
sum of those of $X$ and of $Y$, where all elements of $X$ are 
below all elements of $Y$, and the topology of $X +_\lex Y$ is
Alexandroff if those of $X$ and $Y$ are.
\begin{proposition}
  \label{prop:stature:+lex}
  For all Noetherian spaces $X$ and $Y$,
  $||X +_\lex Y|| = ||X|| + ||Y||$.
\end{proposition}
\proof The closed subsets of $X +_\lex Y$ are the proper closed
subsets $F$ of $X$, plus the sets of the form $X+F'$, where $F'$ ranges
over the closed subsets of $Y$.  By well-founded induction, we see
that $\rk_{X +_\lex Y} (F) = \rk_X (F)$ for subsets of the first kind,
then that $\rk_{X +_\lex Y} (X+F') = ||X|| + \rk_Y (F')$ for subsets of
the second kind.  It follows that
$||X +_\lex Y|| = \rk_{X +_\lex Y} (X+Y) = ||X||+||Y||$.
\qed

\begin{lemma}
  \label{lemma:S:+lex}
  For any two topological spaces $X$ and $Y$, the irreducible closed
  subsets of $X +_\lex Y$ are those of $X$, plus sets of the form
  $X+D$, where $D$ is irreducible closed in $Y$.
\end{lemma}
\proof If $C$ is a closed subset of $X$ that is irreducible in $X$,
then we claim that it is irreducible in $X +_\lex Y$.  By assumption,
$C$ is non-empty.  Let us assume that $C$ is included in the union of
two closed subsets $F_1$ and $F_2$ of $X +_\lex Y$.  If one of them is
of the form $X+F'$ for some closed subset $F'$ of $Y$, then $C$ is
included in that one.  Otherwise, both are closed subsets of $X$, and
therefore $C$ is included in $F_1$ or in $F_2$.

Conversely, if $C$ is irreducible closed in $X +_\lex Y$, and included
in $X$, then it is easy to see that $C$ is irreducible closed in $X$.

Let us consider any set of the form $X+D$, where $D$ is closed in $Y$.
If $D$ is irreducible closed in $Y$, then $D$ is non-empty, hence so
is $X+D$.  Let us assume that $X+D$ is included in the union of two
closed subsets $F_1$ and $F_2$ of $X +_\lex Y$.  It cannot be that
$F_1$ and $F_2$ are both included in $X$, since $X+D$ is not, owing to
the fact that $D$ is non-empty.  If $F_1$ is of the form $X+F'_1$ for
some closed subset of $Y$, and $F_2$ is included in $X$, then $F_1
\cup F_2 = F_1$, so $X+D$ is included in $F_1$.  We omit the
symmetrical case.  If $F_1 = X+F'_1$ and $F_2 = X+F'_2$ for some closed
subsets $F'_1$ and $F'_2$ of $Y$, finally, then $X+D \subseteq (X+F'_1)
\cup (X+F'_2)$ implies that $D$ is included in $F'_1 \cup F'_2$, hence in
one of $F'_1$ or $F'_2$, since $D$ is irreducible.  Therefore $X+D$ is
included in $F_1$ or in $F_2$.

Conversely, let us assume that $X+D$ is irreducible closed in
$X +_\lex Y$.  If $D$ is empty, then $X$ is irreducible closed in
$X +_\lex Y$, and we have seen that $X$ must be irreducible closed in
$X$.  Henceforth, we assume that $D$ is non-empty.  Then, given any
two closed subsets $F'_1$ and $F'_2$ of $Y$ whose union contains $D$,
the union of $X+F'_1$ and of $X+F'_2$ contains $X+D$.  Since $X+D$ is
irreducible, $X+D$ is included in one of them, and therefore $D$ is
included in $F'_1$ or in $F'_2$.  \qed

\begin{proposition}
  \label{prop:dim:+lex}
  For all Noetherian spaces $X$ and $Y$,
  $\sob (X +_\lex Y) = \sob X + \sob Y$.
\end{proposition}
\proof
We leave all references to Lemma~\ref{lemma:S:+lex} implicit here.
By well-founded induction on $C \in \Sober X$,
$\rk_{\Sober (X +_\lex Y)} (C) = \rk_X (C)$.
Given any minimal element $D$ of $\Sober Y$,
hence such that $\rk_{\Sober Y} (D)=0$,
$\rk_{\Sober (X +_\lex Y)}
(X+D)$ is the smallest ordinal strictly larger than
$\rk_{\Sober X} (C)$ for every $C \in \Sober X$, hence is equal to
$\sob X$, by definition.  This is the start of an induction on
$D \in \Sober Y$, showing that
$\rk_{\Sober (X +_\lex Y)} (X+D) = \sob X + \rk_Y (D)$.  By adding one
and taking suprema, we obtain that
$\sob (X +_\lex Y) = \sob X + \sob Y$.
\qed

The special case of \emph{liftings} is of particular importance.
\begin{definition}[Lifting]
  \label{defn:lifting}
  The \emph{lifting} $X_\bot$ of a topological space is the
  lexicographic sum of a one-point space with $X$, in short,
  $\{\bot\} +_\lex X$.
\end{definition}
The open subsets of $X_\bot$ are those of $X$, plus $X_\bot$ itself.
$X_\bot$ is Noetherian if and only if $X$ is.  The following is an
easy consequence of Proposition~\ref{prop:dim:+lex} and of
Proposition~\ref{prop:stature:+lex}.
\begin{proposition}
  \label{prop:lifting}
  For every Noetherian space $X$, $\sob {X_\bot} = 1+\sob X$ and
  $||X_\bot|| = 1+||X||$.
\end{proposition}

\section{Products}
\label{sec:products}


The sobrification rank of a product is easily obtained.  It suffices to observe
that the irreducible closed subsets of $X \times Y$ are exactly the
products $C \times D$ of an irreducible closed subset $C$ of $X$ and
of an irreducible closed subset $D$ of $Y$, and therefore
$\Sober (X \times Y)$ is order-isomorphic to the poset product $\Sober
X \times \Sober Y$; this is originally due to
Hoffmann \cite[Theorem~1.4]{Hoffmann:sobr:rem}, see also
\cite[Proposition~8.4.7]{JGL:topo}.

\begin{proposition}
  \label{prop:dim:prod}
  For all non-empty Noetherian spaces $X$ and $Y$,
  $\sob (X \times Y) = (\sob X \oplus \sob Y)-1$.  If one of them is
  empty, then $\sob (X \times Y) = 0$.
\end{proposition}
\proof
If $X$ is empty, then so is $X \times Y$, and therefore
$\sob (X \times Y) = 0$; similarly if $Y$ is empty.

Let us assume that $X$ and $Y$ are both non-empty.  We write $X$ as a
finite union of irreducible closed subsets $C_1$, \ldots, $C_m$, and
$Y$ as a finite union of irreducible closed subsets $D_1$, \ldots,
$D_n$.  Using Lemma~\ref{lemma:rk:x},
$\rk_{\Sober (X \times Y)} (C_i \times D_j) = \rk_{\Sober X \times
  \Sober Y} (C_i, D_j) = \rk_{\Sober X} (C_i) \oplus \rk_{\Sober Y}
(D_j)$ for all $i$ and $j$.  Then
$\sob (X \times Y) = \max_{i, j} \rk_{\Sober (X \times Y)} (C_i \times
D_j)+1 = \max_{i, j} \rk_{\Sober X} (C_i) \oplus \rk_{\Sober Y}
(D_j)+1$, while
$\sob X \oplus \sob Y = (\max_i \rk_{\Sober X} (C_i) +1) \oplus
(\max_j \rk_{\Sober Y} (D_j) +1) = (\max_i \rk_{\Sober X} (C_i)
\oplus1) \oplus (\max_j \rk_{\Sober Y} (D_j) \oplus 1) =\max_{i, j}
\rk_{\Sober X} (C_i) \oplus \rk_{\Sober Y} (D_j) \oplus 2$, which is
therefore equal to $\sob (X \times Y) + 1$.  \qed


We turn to statures of products.  The corresponding result on wpos is
that the maximal order type $o (P \times Q)$ of the product of two
wpos $P$ and $Q$ is equal to $o (P) \otimes o (Q)$, as shown by de
Jongh and Parikh \cite[Section~3]{dJP:maxot}.  We generalize this to
Noetherian spaces, replacing maximal order types by statures.  The
general outline of the argument resembles de Jongh and Parikh's, but
the details vary considerably: de Jongh and Parikh extensively build
upwards-closed subsets as upward closures of finite sets of points,
and that is a technique that is not available to us in general
Noetherian spaces.


\begin{lemma}
  \label{lemma:stature:C}
  For every Noetherian space $X$ and every subset $F$ of $X$,
  \begin{enumerate}
  \item $||X|| \leq ||F|| \oplus ||X \diff F||$;
  \item if $F$ is closed, then $||F|| + ||X \diff F|| \leq ||X||$.
  \end{enumerate}
\end{lemma}
\proof (1) The identity map is continuous from $F + (X \diff F)$ to
$X$, because every open subset $U$ of $X$ can be written as
$(U \cap F) + (U \diff F)$.  The inequality then follows from
Lemma~\ref{lemma:stature:mono}, together with
Proposition~\ref{prop:stature:sum}.

(2) We assume that $F$ is closed.  We claim that the identity map is
continuous from $X$ to $F +_\lex (X \diff F)$.  In order to verify
this, we consider any open set of the latter space, and we verify that
it is open in $X$.  There are two kinds of open subsets of
$F +_\lex (X \diff F)$.  The open subsets $U$ of $X \diff F$ are
intersections $V \cap (X \diff F)$ of an open subset $V$ of $X$ with
$X \diff F$; this intersection must also be open in $X$ since $F$ is
closed.  The open subsets of the form $U+(X \diff F)$ where $U$ is
open in $F$, namely where $U = V \cap F$ for some open subset $V$ of
$X$, are equal to $V \cup (X \diff F)$, hence are open in $X$, once
again because $F$ is closed.

Now (2) follows from Lemma~\ref{lemma:stature:mono}, together with
Proposition~\ref{prop:stature:+lex}.
\qed

The following is the key result on which the main theorem of this
section is built.  The corresponding result in \cite{dJP:maxot}, in
the special case of wpos, is Corollary~2.17.
\begin{corollary}
  \label{corl:decomp:C}
  For every Noetherian space $X$, whose stature is a decomposable
  ordinal $\alpha$, written in Cantor normal form as
  $\omega^{\alpha_1} + \cdots + \omega^{\alpha_m}$ with
  $\alpha \geq \alpha_1 \geq \cdots \geq \alpha_m$, $m \geq 2$,
  \begin{enumerate}
  \item there is a closed subset $F$ of $X$ such that
    $||F|| = \omega^{\alpha_1} + \cdots + \omega^{\alpha_{m-1}}$;
  \item for each such closed subset $F$,
    $||X \diff F|| = \omega^{\alpha_m}$.
  \end{enumerate}
\end{corollary}
\proof (1) Let $\beta \eqdef \omega^{\alpha_1} + \cdots + \alpha^{\alpha_{m-1}}$.  We
have $\beta < \alpha \leq ||X||+1 = |\Hoare_0 X|$, so there is an
element $F$ of $\Hoare_0 X$ such that $\rk_{\Hoare_0 X} (F) = \beta$.
Therefore $\beta = ||F||$, by Lemma~\ref{lemma:rk:C}.

(2)  By Lemma~\ref{lemma:stature:C},
$||F||+||X \diff F|| \leq ||X|| \leq ||F|| \oplus ||X \diff F||$,
namely
$\beta+||X \diff F|| \leq \beta + \omega^{\alpha_m} \leq \beta \oplus
||X \diff F||$.  The first inequality implies
$||X \diff F|| \leq \omega^{\alpha_m}$.  Let us write $||X \diff F||$
in Cantor normal form as
$\omega^{\beta_1} + \cdots + \omega^{\beta_n}$ with
$\beta_1 \geq \cdots \geq \beta_n$, $n \in \nat$.  The second
inequality tells us that
$\beta + \omega^{\alpha_m} = \omega^{\alpha_1} + \cdots +
\omega^{\alpha_{m-1}} + \omega^{\alpha_m}$ is less than or equal to
the sum of the terms $\omega^{\alpha_1}$, \ldots,
$\omega^{\alpha_{m-1}}$ and $\omega^{\beta_1}$, \ldots,
$\omega^{\beta_n}$, sorted in decreasing order.  Hence the list
$\alpha_1, \cdots, \alpha_m$ is lexicographically smaller than or
equal to the list obtained by merging the two lists
$\alpha_1, \cdots, \alpha_{m-1}$ and $\beta_1, \cdots, \beta_n$ and
sorting the result in descending order.  It follows that $n \geq 1$
and $\alpha_i \leq \beta_1$ for some $i$, $1\leq i\leq m$.  This
entails $\alpha_m \leq \beta_1$, hence
$\omega^{\alpha_m} \leq ||X \diff F||$.  Therefore
$||X \diff F|| = \omega^{\alpha_m}$.  \qed

\begin{lemma}
  \label{lemma:rk:union}
  For every finite list of closed subsets $F_1$, \ldots, $F_n$ of a
  Noetherian space $X$,
  $||\bigcup_{i=1}^n F_i|| \leq \bigoplus_{i=1}^n ||F_i||$.
\end{lemma}
\proof The map from $F_1 + \cdots + F_n$ that sends each element of
$F_i$ to itself in $\bigcup_{i=1}^n F_i$ is surjective, and
continuous.  By Lemma~\ref{lemma:stature:mono},
$||\bigcup_{i=1}^n F_i|| \leq ||F_1 + \cdots + F_n|| = ||F_1|| \oplus
\cdots \oplus ||F_n||$; the second equality follows from
Proposition~\ref{prop:stature:sum}.  \qed

\begin{proposition}
  \label{prop:stature:prod:ubound}
  For all Noetherian spaces $X$ and $Y$,
  $||X \times Y|| \leq ||X|| \otimes ||Y||$.
\end{proposition}
\proof By induction on the pair of ordinals $||X||$ and $||Y||$,
ordered lexicographically.  The claim is clear if $||X||=0$ or if
$||Y||=0$.  Otherwise, let us write $||X||$ in Cantor normal form as
$\omega^{\alpha_1} + \cdots + \omega^{\alpha_m}$, with
$\alpha_1 \geq \cdots \geq \alpha_m$, $m \geq 1$, and $||Y||$ in
Cantor normal form as $\omega^{\beta_1} + \cdots + \omega^{\beta_n}$, with
$\beta_1 \geq \cdots \geq \beta_n$, $n \geq 1$.

If $m \geq 2$, then by Corollary~\ref{corl:decomp:C}, there is a
closed subset $F$ of $X$ such that
$||F|| = \omega^{\alpha_1} + \cdots + \alpha^{\alpha_{m-1}}$, and
$||X\diff F|| = \omega^{\alpha_m}$.  Then $F \times Y$ is a closed
subset of $X \times Y$, and therefore
$||X \times Y|| \leq ||F \times Y|| \oplus ||(X \diff F) \times Y||$
by Lemma~\ref{lemma:stature:C}~(1).  By induction hypothesis,
$||X \times Y|| \leq (||F|| \otimes ||Y||) \oplus (\omega^{\alpha_m}
\otimes ||Y||) = (||F|| \oplus \omega^{\alpha_m}) \otimes ||Y|| =
||X|| \otimes ||Y||$.

The case where $n \geq 2$ is symmetric.

Finally, we examine the case where $m=n=1$.  Then
$||X|| = \omega^{\alpha_1}$ and $||Y|| = \omega^{\beta_1}$.  We wish
to show that
$||X \times Y|| = \rk_{\Hoare_0 (X \times Y)} (X \times Y)$ is smaller
than or equal to
$||X|| \otimes ||Y|| = \omega^{\alpha_1 \oplus \beta_1}$.  To this
end, it suffices to show that
$\rk_{\Hoare_0 (X \times Y)} (F) < \omega^{\alpha_1 \oplus \beta_1}$
for every proper closed subset $F$ of $X \times Y$.  We write $F$ as a
finite union of irreducible closed subsets of $X \times Y$.  Each one
must be of the form $C_k \times D_k$, where $C_k$ is irreducible
closed in $X$ and $D_k$ is irreducible closed in $Y$, where $k$ ranges
from $1$ to $p$, say.  (The fact that they are irreducible will not
matter.  The important thing is that we can write $F$ as a finite
union of products of non-empty closed subsets.)  Additionally, $C_k$
is a proper subset of $X$ or $D_k$ is a proper subset of $Y$.  Hence
we can apply the induction hypothesis, to the effect that
$||C_k \times D_k|| \leq ||C_k|| \otimes ||D_k||$.  Moreover,
$||C_k|| \leq ||X|| = \omega^{\alpha_1}$ and
$||D_k|| \leq ||Y|| = \omega^{\beta_1}$, where at least one equality
is strict.  Whatever the case is, we obtain that
$||C_k \times D_k|| < \omega^{\alpha_1} \otimes \omega^{\beta_1} =
\omega^{\alpha_1 \oplus \beta_1}$, since natural product is monotonic,
and strictly monotonic in each of its arguments, provided the other
one is not zero.
Using Lemma~\ref{lemma:rk:union}, we have
$||F|| \leq \bigoplus_{k=1}^p ||C_k \times D_k||$.  Since
$\omega^{\alpha_1 \oplus \beta_1}$ is $\oplus$-indecomposable,
$||F|| < \omega^{\alpha_1 \oplus \beta_1}$, and therefore,
$\rk_{\Hoare_0 (X \times Y)} (F) < \omega^{\alpha_1 \oplus \beta_1}$,
using Lemma~\ref{lemma:rk:C}.  \qed


We turn to the reverse inequality.

\begin{lemma}
  \label{lemma:decomp:C:*}
  For every Noetherian space $X$, whose stature is written in Cantor
  normal form as $\omega^{\alpha_1} + \cdots + \omega^{\alpha_m}$, with
  $\alpha_1 \geq \cdots \geq \alpha_m$, $m \in \nat$, there are closed
  subsets
  $X=F_m \supseteq \cdots \supseteq F_1 \supseteq F_0 = \emptyset$
  such that $||F_i|| = \omega^{\alpha_1} + \cdots + \omega^{\alpha_i}$
  for every $i$, $0\leq i \leq m$, and
  $||F_i \diff F_{i-1}|| = \omega^{\alpha_i}$ for every $i$,
  $1\leq i\leq m$.
\end{lemma}
\proof
By induction on $m$, using Corollary~\ref{corl:decomp:C}.
\qed

\begin{lemma}
  \label{lemma:stature:C:*}
  For every Noetherian space $X$, and every sequence
  $X=F_m \supseteq \cdots \supseteq F_1 \supseteq F_0=\emptyset$ of
  closed subsets of $X$,
  \[
    ||F_1 \diff F_0|| + ||F_2 \diff F_1|| + \cdots + ||F_m\diff
    F_{m-1}||
    \leq ||X|| \leq \bigoplus_{i=1}^m ||F_i \diff F_{i-1}||.
  \]
\end{lemma}
\proof If $m=0$, then $X$ is empty, and this is clear.  Otherwise, we
use Lemma~\ref{lemma:stature:C}~(1), on $X=F_m$, $F_{m-1}$,
\ldots, $F_1$ in succession in order to obtain
$||X|| \leq ||F_{m-1}|| \oplus ||F_m \diff F_{m-1}||$,
$||F_{m-1}|| \leq ||F_{m-2}|| \oplus ||F_{m-1} \diff F_{m-2}||$,
\ldots,
$||F_2|| \leq ||F_1|| \oplus ||F_2 \diff F_1|| = ||F_1 \diff F_0||
\oplus ||F_2 \diff F_1||$, from which the second inequality follows.
The first inequality is proved similarly, using
Lemma~\ref{lemma:stature:C}~(2), instead.  \qed

\begin{proposition}
  \label{prop:stature:prod:lbound}
  For all Noetherian spaces $X$ and $Y$,
  $||X \times Y|| \geq ||X|| \otimes ||Y||$.
\end{proposition}
\proof
By induction on the pair of ordinals $||X||$ and $||Y||$,
ordered lexicographically.  The claim is clear if $||X||=0$ or if
$||Y||=0$.  Otherwise, let us write $||X||$ in Cantor normal form as
$\omega^{\alpha_1} + \cdots + \omega^{\alpha_m}$, with
$\alpha_1 \geq \cdots \geq \alpha_m$, $m \geq 1$, and $||Y||$ in
Cantor normal form as $\omega^{\beta_1} + \cdots + \omega^{\beta_n}$, with
$\beta_1 \geq \cdots \geq \beta_n$, $n \geq 1$.

\begin{figure}
  \centering
  \includegraphics[scale=0.4]{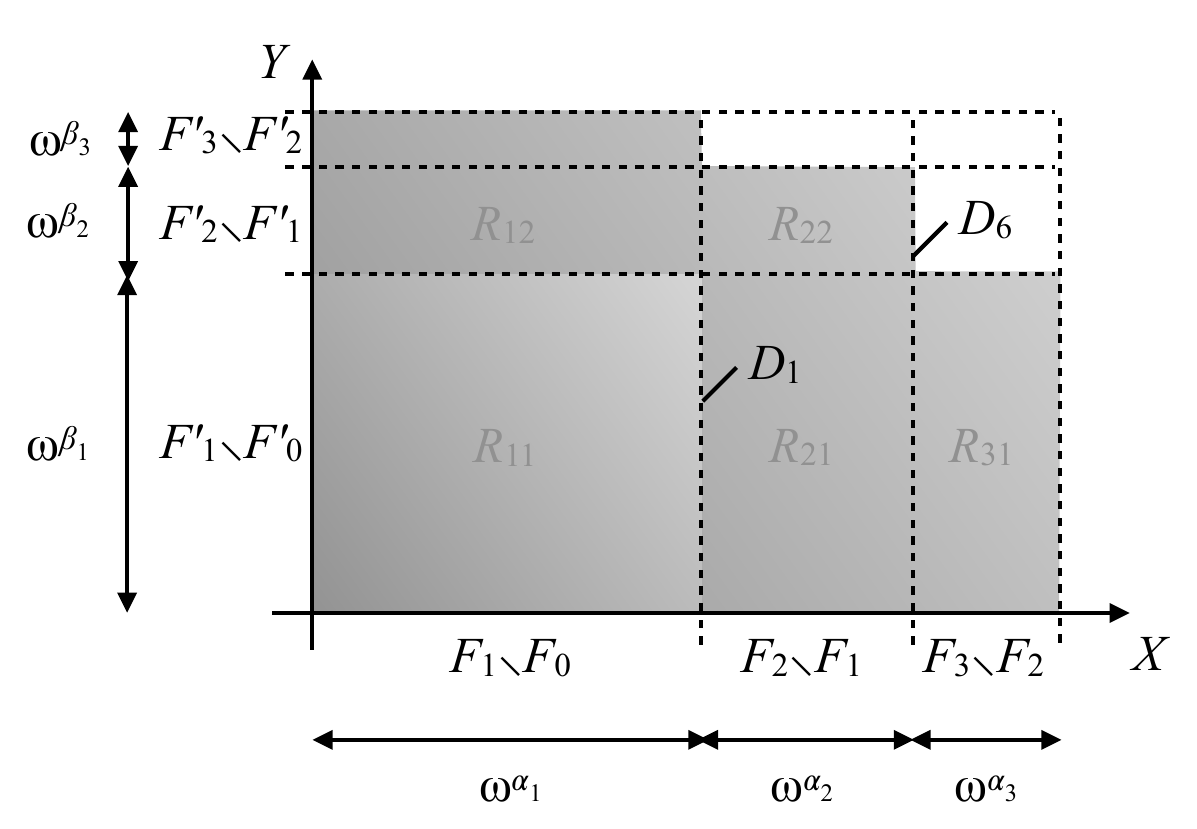}
  \caption{Proving $||X \times Y|| \geq ||X|| \otimes ||Y||$}
  \label{fig:stature:product}
\end{figure}

\emph{The case $m \geq 2$ or $n \geq 2$.}  By
Lemma~\ref{lemma:decomp:C:*}, there are closed subsets
$X=F_m \supseteq \cdots \supseteq F_1 \supseteq F_0 = \emptyset$ such
that $||F_i|| = \omega^{\alpha_1} + \cdots + \omega^{\alpha_i}$ for
every $i$, $0\leq i \leq m$, and
$||F_i \diff F_{i-1}|| = \omega^{\alpha_i}$ for every $i$,
$1\leq i\leq m$.  Similarly, there are closed subsets
$Y=F'_n \supseteq \cdots \supseteq F'_1 \supseteq F'_0 = \emptyset$ such
that $||F'_j|| = \omega^{\beta_1} + \cdots + \omega^{\beta_j}$ for
every $j$, $0\leq j \leq n$, and
$||F'_j \diff F'_{j-1}|| = \omega^{\beta_j}$ for every $j$,
$1\leq j\leq n$.  The situation is illustrated in
Figure~\ref{fig:stature:product}, where $m=n=3$.  The sets $D_1$ and
$D_6$ are instances of a sequence of closed subsets $D_k$ that we will
construct below.

There is a strict ordering $\prec_0$ on pairs $(i,j)$ defined by:
$(i, j) \prec_0 (i', j')$ if and only if
$\alpha_i \oplus \beta_j > \alpha_{i'} \oplus \beta_{j'}$, or
$\alpha_i \oplus \beta_j = \alpha_{i'} \oplus \beta_{j'}$ and
$F_i \times F'_j \subsetneq F_{i'} \times F'_{j'}$.  Let $\prec$ be any
linear extension of $\prec_0$, and let us enumerate the pairs $(i,j)$
as $(i_1, j_1) \prec \cdots \prec (i_{mn}, j_{mn})$.  (In other
words, let us sort the pairs $(i, j)$ with respect to $\prec_0$.)

For every pair $(i, j) \in \{1, \cdots, m\} \times \{1, \cdots, n\}$,
we abbreviate the rectangle
$(F_i \diff F_{i-1}) \times (F'_j \diff F'_{j-1})$ as $R_{ij}$.  For
every $k \in \{1, \cdots, mn\}$, we claim that
$D_k \eqdef \bigcup_{\ell=1}^k R_{i_\ell j_\ell}$ is also equal to
$\bigcup_{\ell=1}^k (F_{i_\ell} \times F'_{j_\ell})$, and is therefore
closed.  The inclusion
$D_k \subseteq \bigcup_{\ell=1}^k (F_{i_\ell} \times F'_{j_\ell})$ is
clear.  In order to show the reverse inclusion, it suffices to show
that $F_{i_\ell} \times F'_{j_\ell}$ is included in $D_k$ for every
$\ell \in \{1, \cdots, k\}$.  One checks easily that every point
$(x, y)$ of $F_{i_\ell} \times F'_{j_\ell}$ lies in some rectangle
$R_{ij}$ with $i \leq i_\ell$ and $j \leq j_\ell$.  If $i=i_\ell$ and
$j=j_\ell$, then $(x,y)$ is in $R_{i_\ell j_\ell}$, hence in $D_k$.
Otherwise, $\alpha_i \geq \alpha_{i_\ell}$,
$\beta_j \geq \beta_{j_\ell}$, $F_i \subseteq F_{i_\ell}$, and
$F'_j \subseteq F'_{j_\ell}$, where one of the two inclusions is strict.
It follows that $(i, j) \prec_0 (i_\ell, j_\ell)$.  Hence $(i,j)$
occurs as $(i_{\ell'}, j_{\ell'})$ for some $\ell'$ such that
$1\leq \ell' < \ell$.  It follows that $(x, y)$ is in
$R_{i_{\ell'} j_{\ell'}}$, hence in $D_k$.

Let $D_0 \eqdef \emptyset$.  We now claim that
$D_k \diff D_{k-1} = R_{i_k j_k}$ for every $k$, $1\leq k\leq mn$.
Indeed, this follows from the fact that the union
$\bigcup_{\ell=1}^k R_{i_\ell j_\ell}$ defining $D_k$ is a disjoint
union, and can be reorganized as the union of
$\bigcup_{\ell=1}^{k-1} R_{i_\ell j_\ell} = D_{k-1}$ and of
$R_{i_k j_k}$.

We now have
$D_{mn} \supseteq \cdots \supseteq D_1 \supseteq D_0 = \emptyset$, and
$D_{mn} = \bigcup_{\ell=1}^{mn} F_{i_\ell} \times F'_{j_\ell} = X
\times Y$.  This allows us to use Lemma~\ref{lemma:stature:C:*}, so
$||X \times Y|| \geq ||D_1\diff D_0|| + ||D_2 \diff D_1|| + \cdots +
||D_{mn} \diff D_{mn-1}||$.  In other words,
$||X \times Y|| \geq ||R_{i_1 j_1}|| + ||R_{i_2 j_2}|| + \cdots +
||R_{i_{mn} j_{mn}}||$.  By induction hypothesis, for every pair
$(i, j)$,
$||R_{ij}|| = ||(F_i \diff F_{i-1}) \times (F'_j \diff F'_{j-1})||$ is
equal to $\omega^{\alpha_i} \otimes \omega^{\beta_j}$, namely to
$\omega^{\alpha_i \oplus \beta_j}$.  (Note that the induction
hypothesis applies because the pair
$(\omega^{\alpha_i}, \omega^{\beta_j})$ is lexicographically smaller
than the pair $(||X||, ||Y||)$, and this is so because $m \geq 2$ or
$n \geq 2$.)  We have obtained that $||X \times Y||$ is larger than or
equal to
$\omega^{\alpha_{i_1} \oplus \beta_{j_1}} + \cdots +
\omega^{\alpha_{i_{mn}} \oplus \beta_{j_{mn}}}$.

Let us observe that
$\alpha_{i_k} \oplus \beta_{j_k} \geq \alpha_{i_{k+1}} \oplus
\beta_{j_{k+1}}$ for every $k$ with $1\leq k < mn$.  Indeed, otherwise
we would have
$\alpha_{i_{k+1}} \oplus \beta_{j_{k+1}} > \alpha_{i_k} \oplus
\beta_{j_k}$, hence $(i_{k+1}, j_{k+1}) \prec_0 (i_k, j_k)$ by
definition of $\prec_0$, and therefore
$(i_{k+1}, j_{k+1}) \prec (i_k, j_k)$, which is impossible.
  
It follows that the list of ordinals
$\omega^{\alpha_{i_1} \oplus \beta_{j_1}}, \cdots,
\omega^{\alpha_{i_{mn}} \oplus \beta_{j_{mn}}}$ is sorted in
decreasing order.  That list enumerates all the ordinals
$\omega^{\alpha_i \oplus \beta_j}$ with $1\leq i\leq m$ and
$1\leq j\leq n$.  Their sum therefore equals
$(\omega^{\alpha_1} + \cdots + \omega^{\alpha_m}) \otimes
(\omega^{\beta_1} + \cdots + \omega^{\beta_n}) = ||X|| \otimes ||Y||$.

\emph{The case $m=n=1$.}  We now assume that
$||X|| = \omega^{\alpha_1}$ and $||Y|| = \omega^{\beta_1}$.  If
$\alpha_1=0$, then $||X||=1$, so there are exactly two closed subsets
of $X$, the empty set and $X$ itself.  The closed subsets of
$X \times Y$ are then exactly the sets of the form $X \times F$, where
$F$ is closed in $Y$.  This implies that $||X \times Y|| = ||Y||$, and
since $||X||=1$, that is equal to $||X|| \otimes ||Y||$.  We reason
similarly if $\beta_1=0$.

If $\alpha_1$ is a successor ordinal, then
$||X|| = \omega^{\alpha_1-1} \times \omega$.  For every $m \in \nat$,
there is a proper closed subset $F_m$ of $X$ such that
$\rk_{\Hoare_0 X} (F_m) = \omega^{\alpha_1-1} \times m$, since
$\omega^{\alpha_1-1} \times m < ||X||$.  By Lemma~\ref{lemma:rk:C},
$||F_m|| = \omega^{\alpha_1-1} \times m$.  By induction hypothesis,
$||F_m \times Y|| \geq (\omega^{\alpha_1-1} \times m) \otimes
\omega^{\beta_1} = (\bigoplus_{i=1}^m \omega^{\alpha_1-1}) \otimes
\omega^{\beta_1} = \bigoplus_{i=1}^m \omega^{(\alpha_1-1) \oplus
  \beta_1} = \omega^{(\alpha_1-1) \oplus \beta_1} \times m$.  Since
$||X \times Y|| \geq ||F_m \times Y||$, and taking suprema over all
$m \in \nat$,
$||X \times Y|| \geq \omega^{(\alpha_1-1) \oplus \beta_1} \times
\omega$.  Now
$\omega^{(\alpha_1-1) \oplus \beta_1} \times \omega =
\omega^{((\alpha_1-1) \oplus \beta_1)+1}$, and
$((\alpha_1-1) \oplus \beta_1) + 1 = ((\alpha_1-1) \oplus \beta_1)
\oplus 1 = ((\alpha_1-1) \oplus 1) \oplus \beta_1 = \alpha_1 \oplus
\beta_1$.

The argument is symmetric if $\beta_1$ is a successor ordinal.

Finally, if $\alpha_1$ and $\beta_1$ are both limit ordinals, then
$\alpha_1 \oplus \beta_1$ is also a limit ordinal, as one can see from
its Cantor normal form.  Hence $\alpha_1 \oplus \beta_1$ is the
supremum of all ordinals $\gamma < \alpha_1 \oplus \beta_1$.  Since
$\gamma \mapsto \omega^\gamma$ is continuous,
$\omega^{\alpha_1} \otimes \omega^{\beta_1} = \omega^{\alpha_1 \oplus
  \beta_1} = \sup_{\gamma < \alpha_1 \oplus \beta_1} \omega^\gamma$.
In order to show that
$||X \times Y|| \geq \omega^{\alpha_1} \otimes \omega^{\beta_1}$, it
therefore suffices to show that $||X \times Y|| \geq \omega^\gamma$
for every $\gamma < \alpha_1 \oplus \beta_1$.

Using the inductive definition of the natural sum,
$\alpha_1 \oplus \beta_1$ is the smallest ordinal strictly larger than
$\alpha \oplus \beta_1$ for every $\alpha < \alpha_1$ and
$\alpha_1 \oplus \beta$ for every $\beta < \beta_1$.  Since
$\gamma < \alpha_1 \oplus \beta_1$, we must therefore have
$\gamma \leq \alpha \oplus \beta_1$ for some $\alpha < \alpha_1$, or
$\gamma \leq \alpha_1 \oplus \beta$ for some $\beta < \beta_1$.  In
the first case, $\omega^\alpha < \omega^{\alpha_1} = ||X||$, so there
is a proper closed subset $F_\alpha$ of $X$ such that
$\rk_{\Hoare_0 X} (F_\alpha) = \omega^\alpha$; namely,
$||F_\alpha|| = \omega^\alpha$, by Lemma~\ref{lemma:rk:C}.  By
induction hypothesis,
$||F_\alpha \times Y|| \geq \omega^{\alpha \oplus \beta_1} \geq
\omega^\gamma$, and therefore $||X \times Y|| \geq \omega^\gamma$.
Similarly, in the second case,
$||X \times Y|| \geq \omega^{\alpha_1 \oplus \beta} \geq
\omega^\gamma$.  \qed

Combining Proposition~\ref{prop:stature:prod:ubound} and
Proposition~\ref{prop:stature:prod:lbound}, we finally obtain the
desired result.
\begin{theorem}
  \label{thm:stature:prod}
  For all Noetherian spaces $X$ and $Y$, $||X \times Y|| = ||X||
  \otimes ||Y||$.
\end{theorem}

\section{Hoare powerspaces and powersets}
\label{sec:hoare-powerspaces}

Here is another example of Noetherian spaces which do not arise from
wqos in general.  The \emph{Hoare powerspace} $\HV X$ of $X$ is just
its space of non-empty closed subsets $\Hoare X$, with the so-called
\emph{lower Vietoris topology}.  We also consider the \emph{lifted
  Hoare powerspace} $\HVz X$, which also includes the empty set.  A
subbase of the lower Vietoris topology is given by sets $\Diamond U$,
defined as the set of those closed sets $F$ that intersect $U$, where
$U$ ranges over the open subsets of $X$.

It was observed in \cite{Gou-lics07} that $\HV X$ and $\HVz X$ are
Noetherian for every Noetherian space $X$.  That may seem surprising
at first, considering that the specialization ordering of each one is
inclusion, and that the inclusion ordering on the downwards-closed
subsets of a wqo $P$ is \emph{not} in general a wqo \cite{Rado:wqo}.

The trick is that the lower Vietoris topology is in general strictly
coarser than the Alexandroff topology.  In fact, the lower Vietoris
topology coincides with the upper topology, since the complement of
$\Diamond U$ is equal to $\dc F$, where $F$ is the complement of $U$.
In the sequel, we will write $\Box F$ instead of $\dc F$ for the set
of closed subsets of $F$.  This 
will dispel any ambiguity, since $\dc F$ is also accepted notation for
the downward closure of $F$ in $X$, not $\Hoare_0 X$.

\subsection{The sobrification rank and stature of $\HVz X$}
\label{sec:dimens-stat-hvz}

Schalk observed that every up-complete sup-semilattice, namely every
poset with suprema of all non-empty families, is sober in its upper
topology \cite[Proposition~1.7]{Schalk:PhD}.  \blue{$\Hoare X$ and
  $\Hoare_0 X$ are up-complete sup-semilattices, and the lower
  Vietoris topology is the upper topology of inclusion.}  Hence both
$\HV X$ and $\HVz X$ are sober.  This makes the following a
triviality.
\begin{theorem}
  \label{thm:dim:H}
  For every Noetherian space $X$, $\sob {\HVz X} = ||X||+1$ and
  $\rsob {\HVz X} = ||X||$.
\end{theorem}
\proof
Since $\HVz X$ is sober, $\sob {\HVz X}$ is just the ordinal rank of
$\Hoare_0 X$, namely $||X||+1$.
\qed

The stature of $\HVz X$ is much more elusive.
\begin{proposition}
  \label{prop:stature:H}
  For every Noetherian space $X$,
  $1+||X|| \leq ||\HVz X|| \leq \omega^{||X||}$.
\end{proposition}
\proof By Proposition~\ref{prop:dim:stature}~(1) and
Theorem~\ref{thm:dim:H}, $1+||X|| = 1+\rsob {\HVz X} \leq ||\HVz X||$.
$\HVz X$ is not empty, and has exactly one component, which is
$\Box X = \HVz X$.  By Proposition~\ref{prop:dim:stature}~(3),
$||\HVz X|| \leq \omega^{\rsob {\HVz X}} = \omega^{||X||}$.  \qed

We claim that those lower and upper bounds are tight in general.
Before we give substance to this claim, we observe that several other
Noetherian spaces are related to $\HVz X$.

\subsection{Powersets, and finitary variants}
\label{sec:powers-finit-vari}

One is the \emph{powerset} $\pow X$ of $X$, with a topology that we
will still call the lower Vietoris topology, whose subbasic open sets
we will still write as $\Diamond U$, and which now denote
$\{A \in \pow X \mid A \cap U \neq \emptyset\}$.  Another one is
$\Hoare_\fin X$, the subspace of $\HVz X$ consisting of all the
finitary closed subsets $\dc \{x_1, \cdots, x_n\}$ of $X$.  Finally,
there is the finitary powerset $\pow_\fin X$, which is the subspace of
$\pow X$ consisting of its finite subsets.

The specialization ordering of $\Hoare_\fin X$, just like $\HVz X$, is
inclusion.  Indeed, the closure of any point $F \in \Hoare_\fin X$ is
$\Box F$, which is also the downward closure of $F$ under inclusion.
The specialization preordering of $\pow_\fin X$, just like $\pow X$,
is inclusion of closures.  This was proved for $\pow X$ in
\cite[Lemma~4.9]{FGL:partI}, and follows from the fact that the
closure of $\{A\}$ in $\pow_\fin X$ (resp., $\pow X$) is
$\Box cl (A)$, where $cl (A)$ denotes the closure of $A$ in $X$.

\begin{remark}
  \label{rem:Pfin}
  While $\HVz X$ and $\pow X$ do not arise from wqos, it so happens
  that $\Hoare_\fin X$ and $\pow_\fin X$ do, provided that $X$ is wqo
  (in its Alexandroff topology).  Indeed, it suffices to verify that
  the lower Vietoris topology on each coincides with the Alexandroff
  topology of their specialization preorderings.  To this end, it is
  enough to show that the upward closure of each point $F$ is open in
  the lower Vietoris topology.  For $\Hoare_\fin X$, the
  specialization preordering is inclusion, and the upward closure of
  $F \eqdef \dc \{x_1, \cdots, x_n\}$ is
  $\Diamond (\upc x_1) \cap \cdots \cap \Diamond (\upc x_n)$.  For
  $\pow_\fin X$, the specialization ordering $\leq^\flat$ is given by
  $A \leq^\flat B$ if and only if the closure of $A$ is included in
  the closure of $B$; when $X$ is Alexandroff, that is equivalent to
  the fact that the downward closure of $A$ is included in that of
  $B$, or equivalently that every element of $A$ is smaller than or
  equal to some element of $B$.  Then the upward closure of
  $\{x_1, \cdots, x_n\}$ is
  $\Diamond (\upc x_1) \cap \cdots \cap \Diamond (\upc x_n)$, as with
  $\Hoare_\fin X$.
\end{remark}

\begin{proposition}
  \label{prop:Pfin}
  For every topological space $X$, the following maps are full,
  continuous, and Skula dense:
  \begin{enumerate}
  \item the functions $cl$ that map every set to its closure in $X$,
    from $\pow X$ to $\HVz X$, and from $\pow_\fin X$ to
    $\Hoare_\fin X$;
  \item the inclusion maps from $\Hoare_\fin X$ into $\HVz X$ and from
    $\pow_\fin X$ into $\pow X$.
  \end{enumerate}

  In particular, for every Noetherian space $X$,
  $\sob {\pow_\fin X} = \sob {\Hoare_\fin X} = \sob {\pow X} = \sob
  {\HVz X} = ||X||+1$, and $||\pow_\fin X|| = ||\Hoare_\fin X|| = ||\pow
  X|| = ||\HVz X||$.
\end{proposition}
\proof We start with $cl$.  Let us reserve the notation $\Diamond U$
for the subbasic open subsets of $\HVz X$ (resp., $\Hoare_\fin X$), and
let us use $\Diamond_\pow U$ to denote the corresponding subbasic open
subset of $\pow X$ (resp., $\pow_\fin X$).

For every open subset $U$ of $X$, $cl^{-1} (\Diamond U)$ is the set of
subsets (resp., finite subsets) $A$ of $X$ such $cl (A)$ intersects
$U$.  Since $cl (A)$ intersects $U$ if and only if $A$ does,
$cl^{-1} (\Diamond U) = \Diamond_\pow U$.  This shows that $cl$ is
continuous.  Also, $cl$ is full since every open subset of $\pow X$
(resp., $\pow_\fin X$) is of the form
$\bigcup_{i \in I} \bigcap_{j \in J_i} \Diamond_\pow U_{ij}$, where
each $J_i$ is a finite set; and that is equal to
$\bigcup_{i \in I} \bigcap_{j \in J_i} cl^{-1} (\Diamond U_{ij}) =
cl^{-1} (\bigcup_{i \in I} \bigcap_{j \in J_i} \Diamond U_{ij})$.
For any two open subsets $\bU$ and $\bV$ of $\HVz
X$ (resp., $\Hoare_\fin X$) such that $cl^{-1} (\bU) = cl^{-1} (\bV)$,
every $F \in \bU$ is also in $cl^{-1} (\bU)$ since $F = cl (F)$,
hence in $cl^{-1} (\bV)$; so $F = cl (F')$ for some $F' \in \bV$.
But $F' \in \bV$ implies that $F'$ is closed in $X$, so $F=F'$, and
therefore $F$ is in $\bV$.  This shows that $\bU$ is included in
$\bV$, and the reverse inclusion is proved similarly.  It follows that
$cl$ is Skula dense.

The inclusion maps from $\Hoare_\fin X$ into $\HVz X$ and from
$\pow_\fin X$ into $\pow X$ are topological embeddings by definition,
hence are full and continuous.  Let us write $i$ for any of those
maps.  In order to show that it is Skula dense, we now reserve the
notation $\Diamond U$ for subbasic open subsets of $\HVz X$, resp.\
$\pow X$, and write $\Diamond_\fin X$ for the corresponding sets in
$\Hoare_\fin X$, resp.\ $\pow_\fin X$.  For any open subset
$\bU \eqdef \bigcup_{i \in I} \bigcap_{j \in J_i} \Diamond U_{ij}$ of
$\HVz X$ (resp., $\pow X$), where each set $J_i$ is finite and each
$U_{ij}$ is open in $X$, $i^{-1} (\bU)$ is equal to
$\bigcup_{i \in I} \bigcap_{j \in J_i} \Diamond_\fin U_{ij}$.  Given
another open subset $\bV$, such that $i^{-1} (\bU) = i^{-1} (\bV)$, we
claim that $\bU=\bV$; by symmetry, it suffices to show
$\bU \subseteq \bV$.  Let $F$ be any element of $\bU$.  There is an
index $i \in I$ such that $F$ intersects $U_{ij}$ for every
$j \in J_i$, say at $x_j$.  Let $F' \eqdef \dc \{x_j \mid j \in J_i\}$
(resp., $F' \eqdef \{x_j \mid j \in J_i\}$).  Then $F'$ is in
$\bigcap_{j \in J_i} \Diamond_\fin U_{ij} \subseteq i^{-1} (\bU) =
i^{-1} (\bV)$, so $i (F') = F'$ is in $\bV$.  We now observe that
$F' \subseteq F$ (resp., $cl (F') \subseteq cl (F)$), and that $\bV$
is upwards-closed in the specialization preordering of $\HVz X$
(resp., $\pow X$), so $F$ is in $\bV$.

The remaining claims follow from
Lemma~\ref{lemma:stature:mono:inj}~(4) and Theorem~\ref{thm:dim:H}.
\qed


\subsection{The bounds on $||\HVz X||$ are tight}
\label{sec:bounds-hvz-x}

We proceed through a series of examples.
\begin{example}
  \label{exa:stature:H:lbound}
  The lower bound of Proposition~\ref{prop:stature:H} is attained.
  Consider $X \eqdef \alpha$, where $\alpha$ is an ordinal, with its
  Alexandroff topology.  Then $||X|| = \alpha$ by
  Lemma~\ref{lemma:ordinal}~(1).  The closed subsets of $X$ are
  themselves totally ordered by inclusion, and every non-empty closed
  subset is irreducible.  Hence
  $\Hoare_0 (\HVz X) = \Sober (\HVz X) \cup \{\emptyset\}$, so that
  $||\HVz X|| = |\Sober (\HVz X) \cup \{\emptyset\}|-1 = (1+|\Sober
  (\HVz X)|)-1$.  Since $\HVz X$ is sober, it is homeomorphic to its
  sobrification, so $||\HVz X|| = 1+|\alpha| = 1+\alpha = 1+||X||$.
  This result can also be obtained by verifying that $\HVz X$ is equal
  to $(1+\alpha)+1$ if $\alpha$ is a limit ordinal, to $1+\alpha$
  otherwise, and with the upper topology in both cases; then apply
  Lemma~\ref{lemma:ordinal}~(2).
\end{example}

\begin{example}
  \label{exa:stature:H:fixpoint}
  Let us consider $X \eqdef \alpha$, where $\alpha$ is any ordinal,
  with its Scott topology.  As in Example~\ref{exa:stature:H:lbound},
  every non-empty closed subset is irreducible; the non-empty closed
  subsets are the subsets $\dc \beta$ with $\beta < \alpha$, plus
  $\alpha$ itself if $\alpha$ is zero or a limit ordinal.  Hence
  $\HVz X$ can be equated to $1+\alpha$ if $\alpha$ is finite, to
  $\alpha+1$ if $\alpha$ is a limit ordinal, and to $\alpha$
  otherwise.  The topology is the upper topology, which coincides with
  the Scott topology.  By Lemma~\ref{lemma:ordinal}~(2), $||\HVz X||$
  is equal to $1+\alpha$ if $\alpha$ is finite, to $\alpha$ if
  $\alpha$ is a limit ordinal, and to $\alpha-1$ otherwise.  Using
  Lemma~\ref{lemma:ordinal}~(2), we check that this is equal to
  $1+||X||$ in all cases.  This has the curious consequence that, if
  $\alpha$ is an infinite ordinal with its Scott topology, then
  $||X|| = ||\HVz X|| = ||\HVz (\HVz X)|| = \cdots = ||\HVz^n X||$ for
  every $n \in \nat$.
\end{example}

\begin{remark}
  \label{rem:stature:H:eps}
  If $||X||$ is a \emph{critical} ordinal $\epsilon$, namely an ordinal such
  that $\epsilon = \omega^\epsilon$, then the lower and
  upper bounds of Proposition~\ref{prop:stature:H} match.  It follows
  that $||\HVz X|| = \epsilon = ||X||$ in this case.  As in
  Example~\ref{exa:stature:H:fixpoint}, this implies that
  $||\HVz^n X|| = ||X|| = \epsilon$ for every $n \in \nat$.
\end{remark}

\begin{example}
  \label{exa:stature:H:fin}
  Here is a case where $||\HVz X||$ is strictly between the lower and
  upper bounds of Proposition~\ref{prop:stature:H}.  Let $X$ be a
  finite set of cardinality $n$, with its discrete topology.  Then
  $||X||=n$ by Lemma~\ref{lemma:finite}.  $\Hoare_0 X$ is the powerset
  of $X$.  The specialization preordering of $\HVz X$ is inclusion,
  which is an ordering, so $\HVz X$ is $T_0$.  By
  Lemma~\ref{lemma:finite}, $||\HVz X|| = 2^n$.
\end{example}

We will use the following to show that the upper bound of
Proposition~\ref{prop:stature:H} is attainable.

\begin{example}
  \label{exa:stature:H:upper}
  Given any space $X$, let $s (X)$ denote the space obtained by
  forming the sum of countably many disjoint copies of $X$, and adding
  a fresh top element $\top$.  Formally, the closed subsets of $s (X)$
  are the disjoint unions $F_1 + \cdots + F_n$ of closed sets from
  finitely many copies of
  $X$, plus $s (X)$ itself.  By Proposition~\ref{prop:stature:sum},
  $||F_1 + \cdots + F_n|| = \bigoplus_{i=1}^n ||F_i||$; with $n$
  fixed, but letting $F_1$, \ldots, $F_n$ vary, the largest value
  taken by
  that sum is $||X|| \otimes n$.  If $||X|| = \omega^\alpha$, that is
  equal to $\omega^\alpha \times n$; taking suprema over all $n$, we
  obtain that $||s (X)|| = \omega^{\alpha+1}$. For every $n \in \nat$,
  there is an obvious embedding of $(\HVz X)^n$ into $\HVz {(s (X))}$,
  which maps every $n$-tuple $(F_1, \cdots, F_n)$ to
  $F_1 + \cdots + F_n$, each $F_i$ being located in the $i$th copy of
  $X$. By Lemma~\ref{lemma:stature:mono:inj}~(3) and
  Theorem~\ref{thm:stature:prod},
  $||\HVz {(s (X))}|| \geq \bigotimes_{i=1}^n ||\HVz X||$. If
  $||\HVz X|| \geq \omega^{\omega^\beta}$, we obtain that
  $||\HVz {(s (X))}|| \geq \omega^{\omega^\beta \times n}$. Since $n$
  is arbitrary, $||\HVz {(s
    (X))}|| \geq \omega^{\omega^{\beta+1}}$.

  We iterate this construction, starting from a Noetherian space $X$
  such that $||X|| = \omega$, for example $\omega$ itself or any
  infinite set with the cofinite topology.  By taking $\alpha \eqdef
  1$ and $\beta \eqdef 0$ (since $||\HVz X|| \geq 1+||X|| = \omega$ by
  Proposition~\ref{prop:stature:H}), we obtain that 
  $||s (X)|| = \omega^2$ and $||\HVz {(s (X))}|| \geq \omega^\omega$.
  We can now take $\alpha \eqdef 2$ and $\beta \eqdef 1$, and
  obtain $||s^2 (X)|| = \omega^3$ and $||\HVz {(s^2 (X))}|| \geq
  \omega^{\omega^2}$.
  In general, $||s^k (X)|| = \omega^{k+1}$ and $||\HVz {(s^k (X))}|| \geq
  \omega^{\omega^k}$ for every $k \in \nat$.

  We now build a form of limit of the spaces $s^k (X)$, $k \in \nat$.
  For every space $Y$, let $i_Y \colon Y \to s (Y)$ map every
  $y \in Y$ to $y$ itself in copy number $0$ of $Y$ inside $s (Y)$. We
  note that $i_Y$ is an embedding, and a closed map. Let
  $s^\omega (X)$ be the quotient of the disjoint sum $\coprod_{k \in
    \nat} s^k (X)$ by the smallest equivalence relation $\equiv$ that
  equates $y$ with $i_{s^k (X)} (y)$ for every $y \in s^k (X)$. For
  each $k \in \nat$ and each closed subset $F$ of $s^k (X)$, let
  $[F]_k$ be the set of equivalence classes of points of
  $F \subseteq s^k (X)$ modulo $\equiv$. We note that
  $[F]_k = [i_{s^k (X)} (F)]_{k+1}$, so that we can always write any
  finite collection of sets $[F]_k$ with the same index $k$. It
  follows that the collection of such sets is closed under finite
  unions. It is also closed under arbitrary non-empty intersections,
  which in fact reduce to finite non-empty intersections in some
  subspace $s^k (X)$. Therefore, the collection of sets $[F]_k$, where
  $k \in \nat$ and $F$ is closed in $X$, plus the whole space
  $s^\omega (X)$ itself, form a Noetherian topology on
  $s^\omega (X)$. Since $s^k (X)$ embeds into $s^\omega (X)$ for every
  $k \in \nat$, $||\HVz {(s^\omega (X))}|| \geq \omega^{\omega^k}$ for
  every $k \in \nat$, and therefore
  $||\HVz {(s^\omega (X))}|| \geq \omega^{\omega^\omega}$.
  The family ${([s^k (X)]_k)}_{k \in \nat}$ forms a cofinal family of
  proper closed sets, so by Proposition~\ref{prop:limit},
  $||s^\omega (X)|| = \sup_{k \in \nat} (\omega^{k+1}+1) =
  \omega^\omega$.
  In particular, $||\HVz {(s^\omega (X))}||$ is larger than or equal
  to, and therefore equal to the upper bound $\omega^{||s^\omega
    (X)||}$ (namely, $\omega^{\omega^\omega}$) given
  in Proposition~\ref{prop:stature:H}.

\end{example}

\section{Finite words}
\label{sec:bf-finite-words}

\subsection{The structure of $\Sober (X^*)$}
\label{sec:structure-sober-x}

Given any set $X$, let $X^*$ be the set of finite words on $X$.
We write $\epsilon$ for the empty word, and $ww'$ for the
concatenation of two words $w$ and $w'$.  Given two subsets $A$ and
$B$ of $X^*$, we also write $AB$ for $\{ww' \mid w \in A, w' \in B\}$.

When $X$ is a topological space, we give $X^*$ the \emph{word
  topology}, which is defined in \cite[Definition~9.7.26]{JGL:topo} as
the topology generated by subsets of the form
$X^* U_1 \allowbreak X^* U_2 \allowbreak X^* \ldots \allowbreak X^*
U_n X^*$, where $n \in \nat$, and $U_1$, $U_2$, \ldots, $U_n$ range
over open subsets of $X$.

For every Noetherian space $X$, $X^*$ is Noetherian
\cite[Theorem~9.7.33]{JGL:topo}.  Let us write $\leq$ for the
specialization preordering of $X$.  The specialization preordering of
$X^*$ is the \emph{word embedding} quasi-ordering $\leq^*$, defined by
$w \leq^* w'$ if and only if one can obtain $w'$ from $w$ by
increasing some of its letters and inserting arbitrarily many letters
at any position \cite[Exercise~9.7.29]{JGL:topo}.  This is the
familiar preordering at the heart of Higman's Lemma, which says that
for every wqo $\leq$, $\leq^*$ is a wqo \cite{Higman:Lemma}.  That
$X^*$ is Noetherian for every Noetherian space $X$ is a topological
generalization of this result, in the sense that if $X$ is Alexandroff
and Noetherian (i.e., wqo), then so is $X^*$
\cite[Exercise~9.7.30]{JGL:topo}.

A \emph{word-product} $\tP$ on a space $X$ is any expression of the
form $\tA_1 \tA_2 \cdots \tA_n$, where $n \in \nat$, and each $\tA_i$
is an \emph{atomic expression}, either $F_i^*$ with $F_i$ closed in
$X$, or $C_i^?$ with $C_i$ irreducible closed in $X$.  The semantics
of an atomic expression $\tA$ will be written $\bA$.  The semantics of
$C_i^?$ is the collection of words containing at most one letter, and
that letter must be in $C_i$.  That is sometimes written
$C_i+\epsilon$, or $C_i \cup \{\epsilon\}$, by abuse of language.
\blue{The semantics of $F_i^*$ is the set of finite words, of
  arbitrary lengths, whose letters are in $F_i$, and the semantics of
  $\tA_1 \tA_2 \cdots \tA_n$ is
  $\mathbf A_1 \mathbf A_2 \cdots \mathbf A_n$.  For example, if we
  take $X$ to be $\{a, b, c\}$ with the discrete topology, then
  $\{a, b\}^* \{c\}^? \{a\}^? \{b\}^? \{a, c\}^*$ is the set of words
  $w_1 w w_2$ where $w_1$ is a word of arbitrary length whose letters
  are among $\{a, b\}$, $w_2$ is any word whose letters are among
  $\{a, c\}$, and $w$ is $cab$ or any of its subwords $ab$, $cb$,
  $ca$, $c$, $a$, $b$ and $\epsilon$.}

When $n=0$, $\tP$ is abbreviated as $\epsilon$, and denotes the
one-element set $\{\epsilon\}$.  In general, we write $\bP$ for the
semantics of $\tP$.  Note that $\bP$ always contains the empty word
$\epsilon$.

We call \emph{word-SRE} any finite sum of word-products, where sum is
interpreted as union.  ``SRE'' stands for ``simple regular
expression'' \cite{ABJ:SRE}.  It turns out that, given any Noetherian
space $X$, the closed subsets of $X^*$ are exactly the semantics
of word-SREs, and the irreducible closed subsets of $X^*$ are exactly
the semantics of word-products \cite[Proposition~6.14]{FGL:partI}.
This extends the corresponding result on wqos, due to Kabil and Pouzet
\cite{KP:age:mots}, which itself extends a previous result of Jullien
on words on a finite set, ordered by equality
\cite[chapitre~VI]{Jullien:PhD}.

The inclusion ordering on $\Sober (X^*)$ is characterized as follows.
First, we observe that inclusion of atomic expressions is
characterized by:
\begin{enumerate}
\item $C^? \subseteq {C'}^?$ if and only if $C \subseteq C'$;
\item $F^* \subseteq {F'}^*$ if and only if $F \subseteq F'$;
\item $C^? \subseteq {F'}^*$ if and only if $C \subseteq F'$;
\item $F^* \subseteq {C'}^?$ if and only if $F$ is empty.
\end{enumerate}
Then we have the following.
\begin{lemma}[Lemmata~7.9 and 7.10, \cite{FGL:partI}]
  \label{lemma:6.10}
  Given two word-products $\tP$ and $\tP'$ on a Noetherian space $X$,
  $\bP \subseteq \bP'$ if and only if $\tP$ is $\epsilon$, or
  $\tP' = \epsilon$ and $\tP$ is a product of atomic expressions all
  equal to $\emptyset^*$, or $\tP$ can be written as $\tA_1 \tQ$ and
  $\tP'$ as $\tA'_1 \tQ'$, where $\tA_1$ and $\tA'_1$ are atomic
  expressions, and one of the following occurs:
  \begin{enumerate}
  \item $\bA_1 \not\subseteq \bA'_1$ and $\bP \subseteq \bQ'$;
  \item or $\tA_1$ is of the form $C^?$, $\tA'_1$ is of the form ${C'}^?$,
    $C \subseteq C'$, and $\bQ \subseteq \bQ'$;
  \item or $\tA'_1$ is of the form ${F'}^*$, $\bA_1 \subseteq \bA'_1$, and
    $\bQ \subseteq \bP'$;
  \item or $\tA_1 = \emptyset^*$ and $\bQ \subseteq \bP'$.
  \end{enumerate}
\end{lemma}
Equivalently, inclusion of word-products is axiomatized in a sound a
complete way by the following rules:
\[
  \begin{array}{ccccc}
    \begin{prooftree}
      \strut
      \justifies
      \epsilon \subseteq \bP'
      \using (0)
    \end{prooftree}
    &
      \begin{prooftree}
        \bA_1 \bQ \subseteq \bQ'
        \justifies
        \bA_1 \bQ \subseteq \bA'_1 \bQ'
        \using (1)
      \end{prooftree}
    &
      \begin{prooftree}
        \bQ \subseteq \bQ'
        \justifies
        C^? \bQ \subseteq {C'}^? \bQ'
        \using (2)
      \end{prooftree}
    &
      \begin{prooftree}
        \bQ \subseteq {F'}^* \bQ'
        \justifies
        \bA_1 \bQ \subseteq {F'}^* \bQ'
        \using (3)
      \end{prooftree}
    &
      \begin{prooftree}
        \bQ \subseteq \bP'
        \justifies
        \emptyset^* \bQ \subseteq \bP'
        \using (4)
      \end{prooftree}
    \\
    & \text{if }\bA_1 \not\subseteq \bA'_1
    & \text{if }C \subseteq C'
    & \text{if }\bA_1 \subseteq {F'}^*
    &
  \end{array}
\]

A word-product $\tA_1 \tA_2 \cdots \tA_n$ is \emph{reduced} if and
only if no $\tA_i$ is equal to $\emptyset^*$, and for every $\tA_i$ of
the form $F^*$, $F^*$ is not included in $\bA_{i+1}$ (if $i<n$) and
not included in $\bA_{i-1}$ (if $i > 1$).  It is easy to see that
every word-product can be rewritten into a reduced word-product with the
same semantics.  One can also show that reduced word-products are
canonical forms for irreducible closed subsets of $X^*$, namely that
two reduced word-products denote the same irreducible closed set if
and only if they are syntactically equal.  The proof is identical to
the corresponding result on wqos, see
\cite[Theorem~4.22]{GHKNS-til2020}; but we will not make use of that
fact.

\subsection{The sobrification rank of $X^*$}
\label{sec:dimension-x}

We will need the following ordinal adjustment operation. We recall
that a critical ordinal is an ordinal $\epsilon$ such that
$\epsilon = \omega^\epsilon$. \begin{definition}[$\adjsw\alpha$]
  \label{defn:adjsw}
  For every ordinal $\alpha$, the ordinal
  $\adjsw\alpha$ is defined as $\alpha+1$ if $\alpha =
  \epsilon+n$ for some critical ordinal
  $\epsilon$ and some natural number $n$, as $\alpha$ otherwise.
\end{definition}
Our goal in this subsection is to show that, for every non-empty
Noetherian space $X$, $\sob {X^*} = \omega^{\adjsw{||X||}}+1$.

\begin{lemma}
  \label{lemma:alpha':inc}
  The map $\alpha \mapsto \adjsw\alpha$ is strictly increasing.
\end{lemma}
\proof Let $\alpha_1 <
\alpha_2$.  We need to consider four cases in order to show that
$\adjsw{\alpha_1} <
\adjsw{\alpha_2}$, but only one is non-trivial, namely when $\alpha_1
= \epsilon+n$ for some critical ordinal $\epsilon$ and some $n \in
\nat$, and $\alpha_2 > \alpha_1$ is not of that form.  Then $\alpha_2
\geq \epsilon + \omega$, so $\adjsw{\alpha_2} = \alpha_2 >
\epsilon+n+1 = \adjsw{\alpha_1}$.  \qed

\begin{lemma}
  \label{lemma:alpha':omega}
  For every ordinal $\alpha$, $\alpha < \omega^{\adjsw\alpha}$.
\end{lemma}
\proof By induction on $\alpha$, we have $\alpha \leq \omega^\alpha$.
The inequality is strict unless $\alpha$ is critical, by definition.
In particular, if $\alpha$ is not of the form $\epsilon+n$ with
$\epsilon$ critical and $n \in \nat$, then
$\alpha < \omega^\alpha = \omega^{\adjsw{\alpha}}$ \blue{(since
  $\alpha^\circ=\alpha$ in this case)}.  When $\alpha = \epsilon+n$
where $\epsilon$ is critical and $n \in \nat$, then
$\omega^{\adjsw\alpha} = \omega^{\epsilon+n+1} = \omega^\epsilon
\times \omega^{n+1} = \epsilon \times \omega^{n+1}$.  Since
$\omega^{n+1} \geq \omega > 2$, this is larger than
$\epsilon \times 2 = \epsilon + \epsilon > \epsilon + n = \alpha$.
\qed

\begin{remark}
  \label{rem:adjsw}
  The map $\alpha \mapsto \adjsw\alpha$ is in fact the smallest
  strictly increasing map such that $\alpha < \omega^{\adjsw\alpha}$
  for every ordinal $\alpha$, as one can check.
\end{remark}

From now on, we use the notation $||F||$ instead of
$\rk_{\Hoare_0 X} (F)$, profiting from Lemma~\ref{lemma:rk:C}.

\begin{lemma}
  \label{lemma:rk:stat}
  Let $X$ be a Noetherian space.  For all $C \in \Sober X$,
  $F, F' \in \Hoare_0 X$,
  \begin{enumerate}
  \item if $C \subseteq F'$ then $1+\rk_{\Sober X} (C) <
    \omega^{\adjsw{||F'||}}$;
  \item if $F \subsetneq F'$ then
    $\omega^{\adjsw{||F||}} < \omega^{\adjsw{||F'||}}$.
  \end{enumerate}
\end{lemma}
\proof (1) The ordinal $1+\rk_{\Sober X} (C)$ is the rank of $C$ in
$\Sober X \cup \{\emptyset\}$.  Since that is included in
$\Hoare_0 X$,
$1+\rk_{\Sober X} (C) \leq \rk_{\Hoare_0 X} (C) \leq \rk_{\Hoare_0 X}
(F') = ||F'||$.  We conclude since $||F'|| < \omega^{\adjsw{||F'||}}$
by Lemma~\ref{lemma:alpha':omega}.

(2) The map $\gamma \mapsto \omega^\gamma$ is strictly monotonic, and
so is the map $\alpha \mapsto \adjsw\alpha$ by
Lemma~\ref{lemma:alpha':inc}.  \qed

\begin{lemma}
  \label{lemma:X*:phi}
  Let $X$ be a Noetherian space, and let us define $\varphi (C^?)$ as
  $1+\rk_{\Sober X} (C)$ for every $C \in \Sober X$, $\varphi (F^*)$
  as $\omega^{\adjsw{||F||}}$ for every $F \in \Hoare_0 X$, and
  $\varphi (\tP)$ as $\bigoplus_{i=1}^n \varphi (\tA_i)$ for every
  reduced word-product $\tP \eqdef \tA_1 \cdots \tA_n$.  For all
  reduced word-products $\tP$ and $\tP'$, $\bP \subseteq \bP'$ implies
  $\varphi (\tP) \leq \varphi (\tP')$, and if the former inclusion is
  strict, then so is the latter inequality.
\end{lemma}
\proof We proceed by induction on the sum of the lengths of $\tP$ and
of $\tP'$.

When $\tP = \epsilon$, the claim is clear since
$\varphi (\epsilon) = 0$, and $\varphi (\tP') \neq 0$ for every
reduced word-product $\tP' \neq \epsilon$. We therefore assume that
$\tP \neq \epsilon$ in the sequel.

In particular, $\tP$ is of the form $\tA_1 \cdots \tA_n$ with
$n \geq 1$.  Each $\bA_i$ contains a non-empty word, using the fact
that no $\tA_i$ is of the form $\emptyset^*$.  Therefore $\bP$
contains a non-empty word.  Since $\bP \subseteq \bP'$, it follows
that $\tP'$ cannot be equal to $\epsilon$.  Let us write $\tP'$ as
$\tA'_1 \tQ'$.  We also write $\tP$ as $\tA_1 \tQ$.

If $\bA_1 \not\subseteq \bA'_1$, then only clause (1) of
Lemma~\ref{lemma:6.10} can have been used to infer
$\bP \subseteq \bP'$, so $\bP \subseteq \bQ'$.  By induction
hypothesis, $\varphi (\tP) \leq \varphi (\tQ')$; $\varphi (\tQ')$ is
strictly less than
$\varphi (\tP') = \varphi (\tA'_1) \oplus \varphi (\tQ')$, so
$\varphi (\tP) < \varphi (\tP')$.  Note that the inequality is strict
because $\varphi (\tA'_1)$ cannot be equal to $0$.

Let us now assume that $\bA_1 \subseteq \bA'_1$.

If $\tA_1=C^?$ and $\tA'_1={C'}^?$, then only clause (2) can have been
used to derive $\bP \subseteq \bP'$.  In this case,
$\bQ \subseteq \bQ'$, so $\varphi (\tQ) \leq \varphi (\tQ')$ by
induction hypothesis.  Additionally, $C \subseteq C'$, so
$\varphi (\tA_1) = 1+\rk_{\Sober X} (C) \leq 1+\rk_{\Sober X} (C') =
\varphi (\tA'_1)$.  It follows that
$\varphi (\tP) = \varphi (\tA_1) \oplus \varphi (\tQ) \leq \varphi
(\tA'_1) \oplus \varphi (\tQ') = \varphi (\tP')$.  If additionally
$\bP \neq \bP'$, then one of the inclusions $C \subseteq C'$ and
$\bQ \subseteq \bQ'$ must be strict.  In the first case,
$1+\rk_{\Sober X} (C) < 1+\rk_{\Sober X} (C')$, and in the second
case, $\varphi (\tQ) < \varphi (\tQ')$ by induction hypothesis, so
that in both cases $\varphi (\tP) < \varphi (\tP')$.

It remains to deal with the case where $\bA_1 \subseteq \bA'_1$ and
$\tA'_1$ is of the form ${F'}^*$ for some $F' \in \Hoare_0 X$.  Since
$\tP$ is reduced, clause (4) does not apply, so only clause (3)
applies.  We apply clause (3) repeatedly until we can no longer.  This
allows us to write $\tP$ as $\tA_1 \tA_2 \cdots \tA_k \tR$, for some
word-product $\tR$, where $k \geq 1$ is largest so that
$\bA_1, \cdots, \bA_k \subseteq {F'}^*$, and $\bR \subseteq \bP'$.  By
the maximality of $k$, the inequality $\bR \subseteq \bP'$ cannot be
obtained by using clause (3), and clauses (2) and (4) do not apply; so
it must have been obtained by using clause (1), or because $\bR$ is
equal to $\epsilon$, to the effect that $\bR \subseteq \bQ'$.

If some $\tA_i$ ($1\leq i\leq k$) is equal to ${F'}^*$, then the fact
that $\tP$ is reduced implies that $k$ cannot be larger than or equal
to $2$. Hence $k=1$. By induction hypothesis,
$\varphi (\tR) \leq \varphi (\tQ')$, and therefore
$\varphi (\tP) = \varphi ({F'}^*) \oplus \varphi (\tR) \leq \varphi
({F'}^*) \oplus \varphi (\tQ') = \varphi (\tP')$. Additionally, if
$\bP \neq \bP'$, then $\bR$ must be different from $\bQ'$, so
$\varphi (\tR) < \varphi (\tQ')$, and therefore
$\varphi (\tP) = \varphi ({F'}^*) \oplus \varphi (\tR) < \varphi
({F'}^*) \oplus \varphi (\tQ') = \varphi (\tP')$.

In the remaining case, every $\tA_i$ is different from ${F'}^*$. We
claim that $\varphi (\tA_i) < \varphi ({F'}^*)$. This follows from
Lemma~\ref{lemma:rk:stat}~(1) if $\tA_i$ is of the form $C^?$, and
from Lemma~\ref{lemma:rk:stat}~(2) if $\tA_i$ is of the form $F^*$.
Since $\varphi ({F'}^*) = \omega^{\adjsw{||F'||}}$ is
$\oplus$-indec\-omposable,
$\varphi (\tA_1) \oplus \cdots \oplus \varphi (\tA_k) < \varphi
({F'}^*)$, and therefore
$\varphi (\tP) = \varphi (\tA_1) \oplus \cdots \oplus \varphi (\tA_k)
+ \varphi (\tR) < \varphi ({F'}^*) \oplus \varphi (\tQ') = \varphi
(\tP')$. \qed

\begin{proposition}
  \label{prop:dim:X*:ubound}
  For every Noetherian space $X$,
  $\sob {X^*} \leq \omega^{\adjsw{||X||}}+1$, or equivalently,
  $\rsob {X^*} \leq \omega^{\adjsw{||X||}}$.
\end{proposition}
\proof As a corollary of Lemma~\ref{lemma:X*:phi}, given any two
word-products $\tP$ and $\tP'$ such that $\bP = \bP'$, we have
$\varphi (\tP) = \varphi (\tP')$.  Hence $\varphi$ defines a strictly
monotonic map from irreducible closed subsets $\bP$ of $X^*$ (not just
word-products $\tP$) to the class of ordinals.  Its largest value is
$\varphi (X^*)$, since $X^*$ itself is a word-product, so the image of
$\varphi$ lies entirely inside the ordinal
$\varphi (X^*)+1 = \omega^{\adjsw{||X||}}+1$.  \qed

We turn to the matching lower bound on $\sob {X^*}$.  The
specialization preordering $\leq$ of a space $X$ induces an
equivalence relation $\equiv$ by $x \equiv y$ if and only if
$x \leq y$ and $y \leq x$.  This allows us to partition $X$ into
equivalence classes.  If $X$ is $T_0$, then those equivalence classes
contain exactly one point, and the following lemma would simply say
that $F$ and $F'$ differ by exactly one point.
\begin{lemma}
  \label{lemma:+1}
  Let $F$, $F'$ be two closed subsets of a Noetherian space $X$, with
  $F \subseteq F'$ and $||F'|| = ||F||+1$.  Then $F' \diff F$ is a
  single equivalence class with respect to the specialization
  preordering of $X$.  Given any point $x$ in that class, $F' = F \cup
  \dc x$.
\end{lemma}
\proof Since $F \subseteq F'$ and $||F'|| \neq ||F||$, there is a
point in $F' \diff F$.  Let us imagine that $F' \diff F$ contains at
least two non-equivalent points $x$ and $y$.  Without loss of
generality, we may assume that $y \not\leq x$.  Then $F$ is strictly
included in the closed set $F \cup \dc x$, which is strictly included
in $F \cup \dc x \cup \dc y \subseteq F'$.  It follows that
$||F'|| \geq ||F \cup \dc x \cup \dc y|| \geq ||F \cup \dc x||+1 \geq
||F||+2$, which is impossible.  Finally, given any $x \in F' \diff F$,
the points of $F'$ consist of those of $F$, plus those that are
equivalent to $x$, hence are in $\dc x$.  Conversely, $x$ is in $F'$,
so $\dc x$ is included in $F'$, and therefore
$F \cup \dc x \subseteq F'$.  \qed

In the situation of Lemma~\ref{lemma:+1}, $F'$ is obtained by adding
an irreducible closed set $C$ to $F$, namely $C \eqdef \dc x$.
Using such sets $F$ and $C$, we build word-products of the form
$(F^* C^?)^m \tP$, where $m \in \nat$ and $\tP$ is a word-product
built from subsets of $F$ (equivalently, such that $\bP \subseteq F^*$).

\begin{lemma}
  \label{lemma:dim:X*:FC}
  Let $F$ be a closed subset of a Noetherian space $X$, $C$ be an
  irreducible closed subset of $X$ that is not included in $F$.  Let
  also $\tP$ and $\tQ$ be two word-products such that
  $\bP, \bQ \subseteq F^*$.  For all $m, n \in \nat$,
  $(F^* C^?)^m \bP \subseteq (F^* C^?)^n \bQ$ if and only if
  $(m, \bP)$ is lexicographically smaller than or equal to $(n, \bQ)$,
  namely if and only if $m<n$, or $m=n$ and $\bP \subseteq \bQ$.
\end{lemma}
\proof Let us assume $(F^* C^?)^m \bP \subseteq (F^* C^?)^n \bQ$. We
pick an element $x$ from $C \diff F$. The word $x^m$ consisting of $m$
copies of $x$ is in $(F^* C^?)^m \bP$ (noting that the empty word is
in the semantics of every word-product, in particular in $\bP$),
hence in $(F^* C^?)^n \bQ$. Since $\bQ \subseteq F^*$ and
$x \not\in F$, $x^m$ must be in $(F^* C^?)^n$, and that in turn
implies that $n \geq m$. If $m< n$, $(m, \bP)$ is lexicographically
smaller than $(n, \bQ)$, so we are left to show that
$(F^* C^?)^n \bP \subseteq (F^* C^?)^n \bQ$ implies
$\bP \subseteq \bQ$, for every $n \in \nat$. This is by induction on
$n$. This is clear if $n=0$. Otherwise, only clause~(3) or~(4) of
Lemma~\ref{lemma:6.10} can have led to this inclusion, and this can
only occur if $C^? (F^* C^?)^{n-1} \bP \subseteq (F^* C^?)^n \bQ$.
Since $C^? \not\subseteq F^*$, only clause (1) can have led to the
latter, so
$C^? (F^* C^?)^{n-1} \bP \subseteq C^? (F^* C^?)^{n-1} \bQ$. This can
only have been obtained through clause~(2), from the inclusion
$ (F^* C^?)^{n-1} \bP \subseteq (F^* C^?)^{n-1} \bQ$, and therefore
$\bP \subseteq \bQ$ by the induction hypothesis.

In the converse direction, we first claim that $m<n$ implies
$(F^* C^?)^m \bP \subseteq (F^* C^?)^n \bQ$.  Indeed, the words of
$(F^* C^?)^m \bP$ are the words $w_1 w_2$ where $w_1 \in (F^* C^?)^m$
and $w_2 \in \bP$.  Every such $w_2$ is in $F^*$.  This shows that
$(F^* C^?)^m \bP \subseteq (F^* C^?)^m F^*$.  Every word $w$ in
$(F^* C^?)^m F^*$ is also in $(F^* C^?)^n \bQ$, as the concatenation
of $w \in (F^* C^?)^m F^*$ and of
$\epsilon \in C^? (F^* C^?)^{n-m-1} \bQ$.

It remains to show that if $m=n$ and $\bP \subseteq \bQ$, then $(F^*
C^?)^m \bP \subseteq (F^* C^?)^n \bQ$, and that is obvious.
\qed

We will also need the following simpler construction.
\begin{lemma}
  \label{lemma:dim:X*:CC}
  Let $X$ be a Noetherian space.  For all $C, C', C'' \in \Sober X$
  such that $C', C'' \subseteq C$, for all $m, n \in \nat$,
  $(C^?)^m {C'}^? \subseteq (C^?)^n {C''}^?$ if and only if $(m, C')$
  is lexicographically smaller than or equal to $(n, C'')$.
\end{lemma}
\proof Let us assume $(C^?)^m {C'}^? \subseteq (C^?)^n {C''}^?$.  Let
$x$ be any point in $C'$, hence also in $C$.  Then $x^{m+1}$ is in the
left-hand side, hence also in the right-hand side, and this implies
that $n \geq m$.  If additionally $n=m$, then
$(C^?)^n {C'}^? \subseteq (C^?)^n {C''}^?$ can only be derived by $n$
applications of clause (2) of Lemma~\ref{lemma:6.10}, implying that
${C'}^? \subseteq {C''}^?$, namely $C' \subseteq C''$.

Conversely, if $m < n$, then
$(C^?)^m {C'}^? \subseteq (C^?)^n {C''}^?$, since every word in
$(C^?)^m {C'}^?$ consists of at most $m+1$ letters, all from $C$ (the
last one possibly being in $C'$, hence in $C$).  If $m=n$ and
$C' \subseteq C''$, then $(C^?)^m {C'}^? \subseteq (C^?)^n {C''}^?$,
too, as one sees easily.  \qed

We will use the latter lemmata in conjunction with the well-known fact
that, given two well-founded posets $P$ and $Q$, and points $p \in P$
and $q \in Q$,
$\rk_{P \times_\lex Q} (p, q) = |Q| \times \rk_P (p) + \rk_Q (q)$,
where $P \times_\lex Q$ denotes the lexicographic product of $P$ and
$Q$.  As a consequence, $|P \times_\lex Q| = |Q| \times |P|$.  It also
follows that every well-founded poset $R$ that admits a strictly
monotonic map from $P \times_\lex Q$ to $R$ must have rank at least
$|Q| \times |P|$.

\begin{proposition}
  \label{prop:dim:X*:lbound}
  For every non-empty Noetherian space $X$,
  $\sob {X^*} \geq \omega^{\adjsw{||X||}}+1$.
\end{proposition}
\proof By well-founded induction on non-empty closed subsets $F$ of
$X$, we show that
$\rk_{\Sober (X^*)} (F^*) \geq \omega^{\adjsw{||F||}}$. The claim will
follow since $\sob {X^*} = \rk_{\Sober (X^*)} (X^*)+1$, as $X^*$ is
the largest element of $\Sober (X^*)$, and using
Lemma~\ref{lemma:dimension}.

Since $F$ is non-empty, $||F|| \geq ||\emptyset||+1 = 1$.

If $||F||=1$, then using Lemma~\ref{lemma:+1} on the strict inclusion
$\emptyset \subsetneq F$, $F$ itself is an irreducible closed subset
of the form $\dc x$, and this contains no proper closed subset except
the empty set. The elements of $\Sober (X^*)$ included in $F^*$ are
$(F^?)^n$, $n \in \nat$, plus $F^*$; the rank of $(F^?)^n$ is $n$, so
$\rk_{\Sober (X^*)} (F^*) = \omega$.

If $||F||$ is a successor ordinal $\alpha+1$ different from $1$, then
by the inductive definition of rank (and Lemma~\ref{lemma:rk:C}),
there is a closed subset $F'$ of $F$ such that $||F'|| = \alpha$.  By
induction hypothesis,
$\rk_{\Sober (X^*)} ({F'}^*) \geq \omega^{\adjsw\alpha}$.  Hence there
are irreducible closed subsets $\bP$ of ${F'}^*$ of arbitrary rank
between $0$ and $\omega^{\adjsw\alpha}$ in $\Sober (X^*)$.  By
Lemma~\ref{lemma:+1}, $F$ is equal to $F' \cup C$, where $C$ is an
irreducible closed subset of the form $\dc x$ with $x \not\in F'$.
Using Lemma~\ref{lemma:dim:X*:FC}, there is a strictly monotonic map
from $\nat \times_\lex Q$ into the poset of proper irreducible closed
subsets of $F^*$, where $Q$ is the poset of irreducible closed subsets
of ${F'}^*$ of rank strictly less than $\omega^{\adjsw\alpha}$,
through $(n, \bP) \mapsto ({F'}^* C^?)^n \bP$.  It follows that
$\rk_{\Sober (X^*)} (F^*) \geq \omega^{\adjsw\alpha} \times \omega =
\omega^{\adjsw\alpha+1}$.  Now $\adjsw\alpha+1 = \adjsw{(\alpha+1)}$,
in both cases of the definition of $\adjsw\alpha$.

If $||F||$ is a limit ordinal $\alpha$, then there are closed subsets
$F'$ of $F$ of arbitrarily high rank $\beta < \alpha$, and by
induction
$\rk_{\Sober (X^*)} (F^*) \geq \rk_{\Sober (X^*)} ({F'}^*) \geq
\omega^{\adjsw\beta}$. When $\alpha$ is not critical (and since
$\alpha$ is a limit ordinal), we realize that $\adjsw\alpha = \alpha$,
and that the latter inequality implies
$\rk_{\Sober (X^*)} (F^*) \geq \omega^\beta$ for every
$\beta < \alpha$, and therefore
$\rk_{\Sober (X^*)} (F^*) \geq \omega^\alpha = \omega^{\adjsw\alpha}$,
by taking suprema over $\beta < \alpha$.

When $\alpha$ is a critical ordinal, we require another argument.  We
first observe that $\sob F \geq \alpha$.  Indeed, if
$\sob F < \alpha$, then using Proposition~\ref{prop:dim:stature}~(2)
and the fact that $\alpha$ is critical,
$||F||+1 \leq \omega^\alpha = \alpha$, which is impossible since
$||F|| = \alpha$.  Using Lemma~\ref{lemma:dimension}, there is an
irreducible closed subset $C$ of $F$ such that
$\rk_{\Sober X} (C) + 1 \geq \alpha$.  As a consequence,
$\rk_{\Sober X} (C) \geq \alpha$.

Using Lemma~\ref{lemma:dim:X*:CC}, there is a strictly monotonic map
from $\nat \times_\lex Q$ into the poset of proper irreducible closed
subsets of $F^*$, where $Q$ is the poset of irreducible closed subsets
of $C$ of rank strictly less than $\alpha$, through
$(n, C') \mapsto (C^?)^n {C'}^?$.  Therefore
$\rk_{\Sober (X^*)} (F^*) \geq \alpha \times \omega$.  Now
$\alpha \times \omega = \omega^\alpha \times \omega =
\omega^{\alpha+1} = \omega^{\adjsw\alpha}$.  \qed

We finally combine Proposition~\ref{prop:dim:X*:ubound} and
Proposition~\ref{prop:dim:X*:lbound}, and we include the case of the
empty Noetherian space in the following theorem.
\begin{theorem}
  \label{thm:dim:X*}
  For every Noetherian space $X$,
  $\sob {X^*} = \omega^{\adjsw{||X||}}+1$ (equivalently,
  $\rsob {X^*} = \omega^{\adjsw{||X||}}$) if $X$ is non-empty, $1$
  otherwise.
\end{theorem}
A special case of this result appears as Proposition~5.5 of
\cite{BFG-lmcs20}, where $X$ was assumed to be non-empty and finite.

\subsection{The stature of $X^*$}
\label{sec:stature-x}

We start with an easy upper bound.
\begin{lemma}
  \label{lemma:stature:X*:ubound:infinite}
  For every non-empty Noetherian space $X$,
  $||X^*|| \leq \omega^{\omega^{\adjsw\alpha}}$, where
  $\alpha \eqdef ||X||$.
\end{lemma}
\proof
$X^*$ is irreducible closed, so the space $X^*$ has exactly one
component.
By Proposition~\ref{prop:dim:stature}~(3), $||X^*|| \leq \omega^{\rsob
  {X^*}}$, and $\rsob {X^*} = \omega^{\adjsw\alpha}$ by
Theorem~\ref{thm:dim:X*}.
%
%
\qed

One can improve upon that upper bound when $||X||$ is finite.  We need
the following simple observation first.
\begin{lemma}
  \label{lemma:finite:conv}
  For every Noetherian space $X$ such that $\alpha \eqdef ||X||$ is
  finite, $X$ has only finitely many open subsets.  Additionally,
  $X$ has exactly $\alpha$ equivalence classes with respect to its
  specialization preordering.
\end{lemma}
In particular, every $T_0$, finite Noetherian space $X$ contains
exactly $||X||$ points.  This is a form of converse to
Lemma~\ref{lemma:finite}.

\proof We prove the second claim by induction on $\alpha$.  If
$\alpha=0$, then $X$ is empty, and the claim is clear.  Otherwise,
there is a closed subset $F$ of $X$ such that
$\rk_{\Hoare_0 X} (F) = \alpha-1$, while $\rk_{\Hoare_0 X} (X) = \alpha$
by definition.  By Lemma~\ref{lemma:+1} with $F' \eqdef X$, there is a
point $x$ such that $X = F \cup \dc x$, and $X \diff F$ is a single
equivalence class with respect to the specialization preordering
$\leq$ of $X$.  By Lemma~\ref{lemma:rk:C}, $||F|| = \alpha-1$, so the
subspace $F$ has exactly $\alpha-1$ equivalence classes with respect
to $\leq_{|F}$, hence $X$ has exactly $\alpha$ equivalence classes.

The first claim follows, since every open subset of $X$ is
upwards-closed with respect to $\leq$, hence is closed under the
associated equivalence relation, and therefore is a union of
equivalence classes.  \qed

\begin{lemma}
  \label{lemma:stature:X*:ubound:finite}
  For every Noetherian space $X$ such that $\alpha \eqdef ||X||$ is
  finite, we have:
  \begin{enumerate}
  \item if $\alpha=0$, then $||X^*||=1$;
  \item if $\alpha=1$, then $||X^*||=\omega$;
  \item if $\alpha \neq 0$, then $||X^*|| \leq \omega^{\omega^{\alpha-1}}$.
  \end{enumerate}
\end{lemma}
\proof
(1) If $\alpha=0$, then $X$ is empty, so $X^* = \{\epsilon\}$, and
therefore $||X^*|| = 1$.

(2) When $\alpha=1$, all the proper closed subsets of $X$ have rank
equal to $0$, hence are empty.  In other words, the topology of $X$
is the indiscrete topology, whose sole closed sets are $X$ and the
empty set.  Then $X$ is the sole element of $\Sober X$.
The only irreducible closed subsets of $X^*$ are $(X^?)^n$, $n \in
\nat$, plus $X^*$.  They form a chain, so any non-empty finite union
of word-products reduces to a single word-product.  This entails
that $\Hoare_0 (X^*)$ is the same chain, with the empty set added as a
new bottom element.  In particular, $||X^*|| = \omega$.

(3) We prove this by induction on $\alpha \geq 1$.
Item~(2) is the base case.  Let $\alpha \geq 2$.  By
Lemma~\ref{lemma:finite:conv}, $X$ has only finitely many closed
subsets.  We enumerate its \emph{proper} closed subsets as $F_1$,
\ldots, $F_m$, and its components as $C_1$, \ldots, $C_p$.  For every
$i$ with $1\leq i \leq m$,
$||F_i|| = \rk_{\Hoare_0 X} (F_i) < \rk_{\Hoare_0 X} (X) = ||X|| =
\alpha$ (using Lemma~\ref{lemma:rk:C}), so
$||F_i^*|| \leq \omega^{\omega^{\alpha-2}}$ for every $i$, by
induction hypothesis.

We build word-products $\tP_n$, $n \in \nat$, by induction on $n$, by
letting $\tP_0 \eqdef \epsilon$ and
$\tP_{n+1} \eqdef F_1^* \cdots F_m^* C_1^? \cdots C_p^? \tP_n$.  The
point of this construction is that, for every word-product $\tP$ such
that $\bP \neq X^*$, the inclusion $\bP \subseteq \bP_n$ holds for $n$
large enough, namely for every $n$ larger than or equal to the length
$|\tP|$ of $\tP$.  This is shown by induction on $|\tP|$.  If
$\tP = \epsilon$, then $\bP \subseteq \bP_n$ for every $n$, by
Lemma~\ref{lemma:6.10} (or rule (0)).  If $\tP$ is of the form
$\tA_1 \tQ$, then $|\tP| = 1 + |\tQ|$.  Let us fix $n \geq |\tP|$.  We
observe that $\bQ \neq X^*$, otherwise $\bP$ would contain every word
on $X$ and therefore be equal to $X^*$.  Hence we can apply the
induction hypothesis, so that $\bQ \subseteq \bP_{n-1}$.  We claim
that $\bP \subseteq \bP_n$.  If $\tA_1$ is of the form $C^?$, then the
words of $\bP$ are those of $\bQ$, which are in $\bP_{n-1}$, hence
also in $\bP_n$, plus those of the form $yw$ with $y \in C$ and
$w \in \bQ \subseteq \bP_{n-1}$; then $C$ is included in some $C_j$,
so $yw$ is in $\bP_n$.  If $\tA_1$ if of the form $F^*$, then $F$
cannot be equal to the whole of $X$, since in that case $\bP$ would
contain every word on $X$ and therefore be equal to $X^*$.  Hence
$F=F_i$ for some $i$, and then the inclusion
$\bP = F_i^* \bQ \subseteq \bP_n = X^? F_1^* \cdots F_m^* \bP_{n-1}$
is immediate.

In particular, every proper closed subset $\bA$ of $X^*$ is included
in $\bP_n$ for $n$ large enough. This follows from the previous claim
by writing $\bA$ as a finite union of word-products, all different
from $X^*$. If the first one is included in $\bP_{n_1}$, the second
one is included in $\bP_{n_2}$, \ldots, and the last one is included
in $\bP_{n_k}$, then their union $\bA$ is included in $\bP_n$, where
$n \eqdef \max (n_1, n_2, \cdots, n_k)$.

We now claim that $||\bP_n|| < \omega^{\omega^{\alpha-1}}$ for every
$n \in \nat$. We prove this by induction on $n$. The base case reduces
to $||\bP_0|| = 1 < \omega^{\omega^{\alpha-1}}$. In the inductive
case, we know that $||\bP_n|| < \omega^{\omega^{\alpha-1}}$, and we
aim to prove that $||\bP_{n+1}|| < \omega^{\omega^{\alpha-1}}$. The
function $j \colon X_\bot \to X^*$ that maps every $x \in X$ to the
one-letter word $x$, and $\bot$ to $\epsilon$ is continuous: for all
open subsets $U_1$, \ldots, $U_k$ of $X$,
$j^{-1} (X^* U_1 X^* \cdots X^* U_k X^*)$ is equal to $U_1$ if
$k = 1$, to $X_\bot$ if $k=0$, and is empty if $k \geq 2$; in any
case, this is open. Let
$f \colon F_1^* \times \cdots \times F_m^* \times X_\bot^p \times
\bP_n \to \bP_{n+1}$ map $(w_1, \cdots, w_m, x_1, \cdots, x_p, w)$ to
$w_1 \cdots w_m j (x_1) \cdots j (x_p) w$. Using the fact that the
concatenation map from $X^* \times X^*$ to $X^*$ is continuous
\cite[Lemma~B.1]{FGL:partI}, $f$ is continuous. It is also clearly
surjective, hence Skula dense by Lemma~\ref{lemma:surj}~(2). By
Lemma~\ref{lemma:stature:mono}, $||\bP_{n+1}||$ is less than or equal
to
$||F_1^* \times \cdots \times F_m^* \times X_\bot^p \times
\bP_n||$. By Theorem~\ref{thm:stature:prod}, the latter is equal to
the natural product of $||F_1^*||$, \ldots, $||F_m^*||$, $||X_\bot||$
$p$ times, and $||\bP_n||$. We have $||X_\bot|| = 1+||X|| = 1+\alpha$
by Proposition~\ref{prop:lifting}, we recall that
$||F_i^*|| \leq \omega^{\omega^{\alpha-2}}$ for every $i$, and that
$||\bP_n|| < \omega^{\omega^{\alpha-1}}$. In particular, all the terms
in the natural product are strictly smaller than
$\omega^{\omega^{\alpha-1}}$. (For the first one, we use the fact that
$\omega^{\omega^{\alpha-1}} \geq \omega$, while $1+\alpha$ is finite.)
Since $\omega^{\omega^{\alpha-1}}$ is $\otimes$-indecomposable,
$||\bP_{n+1}|| < \omega^{\omega^{\alpha-1}}$.

We put everything together.  The family ${(\bP_n)}_{n \in \nat}$ is a
cofinal family of closed subsets of $X^*$.  By
Proposition~\ref{prop:limit},
$||X^*|| \leq \sup_{n \in \nat} (||\bP_n||+1) \leq
\omega^{\omega^{\alpha-1}}$.  \qed

In the search for a lower bound of $||X^*||$, we will need the
following trick.  Given a well-founded poset $(P, \leq)$, let $<$ be
the strict part of $\leq$, and let us call \emph{step} of $P$ any pair
$(p, p^+)$, where $p, p^+ \in P$ and $p < p^+$.  A step $(p, p^+)$ is
an \emph{increment} if and only $\rk_P (p^+) = \rk_P (p)+1$.  We write
$\Step (P)$ for the set of steps of $P$, and $\Inc (P)$ for the subset
of increments of $P$.  Both steps and increments are ordered strictly
by $(p, p^+) < (q, q^+)$ if and only if $p^+ \leq q$.  We reuse the
same notation $<$; no confusion should arise.  We also write
$(p, p^+) \leq (q, q^+)$ if and only if $(p, p^+) < (q, q^+)$ or
$(p, p^+) = (q, q^+)$.
\begin{proposition}
  \label{prop:Inc}
  For every well-founded poset $P$,
  $|\Step (P)| = |\Inc (P)| = |P|-1$.
\end{proposition}
We recall that, when $\alpha$ is not a successor ordinal,
$\alpha-1=\alpha$ by convention.

\proof When $|P|=0$, the sets $P$, $\Step (P)$ and $\Inc (P)$ are
empty, so the claim is clear.  We now assume that $|P|$ is either a
successor or a limit ordinal.  We note that
$|\Inc (P)| \leq |\Step (P)|$, since $\Inc (P) \subseteq \Step (P)$.

The strictly monotonic map $(p, p^+) \mapsto \rk_P (p)$ shows that
$|\Step (P)| \leq |P|$.  When $|P|$ is a successor ordinal $\alpha+1$,
this map takes its values in $\alpha$, since for every
$(p, p^+) \in \Step (P)$, $\rk_P (p) < \rk_P (p^+) < \alpha+1$; in
that case, we can improve the former inequality to
$|\Step (P)| \leq \alpha$, namely $|\Step (P)| \leq |P|-1$.

With an eye toward showing the converse, we prove that
for every ordinal $\alpha$,
\begin{enumerate}[label=(\roman*)]
\item for every $q^+ \in P^\top$ such that
  $\rk_{P^\top} (q^+)=\alpha+1$, there is an element $q < q^+$ such
  that $\rk_{P^\top} (q) = \alpha$ and
  $\rk_{\Inc (P^\top)} (q, q^+) \geq \alpha$;
\item for every $q \in P^\top$ such that $\rk_{P^\top} (q) = \alpha$
  and $\alpha$ is a limit ordinal,
  $\sup \{\rk_{\Inc (P^\top)} (p, p^+) \mid (p, p^+) \in \Inc
  (P^\top), p^+ < q\} \geq \alpha$.
\end{enumerate}
This is by mutual induction on $\alpha$. We start with (i). In this
case, there is an element $q < q^+$ such that
$\rk_{P^\top} (q) = \alpha$, and therefore $(q, q^+)$ is a increment
of $P^\top$. If $\alpha$ is itself a successor ordinal $\beta+1$, then
there is an element $p < q$ such that $\rk_{P^\top} (p) = \beta$ and
$\rk_{\Inc (P^\top)} (p, q) \geq \beta$, by induction hypothesis. We
note that $(p, q) < (q, q^+)$, so
$\rk_{\Inc (P^\top)} (q, q^+) \geq \beta+1 = \alpha$. If $\alpha=0$,
then the claim that $\rk_{\Inc (P^\top)} (q, q^+) \geq \alpha$ is
obvious. If $\alpha$ is a limit ordinal, then for every
$\beta < \alpha$, there is an increment $(p, p^+)$ such that $p^+ < q$
such that $\rk_{\Inc (P^\top)} (p, p^+) \geq \beta$ by the induction
hypothesis, part~(ii). Since $(p, p^+) < (q, q^+)$,
$\rk_{\Inc (P^\top)} (q, q^+) \geq \beta+1$. Taking suprema over
$\beta < \alpha$, $\rk_{\Inc (P^\top)} (q, q^+) \geq \alpha$. We turn
to (ii). For every ordinal $\beta < \alpha$, $\beta+1$ is still
strictly smaller than $\alpha$. Then, there is an element $p^+ < q$
such that $\rk_{P^\top} (p^+) = \beta+1$. By induction hypothesis,
part~(i), there is an element $p < p^+$ such that
$\rk_{P^\top} (p) = \beta$ and
$\rk_{\Inc (P^\top)} (p, p^+) \geq \beta$. Hence
$\sup \{\rk_{\Inc (P^\top)} (p, p^+) \mid (p, p^+) \in \Inc (P^\top),
p^+ < q\} \geq \sup_{\beta < \alpha} \beta = \alpha$.

We now fix $\alpha \eqdef \rk_{P^\top} (\top)$, namely
$\alpha \eqdef |P|$.

If $\alpha$ is a successor ordinal, by (i) there is an increment
$(q, \top)$ such that $\rk_{\Inc (P^\top)} (q, \top) \geq \alpha-1$.
Therefore $|\Inc (P^\top)| \geq \alpha$.  The strictly monotonic map
which maps every element $(p, p^+)$ of $\Inc (P)$ to
$\rk_{\Inc (P)} (p, p^+)$ and the remaining elements $(p, \top)$ of
$\Inc (P^\top)$ to $|\Inc (P)|$ shows that
$|\Inc (P^\top)| \leq |\Inc (P)| + 1$, so $\alpha \leq |\Inc (P)|+1$.
Therefore $|\Inc (P)| \geq |P|-1$.

If $\alpha$ is a limit ordinal, by (ii) the supremum of the ranks
$\rk_{\Inc (P^\top)} (p, p^+)$ over all elements $(p, p^+)$ of $\Inc
(P^\top)$ such that $p^+ < \top$, namely over all elements of $\Inc
(P)$, is at least $\alpha$.  The same therefore holds of the supremum
of the ordinals $\rk_{\Inc (P^\top)} (p, p^+)+1$, so $|\Inc (P)| \geq
\alpha = |P| = |P|-1$.  \qed

We will only require steps here, but increments will be needed in our
study of multisets, in Section~\ref{sec:bf-finite-multisets}.

In the following, we will build closed subsets of $X^*$ by
concatenation.  For this, we need to observe that, for any two closed
subsets $\bA$, $\bB$ of $X^*$, where $X$ is Noetherian, $\bA\bB$ is
also closed.  In order to see this, we write $\bA$ as a finite union
$\bigcup_{i=1}^m \bP_i$ of word-products, and similarly $\bB$ as a
finite union $\bigcup_{j=1}^n \bQ_j$ of word-products.
\blue{We can do this because the closed sets are exactly the semantics
  of word-SREs, namely the finite unions of (semantics of)
  word-products.  Then}
we note that
$\bA\bB = \bigcup_{\substack{1\leq i \leq m\\1\leq j\leq n}} \bP_i
\bQ_j$, a finite union of word-products.  This rests on the easily
checked fact that concatenation distributes over union.

Our main gadget consists of closed sets of the form
$(F^* C^?)^{n+1} \bB \cup \bA C^? \bB^+ \cup (F^* C^?)^n F^*$, where
$F$ and $C$ are as in Lemma~\ref{lemma:dim:X*:FC}, $\bA$ is a closed
subset of $X^*$ and $(\bB, \bB^+)$ is a step of $\Hoare_0 (X^*)$. We
will need to compare them with respect to inclusion, and this will
boil down to comparing pairs of sets of one of the two forms
$(F^* C^?)^{n+1} \bB$ or $\bA C^? \bB^+$. This leads us to examine
four different cases, which are explored in the following lemma.
\begin{lemma}
  \label{lemma:stature:X*:gadget}
  Let $F$ be a closed subset of a Noetherian space $X$, $C$ be an
  irreducible closed subset of $X$ that is not included in $F$, and
  $n \in \nat$.  For all closed subsets $\bA$, $\bB$, $\bB'$, $\bP$,
  $\bQ$ of $X^*$:
  \begin{enumerate}
  \item $(F^* C^?)^n \bB' \subseteq (F^* C^?)^n \bB$ if and only if
    $\bB' \subseteq \bB$;
  \item if $(F^* C^?)^{n+1} \bQ \subseteq \bA C^? \bB'$,
    $\bQ \neq \emptyset$, and $\bB' \subseteq F^*$, then
    $(F^* C^?)^n F^* \subseteq \bA$;
  \item if $\bP$ is non-empty, the conditions
    $\bP C^? \bB' \subseteq (F^* C^?)^{n+1} \bB$, $\bB \subseteq F^*$,
    and [if $n \geq 1$ then $\bP \not\subseteq (F^* C^?)^{n-1} F^*$],
    entail $\bB' \subseteq \bB$;
  \item if $\bP C^? \bB' \subseteq \bA C^? \bB'$,
    $\bB' \neq \emptyset$, and $\bB' \subseteq F^*$, then
    $\bP \subseteq \bA$.
  \end{enumerate}
\end{lemma}
\proof Let $x \in C \diff F$.

(1) The if direction is obvious.  Conversely, let us assume that
$(F^* C^?)^n \bB' \subseteq (F^* C^?)^n \bB$.  For every $w \in \bB'$,
$x^n w$ is in $(F^* C^?)^n \bB'$, hence in $(F^* C^?)^n \bB$.  Since
$x$ belongs to $C$ but not to $F$, $x^m w$ must be in $\bB$ for some
$m \leq n$.  Now $w \leq^* x^m w$, and we recall that $\leq^*$ is the
specialization ordering of $X^*$.  Every closed set is
downwards-closed with respect to the specialization ordering, so $w$
is in $\bB$.

(2) For every $w \in (F^* C^?)^n F^*$, we can write $w$ as
$w_0 c_1 w_1 \cdots w_{n-1} c_n w_n$, where each $w_i$ is in $F^*$ and
each $c_i$ is either in $C$ or is the empty word.  For each $i$, let
$c'_i$ be $x$ if $c_i = \epsilon$, $c_i$ otherwise.  Let us form the
word $w' \eqdef w_0 c'_1 w_1 \cdots w_{n-1} c'_n w_n$.  Since $\bQ$ is
non-empty and (downwards-)closed, $\epsilon$ is in $\bQ$.  Therefore
$w' x \epsilon = w'x$ is in $(F^* C^?)^{n+1} \bQ$.  The assumption
implies that it is in $\bA C^? \bB'$.  Since $\bB' \subseteq F^*$, and
therefore the final $x$ of $w'x$ can only be in $C^?$, not in $\bB'$,
$w'$ is in $\bA$ or $w'x$ is in $\bA$.  Now $w \leq^* w'$ and
$w \leq^* w'x$, and $\bA$ is (downwards-)closed, so $w$ is in $\bA$.

(3) Under the given assumptions, there is a word $w$ in $\bP$, and if
$n \geq 1$, we may assume that $w \not\in (F^* C^?)^{n-1} F^*$.  We
reason by contradiction and we assume that there is also a word $w'$ in
$\bB' \diff \bB$.  Then $wxw'$ is in $\bP C^?  \bB'$, hence in
$(F^* C^?)^{n+1} \bB$.  Since $x$ is in $C \diff F$, and since
$\bB \subseteq F^*$, $w$ must be in $(F^* C^?)^i F^*$ and $w'$ must be
in $(F^* C^?)^{n-i} \bB$ for some $i$, $0\leq i\leq n$.  If $n=i$,
then the latter would imply $w' \in \bB$, which is impossible.
Therefore $n-i \geq 1$.  In particular, $n \geq 1$, and $i \leq n-1$.
Since $w$ is in $(F^* C^?)^i F^*$, and since $i \leq n-1$, $w$ is in
the larger set $(F^* C^?)^{n-1} F^*$, which is impossible.

(4) Let us fix a word $w' \in \bB'$.  For every $w \in \bP$, $wxw'$ is
in $\bP C^? \bB'$, hence in $\bA C^? \bB'$.  Therefore we can write
$wxw'$ as $w_1 w_2$ where $w_1 \in \bA C^?$ and $w_2 \in \bB'$.
Since $\bB' \subseteq F^*$, $w_2$ cannot contain $x$.  Therefore $wx$
is a prefix of $w_1$, in particular $wx \leq^* w_1$, so that $wx$ is
in $\bA C^?$.  Then $w$ is in $\bA$ or $wx$ is in $\bA$, and in any
case $w$ is in $\bA$.  \qed

\begin{lemma}
  \label{lemma:stature:X*:inc:FC}
  Let $F$ be a closed subset of a Noetherian space $X$, and $C$ be an
  irreducible closed subset of $X$ that is not included in $F$.  Let
  $\bC_{n+1} \eqdef (F^* C^?)^n F^*$ for every $n \in \nat$, and
  $\bC_0 \eqdef \emptyset$.  Let also $\mathcal A_n$ be the set of all
  closed subsets of $X^*$ containing $\bC_n$ and strictly included in
  $\bC_{n+1}$.
  
  
  For every $n \in \nat$, the map
  $((\bB, \bB^+), \bA) \mapsto (F^* C^?)^{n+1} \bB \cup \bA C^? \bB^+
  \cup \bC_{n+1}$ is a strictly monotonic map from
  $\Step (\Hoare_0 (F^*)) \times_\lex \mathcal A_n$ to
  $\mathcal A_{n+1}$.
\end{lemma}
\proof
%
%
We first check that for all $(\bB, \bB^+) \in \Step (\Hoare_0 (F^*))$
and $\bA \in \mathcal A_n$,
$(F^* C^?)^{n+1} \bB \cup \bA C^? \bB^+ \cup \bC_{n+1}$ is in
$\mathcal A_{n+1}$.  That set is clearly closed and contains
$\bC_{n+1}$.  Since $\bB, \bB^+ \subseteq F^*$ and
$\bA \subseteq \bC_{n+1}$, it follows that
$(F^* C^?)^{n+1} \bB \cup \bA C^? \bB^+ \cup \bC_{n+1}$ is included in
$(F^* C^?)^{n+1} F^* \cup \bC_{n+1} C^? F^* \cup \bC_{n+1} = \bC_{n+2}
\cup \bC_{n+2} \cup \bC_{n+1} = \bC_{n+2}$.  If it were the whole of
$\bC_{n+2}$, then $\bC_{n+2} = (F^* C^?)^{n+1} F^*$ would be included
in $(F^* C^?)^{n+1} \bB \cup \bA C^? \bB^+ \cup \bC_{n+1}$, hence in
$(F^* C^?)^{n+1} \bB$ or in $\bA C^? \bB^+$ or in $\bC_{n+1}$, since
$\bC_{n+2}$ is irreducible (a word-product).  If
$(F^* C^?)^{n+1} F^* \subseteq (F^* C^?)^{n+1} \bB$, then
$F^* \subseteq \bB$ by Lemma~\ref{lemma:stature:X*:gadget}~(1); that
is impossible, since $\bB \subsetneq \bB^+ \subseteq F^*$.  If
$(F^* C^?)^{n+1} F^* \subseteq \bA C^? \bB^+$, then
$(F^* C^?)^n F^* \subseteq \bA$ by
Lemma~\ref{lemma:stature:X*:gadget}~(2); that is impossible because
$\bA$ is strictly included in $\bC_{n+1}$.  Finally, the inclusion
$(F^* C^?)^{n+1} F^* \subseteq \bC_{n+1} = (F^* C^?)^n F^*$ is also
impossible, since, given any $x \in C \diff F$, $x^{n+1}$ is in the
former but not in the latter.

Let $(\bB, \bB^+), (\bB', {\bB'}^+) \in \Step (\Hoare_0 (F^*))$,
$\bA, \bA' \in \mathcal A_n$, and let us assume that
$((\bB, \bB^+), \bA)$ is lexicographically smaller than
$((\bB', {\bB'}^+), \bA')$.  We verify that
$(F^* C^?)^{n+1} \bB \cup \bA C^? \bB^+ \cup \bC_{n+1}$ is a proper
subset of
$(F^* C^?)^{n+1} \bB' \cup \bA' C^? {\bB'}^+ \cup \bC_{n+1}$.

If $(\bB, \bB^+) < (\bB', {\bB'}^+)$, namely if
$\bB^+ \subseteq \bB'$, then both $(F^* C^?)^{n+1} \bB$ and
$\bA C^? \bB^+$ are included in $(F^* C^?)^{n+1} \bB'$.  This is
obvious for the first one, once we note that $\bB \subseteq \bB'$.
For the second one, we use the inclusions
$\bA \subseteq \bC_{n+1} = (F^* C^?)^n F^*$ and $\bB^+ \subseteq \bB'$
to deduce
$\bA C^? \bB^+ \subseteq (F^* C^?)^n F^* C^? \bB' = (F^* C^?)^{n+1}
\bB'$.  Hence
$(F^* C^?)^{n+1} \bB \cup \bA C^? \bB^+ \cup \bC_{n+1} \subseteq (F^*
C^?)^{n+1} \bB' \cup \bA' C^? {\bB'}^+ \cup \bC_{n+1} $.  If that
inclusion were an equality, then $(F^* C^?)^{n+1} \bB'$ would be
included in $(F^* C^?)^{n+1} \bB \cup \bA C^? \bB^+ \cup \bC_{n+1}$.
We write $\bB'$ as a finite union $\bigcup_{k=1}^p \bQ_k$ of
irreducible closed subsets of $X^*$.  For every $k$,
$(F^* C^?)^{n+1} \bQ_k$ is included in $(F^* C^?)^{n+1} \bB$ or in
$\bA C^? \bB^+$ or in $\bC_{n+1}$, by irreducibility.  If
$(F^* C^?)^{n+1} \bQ_k \subseteq \bA C^? \bB^+$, then
$(F^* C^?)^n F^* \subseteq \bA$ by
Lemma~\ref{lemma:stature:X*:gadget}~(2); that is impossible since
$\bA \subsetneq \bC_{n+1} = (F^* C^?)^n F^*$.  The case
$(F^* C^?)^{n+1} \bQ_k \subseteq \bC_{n+1} = (F^* C^?)^n F^*$ is also
impossible since, given any $x \in C \diff F$, $x^{n+1}$ is in
$(F^* C^?)^{n+1} \bQ_k$ but not in $(F^* C^?)^n F^*$.  Therefore, for
every $k$, $(F^* C^?)^{n+1} \bQ_k$ is included in
$(F^* C^?)^{n+1} \bB$; by Lemma~\ref{lemma:stature:X*:gadget}~(1) once
again, $\bQ_k$ is included in $\bB$.  Since that holds for every $k$,
$\bB'$ is included in $\bB$, which is impossible since
$\bB \subsetneq \bB^+ \subseteq \bB'$.

If $(\bB, \bB^+) = (\bB', {\bB'}^+)$ and $\bA \subsetneq \bA'$, we
need to show that
$(F^* C^?)^{n+1} \bB \cup \bA C^? \bB^+ \cup \bC_{n+1}$ is a proper
subset of $(F^* C^?)^{n+1} \bB \cup \bA' C^? \bB^+ \cup \bC_{n+1}$.
It is clearly a subset.  If the two sets were equal, then
$\bA' C^?  \bB^+$ would be included in
$(F^* C^?)^{n+1} \bB \cup \bA C^? \bB^+ \cup \bC_{n+1}$.  Let us write
$\bA'$ as a finite union of irreducible closed sets.  Amongst the
latter, let $\bP_1$, \ldots, $\bP_m$ be those that are not included in
$\bC_n$.  Since $\bA'$ contains $\bC_n$, $\bA'$ is equal to
$\bC_n \cup \bigcup_{i=1}^m \bP_i$.  For each $i$, $\bP_i C^? \bB^+$
is included in $(F^* C^?)^{n+1} \bB$ or in $\bA C^? \bB^+$ or in
$\bC_{n+1}$, by irreducibility.  If
$\bP_i C^? \bB^+ \subseteq (F^* C^?)^{n+1} \bB$, then we recall that
$\bB \subseteq F^*$, that $\bP_i$ is not included in
$\bC_n = (F^* C^?)^{n-1} F^*$ (if $n \geq 1$), and that $\bB^+$ is not
included in $\bB$; but those statements are contradictory, by
Lemma~\ref{lemma:stature:X*:gadget}~(3).  The assumption
$\bP_i C^? \bB^+ \subseteq \bC_{n+1} = (F^* C^?)^n F^*$ is
contradictory, too.  Indeed, given any $x \in C \diff F$, that would
imply that for every $w \in \bP_i$, $wx$ is in $(F^* C^?)^n F^*$,
hence contains at most $n$ occurrences of $x$; then $n \geq 1$ and $w$
contains at most $n-1$ occurrences of $x$, and since $w$ is arbitrary,
$\bP_i$ would be included in $(F^* C^?)^{n-1} F^* = \bC_n$.  Hence,
for each $i$, only the last assumption remains, namely that
$\bP_i C^? \bB^+$ is included in $\bA C^? \bB^+$.  By
Lemma~\ref{lemma:stature:X*:gadget}~(4), which applies since $\bB^+$
is non-empty, being a proper superset of $\bB$, $\bP_i$ is included in
$\bA$.  Since that holds for every $i$, and since $\bC_n$ is also
included in $\bA$, $\bA' = \bC_n \cup \bigcup_{i=1}^m \bP_i$ is
included in $\bA$.  This contradicts our assumption that
$\bA \subsetneq \bA'$.  \qed

\begin{corollary}
  \label{corl:stature:X*:lbound:F}
  Let $F$ be a closed subset of a Noetherian space $X$, and $C$ be an
  irreducible closed subset of $X$ that is not included in $F$.  Let
  $\alpha$ be any ordinal such that
  $||F^*|| \geq \omega^{\omega^\alpha}$.  Then
  $||(F \cup C)^*|| \geq \omega^{\omega^{\alpha+1}}$.
\end{corollary}
\proof
We have $|\Step (\Hoare_0 (F^*)| = |\Hoare_0 (F^*)|-1 = ||F^*||$ by
Proposition~\ref{prop:Inc} and the definition of stature.

Using the notations of Lemma~\ref{lemma:stature:X*:inc:FC},
$|\mathcal A_0| = ||\bC_1|| = ||F^*||$.  Then by part (3) of the
Lemma,
$|\mathcal A_1| \geq |\mathcal A_0| \times |\Step (\Hoare_0 (F^*))| =
||F^*||^2$.
By an easy induction on $n$, $|\mathcal A_n| \geq ||F^*||^{n+1}$.
Since $||F^*|| \geq \omega^{\omega^\alpha}$,
$|\mathcal A_n| \geq \omega^{\omega^\alpha \times (n+1)}$.

Now every element of $\mathcal A_n$ is clearly 
included in $(F \cup C)^*$.
It follows that
$||(F \cup C)^*|| \geq \omega^{\omega^\alpha \times (n+1)}$ for every
$n \in \nat$.  The supremum of the right-hand side as $n$ varies is
$\omega^{\omega^\alpha \times \omega} = \omega^{\omega^{\alpha+1}}$,
which allows us to conclude.  \qed

Using a notation of Schmidt's
\cite[Theorem~9]{Schmidt:maxot}, we define $\alpha'$ as $\alpha-1$ if
$\alpha$ is finite and non-zero, and as $\adjsw\alpha$ otherwise.
\begin{proposition}
  \label{prop:stature:X*:lbound}
  For every non-empty Noetherian space $X$,
  $||X^*|| \geq \omega^{\omega^{\alpha'}}$, where
  $\alpha \eqdef ||X||$.
\end{proposition}
\proof
By induction on the non-zero ordinal $\alpha$.  When
$\alpha=1$, $||X^*|| = \omega = \omega^{\omega^{\alpha-1}}$ by
Lemma~\ref{lemma:stature:X*:ubound:finite}~(2).

When $\alpha$ is a successor ordinal, there is a closed subset $F$ of
$X$ such that $\rk_{\Hoare_0 X} (F) = \alpha-1$.  By
Lemma~\ref{lemma:rk:C}, $||F|| = \alpha-1$.  We use the induction
hypothesis, and we obtain that
$||F^*|| \geq \omega^{\omega^{(\alpha-1)'}}$.  By
Lemma~\ref{lemma:+1}, there is a point $x$ such that
$X = F \cup \dc x$.  Let $C \eqdef \dc x$.  Applying
Corollary~\ref{corl:stature:X*:lbound:F}, we obtain that
$||X^*|| = ||(F \cup C)^*|| \geq \omega^{\omega^{(\alpha-1)'+1}}$.  A
simple case analysis on the definition of $\alpha'$ shows that
$(\alpha-1)'+1 = \alpha'$.

When $\alpha$ is a limit ordinal, for every $\beta < \alpha$, there is
a closed subset $F$ of $X$ such that
$||F|| = \rk_{\Hoare _0 X} (F) = \beta$, and by induction hypothesis
$||F^*|| \geq \omega^{\omega^{\beta'}}$.  It follows that
$||X^*|| \geq \sup_{\beta < \alpha} \omega^{\omega^{\beta'}} =
\omega^{\omega^{\sup_{\beta < \alpha} \beta'}}$.  Since $\beta$ and
$\beta'$ differ by at most $1$, and $\alpha$ is a limit ordinal,
$\sup_{\beta < \alpha} \beta' = \alpha$.  This shows that
$||X^*|| \geq \omega^{\omega^\alpha}$.  In particular,
$||X^*|| \geq \omega^{\omega^{\alpha'}}$, unless $\alpha$ is a
critical ordinal.  (We recall that $\alpha$ is a limit ordinal, and is
therefore of the form $\epsilon+n$ with $\epsilon$ critical and
$n \in \nat$ only if $n=0$.)

We finally deal with the case where $\alpha$ is a critical ordinal.
Let $\Sober' X$ be the subset of those elements of $\Sober X$ of
non-maximal rank, namely of rank different from $\rsob X = \sob X-1$
(see Lemma~\ref{lemma:dimension}).

For every $n \in \nat$, we consider \emph{$n$-elementary}
word-products, of the form $C_0^? C_1^? \cdots C_n^?$, where each
$C_i$ is an irreducible closed subset of $X$, and not all are of
maximal rank $\rsob X$.  Let $\mathcal E_n$ be the collection of
closed subsets of $X^*$ that one obtains as finite unions of
$n$-elementary word-products.

We note that every element $\bA$ of $\mathcal E_n$ is a subset of the
set $\bB_n$ of words of length at most $n+1$, and we claim that
this inclusion is proper.  Indeed, otherwise $\bB_n$ would be
included in $\bA$.  By Lemma~\ref{lemma:dimension}, one of the
components $C$ of $X$ has maximal rank $\rsob X$.  Then $(C^?)^{n+1}$
would be included in $\bA$.  Since $(C^?)^{n+1}$ is a word-product
hence irreducible, it would be included in one of the $n$-elementary
word-products $C_0^? C_1^? \cdots C_n^?$ whose union comprises $\bA$.
In turn, using clauses (1) and (2) of Lemma~\ref{lemma:6.10}, this
would force $C \subseteq C_i$ for each $i$, contradicting the fact
that at least one $C_i$ is of non-maximal rank.

The set $\bB_n$ is closed, because it is equal to the finite
union of the word-products $C_0^? C_1^? \cdots C_n^?$, where $C_0$,
$C_1$, \ldots, $C_n$ range over the components of $X$.  Let
$f_n \colon \Step (\Sober X) \times \mathcal E_n \to \mathcal E_{n+1}$
map $((C, C^+), \bA)$ to $\bB_n C^? \cup \bA {C^+}^?$.  Writing
$\bB_n$ as a finite union of word-products $\bP_i$,
$\bB_n C^?$ is the union of the word-products $\bP_i C^?$, in
which $C$ has non-maximal rank, being a proper subset of $C^+$.  Using
this, it is easy to check that $f_n ((C, C^+), \bA)$ is indeed in
$\mathcal E_{n+1}$.

We claim that $f_n$ is strictly monotonic from
$\Step (\Sober X) \times_\lex \mathcal E_n$ to $\mathcal E_{n+1}$.

If $(C, C^+) < (C', {C'}^+)$, namely if $C^+ \subseteq C'$, then both
$\bB_n C^?$ and $\bA {C^+}^?$ are included in
$\bB_n {C'}^?$, for every $\bA \in \mathcal E_n$.  Let
$w \in \bB_n \diff \bA$ and $x \in C^+ \diff C$.  By
concatenating $w$ with enough copies of $x$, we can assume that $w$
has length exactly $n+1$; the result is still not in $\bA$, since
$\bA$ is downwards-closed with respect to $\leq^*$.  Then $wx$ is in
$\bB_n {C'}^?$, but neither in $\bB_n C^?$ (since
$x \not\in C$) nor in $\bA {C^+}^?$ (since $w \not\in \bA$).
Therefore the inclusion is strict.

If $(C, C^+) = (C', {C'}^+)$ and $\bA \subsetneq \bA'$, then
$\bA {C^+}^? \subseteq \bA' {{C'}^+}^?$, so
$f_n ((C, C^+), \bA) \subseteq f_n ((C', {C'}^+), \bA')$. Let
$w \in \bA' \diff \bA$ and $x \in C^+ \diff C$. Since $w$ is in an
$n$-elementary word-product $C_0^? C_1^? \cdots C_n^?$ included in
$\bA'$, it is of the form $w_0 w_1 \cdots w_n$ where each $w_i$ is in
$C_i^?$; for those indices $i$ such that $w_i=\epsilon$, we may replace
$w_i$ by some letter from $C_i$, and we will obtain a word that is
still in $C_0^? C_1^? \cdots C_n^?$, hence in $\bA'$, but not in $\bA$
since $\bA$ is downwards-closed with respect to $\leq^*$. Hence,
without loss of generality, we may assume that $w$ has length exactly
$n+1$. Then $wx$ is in $f_n((C', {C'}^+), \bA')$ and neither in
$\bB_n C^?$ (since $x \not\in C$ and $w$ is of length $n+1$) nor
in $\bA {C^+}^?$ (since $w \not\in \bA$, $w$ is of length $n+1$, and
$\bA$ contains words of lengths at most $n+1$ only), hence not in
$f_n ((C, C^+), \bA)$. Hence the inclusion is strict.

It follows that
$|\mathcal E_{n+1}| \geq |\mathcal E_n| \times |\Step (\Sober X)|$,
for every $n \in \nat$. Now $\sob X^* = \omega^{||X||+1}+1$ by
Theorem~\ref{thm:dim:X*}, namely, $|\Sober X| =
\omega^{\alpha+1}+1$. Using Proposition~\ref{prop:Inc},
$|\Step (\Sober X)| = \omega^{\alpha+1}$. Hence
$|\mathcal E_{n+1}| \geq |\mathcal E_n| \times \omega^{\alpha+1}$, for
every $n \in \nat$. Since $\mathcal E_0$ is non-empty,
$|\mathcal E_0| \geq 1$. By induction on $n$, then,
$|\mathcal E_n| \geq (\omega^{\alpha+1})^n$. In particular,
$|\mathcal E_n| \geq (\omega^\alpha)^n = \alpha^n$, since $\alpha$ is
critical. Since each element of $\mathcal E_n$ is a proper closed
subset of $X^*$, it follows that
$||X^*|| \geq \sup_{n \in \nat} \alpha^n = \alpha^\omega$. But
$\alpha^\omega = (\omega^\alpha)^\omega = \omega^{\alpha \times
  \omega}$ and
$\alpha \times \omega = \omega^\alpha \times \omega =
\omega^{\alpha+1}$, so
$||X^*|| \geq \omega^{\omega^{\alpha+1}} = \omega^{\omega^{\alpha'}}$.
\qed

We put together the results of
Proposition~\ref{prop:stature:X*:lbound},
Lemma~\ref{lemma:stature:X*:ubound:infinite}, and
Lemma~\ref{lemma:stature:X*:ubound:finite}.
\begin{theorem}
  \label{thm:stature:X*}
  For every Noetherian space $X$, letting $\alpha \eqdef ||X||$, the
  stature $||X^*||$ is equal to $1$ if $\alpha=0$, and to
  $\omega^{\omega^{\alpha'}}$ otherwise, where $\alpha' = \alpha-1$ if
  $\alpha$ is finite and non-zero, $\alpha' = \alpha+1$ if
  $\alpha=\epsilon+n$ for some critical ordinal $\epsilon$ and some
  natural number $n$, and $\alpha' = \alpha$ otherwise.
\end{theorem}

\section{Heterogeneous Words and the Prefix Topology}
\label{sec:bf-heter-words}

The prefix topology is another interesting topology on $X^*$
\cite[Exercise~9.7.36]{JGL:topo}, which has no equivalent in the wqo
world.  Its specialization quasi-ordering is a form of the prefix
ordering, which is almost never a wpo.  We would like to mention that
this is the topology needed to decide reachability of sets defined by
forbidden patterns in the so-called oblivious $k$-stack system model
of \cite[Section~5]{JGL-icalp10}.

In general, the prefix topology makes sense not just on $X^*$, but on
spaces of sequences of elements taken from possibly different spaces,
see \cite[Section~9]{FGL:partI}.  Henceforth, let $X_1$, $X_2$,
\ldots, $X_n$, \ldots{} be countably many topological spaces.  A
\emph{heterogeneous word} over these spaces is any tuple
$(x_1, x_2, \cdots, x_m)$ in
$X_1 \times X_2 \times \cdots \times X_m$, $m \in \nat$.  We write it
as $x_1 x_2 \cdots x_m$, and we call $m = |w|$ the \emph{length} of
the word $w = x_1 x_2 \cdots x_m$.

We write $\bigtriangleright_{n=1}^{+\infty} X_n$ for the space of all
such heterogeneous words, with the so-called \emph{prefix topology}.
The latter is defined as follows.  A \emph{telescope} on
${(X_n)}_{n \geq 1}$ is a sequence
$\mathcal U \eqdef U_0, U_1, \cdots, U_n, \cdots$ of open sets, where
$U_n$ is open in $\prod_{i=1}^n X_i$ for each $n \in \nat$, and such
that $U_n X_{n+1} \subseteq U_{n+1}$ for every $n \in \nat$.  (We
write $UX$ instead of $U \times X$, for convenience, and $\epsilon$
for the empty tuple.  When $n=0$, $\prod_{i=1}^n X_i = \{\epsilon\}$,
and $U_0$ can only be the empty set or $\{\epsilon\}$.)  A \emph{wide
  telescope} is a telescope such that $U_n = \prod_{i=1}^n X_i$ for
some $n \in \nat$, or equivalently for all $n$ large enough.  The open
sets of the prefix topology are the empty set and all the sets
$\lfloor \mathcal U \rangle \eqdef \{w \in
\bigtriangleright_{n=1}^{+\infty} X_n \mid w \in U_{|w|}\}$, where
$\mathcal U$ is a wide telescope.

The specialization preordering of
$\bigtriangleright_{n=1}^{+\infty} X_n$ is the \emph{prefix
  preordering} $\leq^\triangleright$, defined by
$x_1 \cdots x_m \leq^\triangleright y_1 \cdots y_n$ if and only if
$m \leq n$ and $x_i \leq y_i$ for every $i$, $1\leq i\leq m$.
Moreover, when all the spaces $X_n$ are Noetherian, so is
$\bigtriangleright_{n=1}^{+\infty} X_n$
\cite[Proposition~9.1]{FGL:partI}.

The irreducible closed subsets of
$\bigtriangleright_{n=1}^{+\infty} X_n$ are exactly the sets of the
form
$\lceil C_1 \cdots C_n\rangle \eqdef \{x_1 \cdots x_m \mid m \leq n,
x_1 \in C_1, \ldots, x_m \in C_m\}$, where $n \in \nat$ and each $C_i$
is irreducible closed in $X_i$ for each $i$, plus
$\bigtriangleright_{n=1}^{+\infty} X_n$ itself when all the spaces
$X_n$ are non-empty \cite[Lemma~9.6]{FGL:partI}.  The set
$\bigtriangleright_{n=1}^{+\infty} X_n$ contains all the other
irreducible subsets, and the inclusion relation between the sets of
the former kind is given by
$\lceil C_1 \cdots C_m\rangle \subseteq \lceil C'_1 \cdots
C'_n\rangle$ if and only if $m \leq n$ and $C_i \subseteq C'_i$ for
every $i$, $1\leq i\leq m$, namely if and only if
$C_1 \cdots C_m \subseteq^\triangleright C'_1 \cdots C'_n$; in
general, if no $X_n$ is empty, the map
$i \colon (\bigtriangleright_{n=1}^{+\infty} \Sober {X_n})^\top \to
\Sober (\bigtriangleright_{n=1}^{+\infty} X_n)$ that sends $\top$ to
$\bigtriangleright_{n=1}^{+\infty} X_n$ and $C_1 \cdots C_m$ to
$\lceil C_1 \cdots C_m \rangle$ is an order-isomorphism
\cite[Proposition~9.7]{FGL:partI}.

Computing the sobrification rank of $\bigtriangleright_{n=1}^{+\infty} X_n$
therefore boils down to computing the ordinal rank of a product of
well-founded posets in the prefix ordering.

For posets $P_1$, \ldots, $P_n$, \ldots{}, we define
$\bigtriangleright_{n=1}^{+\infty} P_n$ as the set of words
$p_1 \cdots p_n$ where each $p_i$ is in $P_i$, ordered by the prefix
ordering $\leq^\triangleright$.  This is also the poset underlying the
space $\bigtriangleright_{n=1}^{+\infty} P_n$, where each $P_n$ is
given the Alexandroff topology of its ordering.  We also let
$\bigtriangleright_{i=1}^n P_i$ be the subposet of words of length at
most $n$.

In the following, let $\leq_P$ denote the ordering on a poset $P$.
\begin{lemma}
  \label{lemma:triangleright}
  For any two posets $P$, $Q$, let $P \triangleright Q$ be the set of
  elements that are either elements of $P$ or pairs in $P \times Q$,
  ordered by: $(p, q) \leq_{P \triangleright Q} (p', q')$ if and only
  if $p\leq_P p'$ and $q \leq_Q q'$; $p \leq_{P \triangleright Q} p'$
  if and only if $p \leq_P p'$; $p \leq_{P \triangleright Q} (p', q')$
  if and only if $p \leq_P p'$;
  $(p, q) \not\leq_{P \triangleright Q} p'$ for any $p$, $p'$, $q$.

  If $P$ and $Q$ are well-founded and non-empty, then for all $p \in
  P$ and $q \in Q$,
  \begin{enumerate}
  \item $\rk_{P \triangleright Q} (p) = \rk_P (p)$;
  \item $\rk_{P \triangleright Q} (p, q) = \rk_P (p) \oplus (1+\rk_Q
    (q))$.
  \end{enumerate}
\end{lemma}
\proof Let $Q_\bot$ be $Q$ with a fresh element $\bot$ added below all
others.  For every $q \in Q$, $\rk_{Q_\bot} (q) = 1+\rk_Q (q)$, by an
easy induction on $\rk_Q (q)$.  The map
$f \colon P \triangleright Q \to P \times_\lex Q_\bot$ defined by
$f (p, q) \eqdef (p, q)$, $f (p) \eqdef (p, \bot)$ is an
order-isomorphism.  The claim then follows from
Lemma~\ref{lemma:rk:x}.  \qed

\begin{lemma}
  \label{lemma:triangle:rank:finite}
  Let $P_1$, \ldots, $P_n$, \ldots{} be non-empty well-founded posets.
  For every word $w \eqdef p_1 \cdots p_m$ in
  $P \eqdef \bigtriangleright_{n=1}^{+\infty} P_n$ (resp., in
  $P \eqdef \bigtriangleright_{i=1}^n P_i$, for any $n \geq m$),
  $\rk_P (w) = 1+ (\rk_{P_1} (p_1) \oplus (1 + (\rk_{P_2} (p_2) \oplus
  (1+ \cdots \oplus (1+ \rk_{P_m} (p_m))))))$ if $m \neq 0$, $0$ otherwise.
\end{lemma}
\proof The rank of $w$ in $\bigtriangleright_{n=1}^{+\infty} P_n$ or
in $\bigtriangleright_{i=1}^n P_i$, for any $n \geq m$, is the same,
because the words smaller than $w$ in any of those spaces are the
same.  Hence it suffices to prove that, given
$P = \bigtriangleright_{i=1}^n P_i$, for some fixed, but arbitrary
natural number $n$, the rank of any word $w \eqdef p_1 \cdots p_m$
with $m \leq n$ in $P$ is as indicated.  We realize that $P$ is
order-isomorphic to
$\{\epsilon\} \triangleright (P_1 \triangleright (P_2 \triangleright
\cdots \triangleright P_n))$, where $w$ is mapped to
$(\epsilon, (p_1, \cdots, (p_{m-1}, p_m)))$.  The result then follows
from Lemma~\ref{lemma:triangleright}, by showing that the rank of
$(p_i, (p_{i+1}, \cdots, (p_{m-1}, p_m)))$ in the appropriate space is
$\rk_{P_i} (p_i) \oplus (1+ (\rk_{P_{i+1}} (p_{i+1}) \oplus \cdots
\oplus (1+ (\rk_{P_{m-1}} (p_{m-1}) \oplus (1+ \rk_{P_m} (p_m))))))$,
by induction on $m-i$.  \qed

The baroque formula of Lemma~\ref{lemma:triangle:rank:finite}
simplifies as follows.
\begin{lemma}
  \label{lemma:triangle:rank:finite:alt}
  In the situation of Lemma~\ref{lemma:triangle:rank:finite},
  $\rk_P (w)=\bigoplus_{i=1}^m \rk_{P_i} (p_i) + (m-k)$, where $k$ is
  the largest number between $1$ and $m$ such that $\rk_{P_k} (p_k)$
  is infinite, or $0$ if there is none.
\end{lemma}
\proof The claim is obvious if $m=0$.  Let us therefore assume
$m \neq 0$.  Let $\alpha_i$ be defined by
$\alpha_m \eqdef 1+\rk_{P_m} (p_m)$,
$\alpha_i \eqdef 1+ (\rk_{P_i} (p_i) \oplus \alpha_{i+1})$ for every
$i$ with $1\leq i< m$.  Lemma~\ref{lemma:triangle:rank:finite} states
that $\rk_P (w) = \alpha_1$.

Since $\rk_{P_{k+1}} (p_{k+1})$, \ldots, $\rk_{P_m} (p_m)$ are all
finite, we have $\alpha_m = \rk_{P_m} (p_m)+1$,
$\alpha_{m-1} = \rk_{P_{m-1}} (p_{m-1}) + \rk_{P_m} (p_m) + 2$,
\ldots, $\alpha_{k+1} = \sum_{i=k+1}^m \rk_{P_i} (p_i) + (m-k)$.  We
may rewrite the latter as
$\alpha_{k+1} = \bigoplus_{i=k+1}^m \rk_{P_i} (p_i) + (m-k)$.  If $k=0$,
namely if every number $\rk_{P_i} (p_i)$ is finite, this proves the
claim directly.

Otherwise,
$\alpha_k = 1+ (\rk_{P_k} (p_k) \oplus \bigoplus_{i=k+1}^m \rk_{P_i}
(p_i) + (m-k))$.  Since $\rk_{P_k} (p_k)$ is infinite, this is also
equal to
$\rk_{P_k} (p_k) \oplus \bigoplus_{i=k+1}^m \rk_{P_i} (p_i) + (m-k) =
\bigoplus_{i=k}^m \rk_{P_i} (p_i) + (m-k)$.  Also, $\alpha_k$ is
itself infinite, so
$\alpha_{k-1} = 1+(\rk_{P_{k-1}} (p_{k-1}) \oplus \alpha_k)$ is equal
to
$1+\bigoplus_{i=k-1}^m \rk_{P_i} (p_i) + (m-k) = \bigoplus_{i=k-1}^m
\rk_{P_i} (p_i) + (m-k)$, and is also infinite.  Similarly, we show
that $\alpha_j = \bigoplus_{i=j}^m \rk_{P_i} (p_i) + (m-k)$ for every
$j$ with $1\leq j \leq k$, by descending induction on $j$.  The claim
follows by taking $j \eqdef 1$.  \qed

We define the infinite natural sum $\bigoplus_{n=1}^{+\infty} \alpha_n$
as the least upper bound of the increasing sequence of ordinals
$\bigoplus_{i=1}^n \alpha_i$, $n \in \nat$.  The following is more
easily stated using reduced sobrification ranks.
\begin{theorem}
  \label{thm:dim:prefix}
  Given countably many non-empty Noetherian spaces $X_1$, $X_2$,
  \ldots, $X_n$, \ldots{}, $\rsob (\bigtriangleright_{n=1}^{+\infty}
  X_n)$ is equal to:
  \begin{enumerate}
  \item $\bigoplus_{n=1}^k \rsob {X_n} + \omega$ if $\rsob {X_n}$ is
    finite for $n$ large enough, where $k$ is the largest index such
    that $\rsob {X_k}$ is infinite, or $0$ if there is none;
  \item $\bigoplus_{n=1}^{+\infty} \rsob {X_n}$ otherwise.
  \end{enumerate}
\end{theorem}
\proof For all $C_1 \in \Sober {X_1}$, \ldots, $C_m \in \Sober {X_m}$,
the rank of $\lceil C_1 \cdots C_m \rangle$ in
$\Sober (\bigtriangleright_{n=1}^{+\infty} X_n)$ is the same as the
rank of the word $C_1 \cdots C_m$ in the isomorphic poset in
$(\bigtriangleright_{n=1}^{+\infty} \Sober {X_n})^\top$, hence also in
the latter minus its top element $\top$, which is
$\bigtriangleright_{n=1}^{+\infty} \Sober {X_n}$.  By
Lemma~\ref{lemma:triangle:rank:finite:alt}, this is
$\bigoplus_{i=1}^m \rk_{\Sober {X_i}} (C_i) + (m-k)$, where $k$ is the
largest number between $1$ and $m$ such that
$\rk_{\Sober {X_i}} (C_i)$ is infinite, or $0$ if there is none.

When $C_1$, \ldots, $C_m$ vary (with $m$ fixed), the latter reaches a
maximum at $\alpha_m \eqdef \bigoplus_{i=1}^m \rsob {X_i} + (m-k)$,
using Definition~\ref{defn:rdim}, where $k$ is the largest number
between $1$ and $m$ such that $\sob {X_i}$ is infinite, $0$ otherwise.
Then $\rsob (\bigtriangleright_{n=1}^{+\infty} X_n)$ is the rank of
the additional top element of
$(\bigtriangleright_{n=1}^{+\infty} \Sober {X_n})^\top$, which is the
supremum of those values as $m$ varies.

\emph{Case 1.}  If there are infinitely many indices $m$ such that
$\sob {X_m}$ is infinite, then the rank of the top element is also the
supremum of the cofinal subfamily of the values $\alpha_m$, where $m$
ranges over the indices such that $\sob {X_m}$ is infinite; those
values are equal to $\bigoplus_{i=1}^m \rk_{\Sober {X_i}} (C_i)$, by
definition and the condition on $\sob {X_m}$, so the rank of the top
element is $\bigoplus_{n=1}^{+\infty} \rsob {X_n}$, and this is the
desired reduced sobrification rank of
$\bigtriangleright_{n=1}^{+\infty} X_n$.

\emph{Case 2.}  Otherwise, there is a number $m_0$ such that
$\rsob {X_m}$ is finite for every $m \geq m_0$.  Let $m_0$ be the
least one, and $k \eqdef m_0-1$.  For every $m \geq m_0$,
$\alpha_m = \bigoplus_{i=1}^m \rsob {X_i} + (m-k)$.  We write the
latter as $\bigoplus_{i=1}^k \rsob {X_i} + a_m$, where $a_m$ is the
natural number $\sum_{i=m_0}^m \rsob {X_i} + (m-k)$.  The first summand
$\bigoplus_{i=1}^k \rsob {X_i}$ is independent of $m$, and the
supremum of the numbers $a_m$, $m \geq m_0$, is equal to $\omega$.
Therefore the rank of the top element is
$\bigoplus_{n=1}^k \rsob {X_n} + \omega$.  \qed

One may simplify the formulae of Theorem~\ref{thm:dim:prefix} as
follows.
\begin{corollary}
  \label{corl:dim:prefix:>=2}
  Given countably many non-empty Noetherian spaces $X_1$, $X_2$,
  \ldots, $X_n$, \ldots{}, such that $\rsob {X_n} \geq 1$ for
  infinitely many indices $n \geq 1$, we have
  $\rsob (\bigtriangleright_{n=1}^{+\infty} X_n) =
  \bigoplus_{n=1}^{+\infty} \rsob {X_n}$.
\end{corollary}
\proof We only have to prove this in case~1 of
Theorem~\ref{thm:dim:prefix}.  Then
$\bigoplus_{n=k+1}^{+\infty} \rsob {X_n}$ is an infinite sum of
natural numbers, infinitely of which are non-zero, and is therefore
equal to $\omega$.  It follows that
$\bigoplus_{n=1}^k \rsob {X_n} + \omega = \bigoplus_{n=1}^{+\infty}
\rsob {X_n}$.
\qed

\begin{remark}
  \label{rem:dim:prefix:fin}
  When some $X_n$ is empty, let $n \in \nat$ be smallest such that
  $X_{n+1}=\emptyset$; then, the space
  $\bigtriangleright_{i=1}^{+\infty} X_i$ still makes sense, and
  coincides with $\bigtriangleright_{i=1}^n X_i$.  The same argument
  as in Theorem~\ref{thm:dim:prefix} then shows that
  $\rsob (\bigtriangleright_{n=1}^{+\infty} X_n) = \rsob
  (\bigtriangleright_{i=1}^n X_i) = \bigoplus_{i=1}^n \rsob
  {X_i}+(n-k)$, where $k$ is the largest number between $1$ and $n$
  such that $\rsob {X_k}$ is infinite, or $0$ if there is none.
  Indeed, in that case $\Sober (\bigtriangleright_{i=1}^n X_i)$ is
  isomorphic to $\bigtriangleright_{i=1}^n \Sober {X_i}$, without any
  additional top element $\top$ \cite[Proposition~9.9]{FGL:partI}; so
  the largest rank of an element of $\bigtriangleright_{i=1}^n X_i$ is
  $\bigoplus_{i=1}^n \rsob {X_i} + (n-k)$, by
  Lemma~\ref{lemma:triangle:rank:finite:alt}.
\end{remark}

When all the space $X_n$ are equal to the same space $X$, we write
$X^\triangleright$ for $\bigtriangleright_{n=1}^{+\infty} X_n$.  This
is the space of finite words over $X$, with the prefix topology.
\begin{corollary}
  \label{corl:dim:prefix:homog}
  For every non-empty Noetherian space $X$,
  $\rsob {X^\triangleright} = \omega^{\alpha_1+1}$, where $\rsob X$ is
  written in Cantor normal form as
  $\omega^{\alpha_1} + \cdots + \omega^{\alpha_m}$,
  $\alpha_1 \geq \cdots \geq \alpha_m$.  When $X$ is empty,
  $X^\triangleright=\{\epsilon\}$ and $\rsob {X^\triangleright} = 0$.
\end{corollary}

\proof The second claim is obvious.  For the first claim, either
$\rsob X=0$, in which case, by Theorem~\ref{thm:dim:prefix}, item~1,
$\rsob {X^\triangleright} = \omega$, which shows the claim; or
$\alpha \eqdef \rsob X \geq 1$, in which case
$\rsob {X^\triangleright}$ is the infinite sum
$\bigoplus_{n=1}^{+\infty} \alpha = \sup_{n \in \nat}
\bigoplus_{i=1}^n \alpha$ of the same ordinal $\alpha$, by
Corollary~\ref{corl:dim:prefix:>=2}.  We write $\alpha$ in Cantor
normal form as $\omega^{\alpha_1} + \cdots + \omega^{\alpha_m}$.  Then
$\bigoplus_{i=1}^n \alpha = \omega^{\alpha_1} \times n + \cdots +
\omega^{\alpha_m} \times n$.  In order to compute the supremum of the
latter values as $n$ varies, we note that this supremum is larger than
or equal to
$\sup_{n \in \nat} \omega^{\alpha_1} \times n = \omega^{\alpha_1}
\times \omega = \omega^{\alpha_1+1}$.  It is equal to it, because
$\omega^{\alpha_1+1}$ is already an upper bound of the family
$\{\omega^{\alpha_1} \times n + \cdots + \omega^{\alpha_m} \times n
\mid n \in \nat\}$.  \qed

As far as stature is concerned, we define the infinite natural product
$\bigotimes_{n=1}^{+\infty} \alpha_n$ of \emph{non-zero} ordinals
$\alpha_n$ as the least upper bound of the monotonic sequence of
ordinals $\bigotimes_{i=1}^n \alpha_i$, $n \in \nat$.
\begin{theorem}
  \label{thm:stature:prefix}
  Given countably many non-empty Noetherian spaces $X_1$, $X_2$,
  \ldots, $X_n$, \ldots{},
  \begin{enumerate}
  \item If $||X_n||$ is infinite for infinitely many values of
    $n \in \nat$, then
    $||\bigtriangleright_{n=1}^{+\infty} X_n|| =
    \bigotimes_{n=1}^{+\infty} ||X_n||$.
  \item If $||X_n||$ is finite for $n$ large enough, then letting $k$
    be the largest index such that $\sob {X_k}$ is infinite, or $0$ if
    there is none, then $||\bigtriangleright_{n=1}^{+\infty} X_n||$ is
    equal to $\bigotimes_{m=1}^k ||X_m|| \times \omega$.
  \end{enumerate}
\end{theorem}
The formula in (2) is not a special case of (1), unless
$||X_n|| \geq 2$ for infinitely many values of $n$.  Indeed, if
$||X_n||=1$ for all $n \geq \ell+1$ (where $\ell \geq k$), then
$\bigotimes_{n=1}^{+\infty} ||X_n|| = \bigotimes_{n=1}^\ell ||X_n|| =
\bigotimes_{n=1}^k ||X_n|| \times p$, where $p$ is the natural number
$\prod_{n=k+1}^\ell ||X_n||$; that is different from
$\bigotimes_{m=1}^k ||X_m|| \times \omega$.

\proof For all $i, j \in \nat$ with $i \leq j$, let $\bA_{ij}$ be the
subset of $\bigtriangleright_{n=i+1}^{+\infty} X_n$ of those words of
length at most $j-i$.  This is the complement of
$\lfloor \mathcal U \rangle$ where
$\mathcal U \eqdef U_i, U_{i+1}, \cdots, U_n, \cdots$ is the wide
telescope defined by $U_i \eqdef \cdots \eqdef U_j \eqdef \emptyset$,
and $U_k \eqdef \prod_{n=i+1}^k X_n$ for every $k \geq j+1$.
Therefore $\bA_{ij}$ is a closed subset of
$\bigtriangleright_{n=i+1}^{+\infty} X_n$.

When $i=j$, $\bA_{jj} = \{\epsilon\}$, so $||\bA_{jj}|| = 1$.

If $i < j$, then we claim that the map $f$ defined by
$f (\bot) \eqdef \epsilon$, $f (x_i, w) \eqdef x_i w$ is a
homeomorphism from $(X_{i+1} \times \bA_{(i+1)j})_\bot$ onto
$\bA_{ij}$.  In order to see that $f$ is continuous, we show that the
inverse image of any closed set by $f$ is closed.  Since $\bA_{ij}$ is
a proper closed subset of $\bigtriangleright_{n=i+1}^{+\infty} X_n$,
every closed subset of $\bA_{ij}$ is a proper closed subset $\bB$ of
$\bigtriangleright_{n=i+1}^{+\infty} X_n$, which is included in
$\bA_{ij}$.  Since $\bigtriangleright_{n=i+1}^{+\infty} X_n$ is
Noetherian, $\bB$ is a finite union of (proper) irreducible closed
sets.  Each one is of the form $\lceil C_{i+1} \cdots C_n \rangle$,
where $n \geq i$ and $C_{i+1}$, \ldots, $C_n$ are irreducible closed
in $X_{i+1}$, \ldots, $X_n$ respectively.  We may also assume that
$n \leq j$, since otherwise
$\lceil C_{i+1} \cdots C_n \rangle \cap \bA_{ij} = \lceil C_{i+1}
\cdots C_j \rangle$.  Now
$f^{-1} (\lceil C_{i+1} \cdots C_n \rangle) = \{\bot\} \cup (C_{i+1}
\times \lceil C_{i+2} \cdots C_n \rangle)$ if $n\geq i+1$, or
$\{\bot\}$ if $n=i$, and those sets are closed.  (Note that $\{\bot\}$
is the downwards-closure, hence the closure, of $\bot$.)  This shows
that $f$ is continuous.  In order to see that its inverse is
continuous, too, it suffices to show that every irreducible closed
subset of $(X_{i+1} \times \bA_{(i+1)j})_\bot$ is of one of the above
forms.  The irreducible closed subsets of
$(X_{i+1} \times \bA_{(i+1)j})_\bot = \{\bot\} +_\lex (X_{i+1} \times
\bA_{(i+1)j})$ are $\{\bot\}$ itself, and the sets of the form
$\{\bot\} \cup (C_{i+1} \times \lceil C_{i+2} \cdots C_n \rangle)$,
where $C_{i+1}$, \ldots, $C_n$ are irreducible closed and
$n \geq i+1$, using Lemma~\ref{lemma:S:+lex}, and the fact that the
irreducible closed subsets of a product are the products of
irreducible closed subsets (see Section~\ref{sec:products}).

Using Proposition~\ref{prop:lifting} and
Theorem~\ref{thm:stature:prod}, we obtain that
$||\bA_{ij}|| = 1+(||X_{i+1}|| \otimes ||\bA_{(i+1)j}||)$ for all
$i < j$.  By an easy induction on $j-i$, which proceeds along similar
principles as the proof of Lemma~\ref{lemma:triangle:rank:finite:alt},
we obtain that for all $i \leq j$,
$||\bA_{ij}|| = \bigotimes_{m=i+1}^k ||X_m|| \times ( 1+||X_{k+1}||+
||X_{k+1}|| \times ||X_{k+2}|| + \cdots + \prod_{p=k+1}^j ||X_p||)$,
where $k$ is the largest number between $i+1$ and $j$ such that
$||X_k||$ is infinite, and $i$ if there is no such number.  (We take
products over an empty family to be equal to $1$.)
% ...

In particular, the stature of the closed subspace
$\bA_n \eqdef \bA_{0n}$ of $\bigtriangleright_{i=1}^{+\infty} X_i$ is
$\bigotimes_{m=1}^k ||X_m|| \times ( 1+||X_{k+1}||+ ||X_{k+1}|| \times
||X_{k+2}|| + \cdots + \prod_{p=k+1}^n ||X_p||)$, where $k$ is the
largest number between $1$ and $n$ such that $||X_k||$ is infinite,
and $0$ if there is no such number.

Clearly, $\bA_n$ is a proper subset of
$\bigtriangleright_{i=1}^{+\infty} X_i$.  Every proper irreducible
closed subset of $\bigtriangleright_{i=1}^{+\infty} X_i$ is of the
form $\lceil C_1 \cdots C_n \rangle$, where $n \in \nat$ and each
$C_i$ is irreducible closed; then, $\lceil C_1 \cdots C_n \rangle$ is
included in $\bA_n$.  Every proper closed subset $\bB$ of
$\bigtriangleright_{i=1}^{+\infty} X_i$, which one can write as a
finite union of proper irreducible closed sets, is therefore also
included in some $\bA_n$.  Proposition~\ref{prop:limit} then implies
that
$||\bigtriangleright_{m=1}^{+\infty} X_m|| = \sup_{n \in \nat}
||\bA_n||$.

We now make two cases, as in the proof of
Theorem~\ref{thm:dim:prefix}.

(1) If there are infinitely many indices $m$ such that
$||X_m||$ is infinite, then the supremum on the right-hand side can
equivalently be taken over those $n \in \nat$ such that $||X_n||$ is
infinite.  For each of those,
$||\bA_n|| = \bigotimes_{m=1}^n ||X_m||$, so
$||\bigtriangleright_{m=1}^{+\infty} X_m|| =
\bigotimes_{m=1}^{+\infty} ||X_m||$.

(2) Otherwise, let $k$ be the largest number such that $||X_k||$ is
infinite, or $0$ if there is no such number.  For every $n \geq k$, we
have
$||\bA_n|| = \bigotimes_{m=1}^k ||X_m|| \times ( 1+||X_{k+1}||+
||X_{k+1}|| \times ||X_{k+2}|| + \cdots + \prod_{p=k+1}^n ||X_p||)$.
Since $||X_p|| \geq 1$ for every $p$, the least upper bound of the
latter values is $\bigotimes_{m=1}^k ||X_m|| \times \omega$.  \qed

\begin{corollary}
  \label{corl:stature:prefix:homog}
  For every Noetherian space $X$,
  \begin{enumerate}
  \item If $||X||=0$, then $||X^\triangleright||=1$;
  \item If $||X||$ is finite and non-zero, then $||X^\triangleright||=\omega$;
  \item Otherwise, $||X^\triangleright||$ is equal to
    $\omega^{\omega^{\beta_1+1}}$, where $||X||$ is written in Cantor
    normal form as $\omega^{\alpha_1} + \cdots + \omega^{\alpha_m}$
    ($\alpha_1 \geq \cdots \geq \alpha_m$), and $\alpha_1$ is written
    in Cantor normal form as
    $\omega^{\beta_1} + \cdots + \omega^{\beta_n}$
    ($\beta_1 \geq \cdots \geq \beta_n$).
  \end{enumerate}
 \end{corollary}
\proof
(1) If $||X||=0$, then $X$ is empty, so $X^\triangleright =
\{\epsilon\}$.

(2) follows from Theorem~\ref{thm:stature:prefix}~(2).

(3) By Theorem~\ref{thm:stature:prefix}~(1),
$||X^\triangleright|| = \otimes_{n=1}^{+\infty} ||X||$.  We write
$||X||$ in Cantor normal form as
$\omega^{\alpha_1} + \cdots + \omega^{\alpha_m}$
($\alpha_1 \geq \cdots \geq \alpha_m$).  Necessarily,
$\alpha_1 \geq 1$, since $||X||$ is infinite.  For every $k \in \nat$,
$\otimes_{n=1}^k ||X||$ is a (natural) sum of terms of the form
$\omega^\alpha$, where each $\alpha$ is equal to
$\alpha_{i_1} \oplus \cdots \oplus \alpha_{i_k}$ for some tuple
$(i_1, \cdots, i_k) \in \{1, \cdots, m\}^k$.  Every such $\alpha$ is
smaller than or equal to $\alpha_1 \otimes k$.  We write $\alpha_1$ in
Cantor normal form as $\omega^{\beta_1} + \cdots + \omega^{\beta_n}$
($\beta_1 \geq \cdots \geq \beta_n$).  Then
$\alpha_1 \otimes k = \omega^{\beta_1} \times k + \cdots +
\omega^{\beta_n} \times k$ is strictly smaller than
$\omega^{\beta_1+1}$.  Hence every summand $\omega^\alpha$ of
$\otimes_{n=1}^k ||X||$ is strictly smaller than
$\omega^{\omega^{\beta_1+1}}$.  It follows that
$\otimes_{n=1}^k ||X|| < \omega^{\omega^{\beta_1+1}}$.  Taking suprema
over $k \in \nat$,
$||X^\triangleright|| \leq \omega^{\omega^{\beta_1+1}}$.

As far as the reverse inequality is concerned, for every $k \in \nat$,
$\omega^{\alpha_1 \otimes k}$ is a summand (the largest one) in
$\otimes_{n=1}^k ||X||$, so
$||X^\triangleright|| \geq \omega^{\alpha_1 \otimes k}$ for every
$k \in \nat$.  In turn,
$\alpha_1 \otimes k \geq \omega^{\beta_1} \times k$.  Taking suprema
over $k \in \nat$,
$||X^\triangleright|| \geq \omega^{\omega^{\beta_1+1}}$.  \qed

\section{Finite Multisets}
\label{sec:bf-finite-multisets}

A (finite) \emph{multiset} over a set $X$ is a map $m$ from $X$ to
$\nat$ such that $m (x)=0$ for all but finitely many elements of $X$.
We write $X^\circledast$ for the set of all (finite) multisets over
$X$.  The \emph{Parikh mapping} $\Psi \colon X^* \to X^\circledast$
maps every word $w$ to the multiset $m$ such that $m (x)$ is the
number of occurrences of $x$ in $w$ \cite{Parikh:cfl}.  We write
$\mopen x_1, \cdots, x_n \mclose$ for $\Psi (x_1 \cdots x_n)$, and
$\mempty$ for $\Psi (\epsilon)$.  Multiset union $m \uplus m'$ maps
every element $x$ to $m (x) + m' (x)$.  Clearly,
$\Psi (ww') = \Psi (w'w) = \Psi (w) \uplus \Psi (w')$ for all
$w, w' \in X^*$.

If $X$ is equipped with a preordering $\leq$, the \emph{sub-multiset}
preordering $\leq^\circledast$ is defined by
$\mopen x_1, \cdots, x_m \mclose \leq^\circledast \mopen y_1, \cdots,
y_n\mclose$ if and only if there is an injective map
$r \colon \{1, \cdots, m\} \to \{1, \cdots, n\}$ such that
$x_i \leq y_{r (i)}$ for every $i$, $1\leq i\leq m$.  If $X$ is wqo
under $\leq$, then $X^\circledast$ is wqo under $\leq^\circledast$.
We need to mention the folklore lemma.
\begin{lemma}
  \label{lemma:mset:order}
  If $\leq$ is a partial ordering on a set $P$, then
  $\leq^\circledast$ is a partial ordering on $P^\circledast$.
\end{lemma}
\proof Reflexivity and transitivity are obvious.  Let
$\mopen p_1, \cdots, p_m \mclose$ and
$\mopen q_1, \cdots, \allowbreak q_n \mclose$ be two multisets that
are less than or equal to each other with respect to
$\leq^\circledast$.  There is an injective map
$r \colon \{1, \cdots, m\} \to \{1, \cdots, n\}$ such that
$p_i \leq q_{r (i)}$ for every $i$, $1\leq i\leq m$, and there is an
injective map $s \colon \{1, \cdots, n\} \to \{1, \cdots, m\}$ such
that $q_j \leq p_{s (j)}$ for every $j$, $1\leq j\leq n$.  In
particular, $m\leq n$ and $n \leq m$, so $m=n$ and both $r$ and $s$
are permutations.  Let $f \eqdef s \circ r$.  This is also a
permutation, and $p_{f (i)} \geq p_i$ for every $i$, $1\leq i\leq n$.
By a classic trick of finite group theory, for each $i$, there is a
number $k \geq 1$ such that $f^k (i)=i$.  Then,
$p_i = p_{f^k (i)} \geq p_{f^{k-1} (i)} \geq \cdots \geq p_{f (i)}
\geq p_i$, so that all those elements are equal.  It follows that
$p_{f (i)}=p_i$ for every $i$.  Hence, for every $i$, $1\leq i\leq m$,
$p_i = p_{f (i)} = p_{s (r (i))} = q_{r (i)}$.  Since $r$ is a
permutation, we have just shown that $p_1, \cdots, p_m$ and
$q_1, \cdots, q_n$ are the same list up to permutation, namely, that
$\mopen p_1, \cdots, p_m \mclose = \mopen q_1, \cdots, q_n \mclose$.
\qed

It turns out that, if $X$ is wpo, and $\alpha$ is its maximal order
type, written in Cantor normal form as
$\omega^{\alpha_1} + \cdots + \omega^{\alpha_n}$, then the maximal
order type of $X^\circledast$ is
$\omega^{\omega^{\adjsw{\alpha_1}} + \cdots +
  \omega^{\adjsw{\alpha_n}}}$, where $\adjsw\alpha$ was defined in
Definition~\ref{defn:adjsw}.  This was proved by van der Meeren,
Rathjen and Weiermann in \cite[Theorem~5]{vdMRW:bigVeblen}, correcting
a previous statement \cite[Theorem~2]{Weiermann:multiset}.
\begin{definition}[$\adjmul\alpha$]
  \label{defn:adjmul}
  For every ordinal $\alpha$, written in Cantor normal form as
  $\omega^{\alpha_1} + \cdots + \omega^{\alpha_m}$ with
  $\alpha \geq \alpha_1 \geq \cdots \geq \alpha_m$, $m \in \nat$, the
  ordinal $\adjmul\alpha$ is defined as
  $\omega^{\adjsw{\alpha_1}} + \cdots + \omega^{\adjsw{\alpha_m}}$.
\end{definition}
Hence van der Meeren, Rathjen and Weiermann's result is that
$o (X^\circledast) = \omega^{\adjmul{o (X)}}$.

As with other constructions, we wish to extend this result to
Noetherian spaces, computing the stature $||X^\circledast||$, but also
the sobrification rank $\sob {X^\circledast}$.

For every topological space $X$, and given that $X^*$ has the word
topology, we topologize $X^\circledast$ with the quotient topology
induced by the Parikh mapping $\Psi$ \cite[Exercise~9.7.35]{JGL:topo}.
The following can be found in \cite[Proposition~8.2]{FGL:partI}.  If
$X$ is Noetherian, then so is $X^\circledast$, a base of the topology
on $X^\circledast$ is given by the sets
$\langle U_1, \cdots, U_n \rangle \eqdef \{m \uplus \mopen x_1,
\cdots, x_n\mclose \mid m \in X^\circledast, x_1 \in U_1, \cdots, x_n
\in U_n\}$, where $U_1$, \ldots, $U_n$ range over the open subsets of
$X$ and $n \in \nat$.  The specialization preordering of
$X^\circledast$ is $\leq^\circledast$, where $\leq$ is the
specialization preordering of $X$.  Moreover, if $X$ has the
Alexandroff topology of $\leq$, then $X^\circledast$ has the
Alexandroff topology of $\leq^\circledast$.

Still assuming $X$ Noetherian, the irreducible closed subsets of
$X^\circledast$ are the \emph{m-products}
$F \mid C_1, \cdots, C_n \eqdef \Psi (F^* C_1^? \cdots C_n^?)$, where
$F$ is a closed subset of $X$, $n \in \nat$, and $C_1$, \ldots, $C_n$
are irreducible closed subsets of $X$,
\cite[Proposition~8.4]{FGL:partI}; $F \mid C_1, \cdots, C_n$ is the
set of multisets that contain as many elements from $F$ as one wishes,
plus at most one from each $C_i$, $1\leq i\leq n$.  When $n=0$, this
can be written $F \mid$; alternatively, this is just $F^\circledast$.

We will later use the notation $F \mid C_1, \cdots, C_n$ even when
$C_1$, \ldots, $C_n$ are closed but not necessarily irreducible, to
denote the set of multisets containing an arbitrary number of elements
of $F$, plus at most one from each $C_i$ (e.g., in the proof of
Proposition~\ref{prop:stature:mset:ubound}, and again later).  In that
case, $F \mid C_1, \cdots, C_n$ is not necessarily irreducible.

The inclusion relation between m-products is characterized as follows
\cite[Lemma~8.6]{FGL:partI}: $F \mid C_1, \cdots, C_m$ is included in
$F' \mid C'_1, \cdots, C'_n$ if and only if $F \subseteq F'$ and
$\mopen C_{i_1}, \cdots, \allowbreak C_{i_k} \mclose
\subseteq^\circledast \mopen C'_1, \cdots, C'_n \mclose$, where
$C_{i_1}$, \ldots, $C_{i_k}$ enumerates those irreducible closed sets
$C_i$ that are not included in $F'$.

It is easy that every m-product can be written in \emph{reduced} form
as $F \mid C_1, \cdots, C_m$, where no $C_i$ is included in $F$.
Indeed, if $C \subseteq F$, then we can remove $C$ from the notation
$F \mid C, C_1, \cdots, C_m$, obtaining $F \mid C_1, \cdots, C_m$, and
that denotes the same set of multisets.  Reduced forms are canonical
forms, as we now see.
\begin{lemma}
  \label{lemma:mset:reduced}
  If $F \mid C_1, \cdots, C_m$ and $F' \mid C'_1, \cdots, C'_n$ are
  two equal m-products in reduced form, then $F=F'$, $m=n$, and the
  list $C_1, \cdots, C_m$ is equal to $C'_1, \cdots, C'_n$ up to
  permutation.
\end{lemma}
\proof Since each one is included in the other, we have
$F \subseteq F'$ and $F' \subseteq F$, hence $F=F'$.  Additionally, we
have
$\mopen C_{i_1}, \cdots, \allowbreak C_{i_k} \mclose
\subseteq^\circledast \mopen C'_1, \cdots, C'_n \mclose$, where
$C_{i_1}$, \ldots, $C_{i_k}$ enumerates those irreducible closed sets
$C_i$ that are not included in $F'$.  Since $F=F'$ and the first
m-product is in reduced form, that simply means that
$\mopen C_1, \cdots, \allowbreak C_m \mclose \subseteq^\circledast
\mopen C'_1, \cdots, C'_n \mclose$.  We obtain the reverse inequality
in a similar fashion, and then we conclude since
$\subseteq^\circledast$ is antisymmetric, by
Lemma~\ref{lemma:mset:order}.  \qed

\begin{lemma}
  \label{lemma:mset:->}
  Let $X$ be a Noetherian space, and $\leadsto$ be the smallest
  relation such that, for every m-product $F \mid C_1, \cdots, C_m$ in
  reduced form (up to permutation of $C_1$, \ldots, $C_m$):
  \begin{enumerate}
  \item
    $(F \mid C_1, \cdots, C_m) \leadsto (F \mid C_1, \cdots, C_m, C)$,
    for every irreducible closed subset $C$ of $X$ not included in
    $F$;
  \item $(F \mid C_1, \cdots, C_i, \cdots, C_m) \leadsto (F \mid C_1,
    \cdots, C'_i, \cdots, C_m)$ for every irreducible closed subset
    $C'_i$ that contains $C_i$ strictly, for every $i$ with $1\leq
    i\leq m$;
  \item
    $(F \mid C_1, \cdots, C_m) \leadsto (F \cup C \mid C_1, \cdots,
    C_k)$ for every $k$ with $0 \leq k < m$, provided that
    $C_{k+1} = \cdots = C_m = C$ and $C_1$, \ldots, $C_k$ are not
    included in $F \cup C$.
  \end{enumerate}
  Let also $\leadsto^+$ be the transitive closure of $\leadsto$.  For
  all m-products $P$ and $P'$ in reduced form, $P \subsetneq P'$ if
  and only $P \leadsto^+ P'$.
\end{lemma}
\proof We note that if $P$ is in reduced form and if $P \leadsto P'$,
then $P'$ is in reduced form as well, and $P \subseteq P'$.
Additionally, $P \neq P'$, using Lemma~\ref{lemma:mset:reduced}.

Conversely, let $P \eqdef (F \mid C_1, \cdots, C_m)$ and
$P' \eqdef (F' \mid C'_1, \cdots, C'_n)$ be in reduced form, and let
us assume that $P \subseteq P'$.  In particular, $F \subseteq F'$.
We show that $P \leadsto^* P'$ by induction on $p$, where $\leadsto^*$
is the reflexive transitive closure of $\leadsto$, and where $p$ is
the number of components of $F'$ that are not included in $F$.  This
will prove our claim, since if $P \neq P'$, then the number of
$\leadsto$ steps from $P$ to $P'$ cannot be $0$.

If $p=0$, then $F'$ is included in $F$, so $F=F'$.  There is an
injective map $r \colon \{1, \cdots, m\} \to \{1, \cdots, n\}$ such
that $C_i \subseteq C'_{r (i)}$ for every $i \in \{1, \cdots, m\}$,
using the fact that $P$ is in reduced form.  We then have:
\begin{align*}
  P & = (F \mid C_1, \cdots, C_m) \\
  & \leadsto^* (F \mid C'_{r (1)}, \cdots, C'_{r (m)})
  & \text{by (2), used at most $m$ times} \\
    & \leadsto^* (F' \mid C'_1, \cdots, C'_{r (1)}, \cdots, C'_{r
      (m)}, \cdots, C'_n) 
  & \text{by (1), used $n-m$ times} \\
  & = P'.
\end{align*}

If $p \geq 1$, then let us pick one component $C$ of $F'$ that is not
included in $F$.  Up to permutation, we may assume that the indices
$i \in \{1, \cdots, m\}$ such that $C_i \subseteq C$ are exactly those
between $k+1$ and $m$, where $0 \leq k \leq m$.  Then:
\begin{align*}
  P
  & = (F \mid C_1, \cdots, C_k, C_{k+1}, \cdots, C_m) \\
  & \leadsto^* (F \mid C_1, \cdots, C_k, \underbrace{C, \cdots, C,
    C}_{m-k \text{ copies}})
  & \text{using (2) at most $m-k$ times} \\
  & \leadsto (F \mid C_1, \cdots, C_k, \underbrace{C, \cdots, C,
    C}_{m-k+1 \text{ copies}})
  & \text{using (1)} \\
  & \leadsto (F \cup C \mid C_1, \cdots, C_k)
  & \text{by (3) (see below)} \\
  & \leadsto^* P' & \text{by induction hypothesis.}
\end{align*}
The use of (3) is justified by the fact that $C_1$, \ldots, $C_k$ are
not included in $C$, and are not included in $F$ since $P$ is in
reduced form; so they are not included in $F \cup C$ either, by
irreducibility.  The use of (1) in the previous step guarantees that
the number of copies of $C$ that we will move from the right-hand side
to the left-hand side of $\mid$ in the last step is at least $1$, as
is also required in order to apply (3).  \qed

\subsection{Bounds on the sobrification rank of $X^\circledast$}
\label{sec:dimension-x-ast}

We recall that, for every ordinal $\alpha$, written in Cantor normal
form as $\omega^{\alpha_1} + \cdots + \omega^{\alpha_m}$, the ordinal
$\alpha \otimes \omega$ is equal to
$\omega^{\alpha_1+1} + \cdots + \omega^{\alpha_m+1}$.

\begin{lemma}
  \label{lemma:mset:+}
  For all ordinals $\alpha$, $\beta$ and natural numbers $k$, if
  $\beta \neq 0$ then
  $(\alpha \otimes \omega) \oplus (\beta \otimes k) < (\alpha+\beta)
  \otimes \omega$.
\end{lemma}
\proof We write $\alpha$ in Cantor normal form as
$\omega^{\alpha_1} + \cdots + \omega^{\alpha_m}$, with
$\alpha \geq \alpha_1 \geq \cdots \geq \alpha_m$, $m \in \nat$, and
similarly $\beta$ as $\omega^{\beta_1} + \cdots + \omega^{\beta_n}$,
with $\beta \geq \beta_1 \geq \cdots \geq \beta_n$.  Since
$\beta \neq 0$, it follows that $n \geq 1$.  Let $i$ be the largest
index in $\{1, \cdots, m\}$ such that $\alpha_i \geq \beta_1$, or $0$
if there is none.  Then $\omega^{\alpha_m} + \omega^{\beta_1}$,
\ldots, $\omega^{\alpha_{i+1}} + \omega^{\beta_1}$ are all equal to
$\omega^{\beta_1}$, so
$\alpha + \beta = \omega^{\alpha_1} + \cdots + \omega^{\alpha_i} +
\omega^{\beta_1} + \cdots + \omega^{\beta_n}$, and therefore
$(\alpha + \beta) \otimes \omega = \omega^{\alpha_1+1} + \cdots +
\omega^{\alpha_i+1} + \omega^{\beta_1+1} + \cdots +
\omega^{\beta_n+1}$.  We compare this to
$(\alpha \otimes \omega) \oplus (\beta \otimes k)$, which is equal to
$\omega^{\alpha_1+1} + \cdots + \omega^{\alpha_i+1} + \gamma$, where
$\gamma \eqdef (\omega^{\alpha_{i+1}+1} + \cdots \omega^{\alpha_m+1})
\oplus (\omega^{\beta_1} \times k + \cdots + \omega^{\beta_n} \times
k)$.  The latter is a sum of terms of the form $\omega^\delta$ with
$\delta < \beta_1+1$, so
$\gamma < \omega^{\beta_1+1} \leq \omega^{\beta_1+1} + \cdots +
\omega^{\beta_n+1}$.  The claim follows.  \qed

\begin{lemma}
  \label{lemma:rk:chasles}
  Let $X$ be a topological space, $F$ be a closed subset of $X$, and
  $C$, $C'$ be two irreducible closed subsets of $X$.  Then:
  \begin{enumerate}
  \item $C \diff F$ is empty or irreducible closed in $X \diff F$;
  \item $C \diff F \subseteq C' \diff F$ if and only if $C \subseteq
    F$ or $C \subseteq C'$;
  \item If $X$ is Noetherian, then
    $\rk_{\Sober X \cup \{\emptyset\}} (C) \leq \rk_{\Sober X \cup
      \{\emptyset\}} (C') + \rk_{\Sober (X \diff C') \cup
      \{\emptyset\}} (C \diff C')$.
  \end{enumerate}
\end{lemma}
\proof (1) $C \diff F$ is closed in $X \diff F$.  Let us assume that
it is non-empty, namely that $C$ is not included in $F$.  We verify
that $C \diff F$ is irreducible.  If $C \diff F$ is included in the
union of two closed subsets $F_1 \diff F$ and $F_2 \diff F$ of
$X \diff F$ (where $F_1$ and $F_2$ are closed in $X$), then $C$ is
included in $F \cup F_1 \cup F_2$.  Since $C$ is irreducible in $X$
and not included in $F$, it is included in $F_1$ or in $F_2$.  Then,
$C \diff F$ is also included in $F_1 \diff F$, or in $F_2 \diff F$.

(2) The if direction is clear.  Conversely, if $C \diff F \subseteq C'
\diff F$, then $C \subseteq F \cup C'$, and the claim follows from the
irreducibility of $C$.

(3) We fix $C'$, and we prove the inequality by well-founded induction
on $C \diff C'$.  If $C \diff C'$ is empty, then
$\rk_{\Sober (X \diff C') \cup \{\emptyset\}} (C \diff C') = 0$, and
$C \subseteq C'$, so
$\rk_{\Sober X \cup \{\emptyset\}} (C) \leq \rk_{\Sober X \cup
  \{\emptyset\}} (C')$.  Otherwise, it suffices to show that
$\alpha < \rk_{\Sober X \cup \{\emptyset\}} (C') + \rk_{\Sober (X
  \diff C') \cup \{\emptyset\}} (C \diff C')$ for every ordinal
$\alpha < \rk_{\Sober X \cup \{\emptyset\}} (C)$.  By the inductive
definition of rank, there is a $C'' \in \Sober X \cup \{\emptyset\}$
such that $C'' \subsetneq C$ and
$\rk_{\Sober X \cup \{\emptyset\}} (C'') = \alpha$.  By induction
hypothesis,
$\alpha = \rk_{\Sober X \cup \{\emptyset\}} (C'') \leq \rk_{\Sober X
  \cup \{\emptyset\}} (C') + \rk_{\Sober (X \diff C') \cup
  \{\emptyset\}} (C'' \diff C')$.  Now $C'' \diff C'$ is a subset of
$C \diff C'$, and a proper one: otherwise
$C \diff C' \subseteq C'' \diff C'$, and since $C \not\subseteq C'$,
this would imply $C \subseteq C''$ by (2).  Hence, using the fact that
$+$ is strictly monotonic in its second argument,
$\alpha < \rk_{\Sober X \cup \{\emptyset\}} (C') + \rk_{\Sober (X
  \diff C') \cup \{\emptyset\}} (C \diff C')$, as claimed.  \qed

\begin{lemma}
  \label{lemma:dim:mset:phi}
  Let $X$ be a Noetherian space.  For every m-product
  $P \eqdef F \mid C_1, \cdots, C_m$ in reduced form, let
  $\varphi (P) \eqdef (||F|| \otimes \omega) \oplus \bigoplus_{i=1}^m
  \rk_{\Sober (X \diff F) \cup \{\emptyset\}} (C_i \diff F)$.  Then
  $\varphi$ is a strictly monotonic ordinal-valued map with domain
  $\Sober (X^\circledast)$.
\end{lemma}
\proof We first note that $\varphi (P)$ is well defined, as
$C_i \diff F$ is an element of
$\Sober (X \diff F) \cup \{\emptyset\}$, by
Lemma~\ref{lemma:rk:chasles}~(1).

Using Lemma~\ref{lemma:mset:->}, it suffices to show that, for all
m-products $P$ and $P'$ in reduced form such that $P \leadsto P'$,
$\varphi (P) < \varphi (P')$.  This is obvious for steps of the form
(1) or (2); note that in the case of form (1),
$C \diff F \neq \emptyset$ by assumption, so that
$\rk_{\Sober (X \diff F) \cup \{\emptyset\}} (C \diff F) \neq 0$.

For steps of the form (3), $P$ is of the form
$F \mid C_1, \cdots, C_m$, $P' = (F \cup C \mid C_1, \cdots, C_k)$,
$0\leq k< m$, and $C_{k+1} = \cdots = C_m = C$.  (Note also that the
condition that $C_1$, \ldots, $C_k$ are not included in $F \cup C$
implies that $P'$ is written in reduced form.)  Then, letting
$\beta \eqdef \rk_{\Sober (X \diff F) \cup \{\emptyset\}} (C \diff
F)$, we have:
\begin{align*}
  \varphi (P)
  & = (||F|| \otimes \omega) \oplus
    (\beta \otimes (m-k)) \oplus
    \bigoplus_{i=1}^k \rk_{\Sober (X \diff F) \cup \{\emptyset\}} (C_i
    \diff F).
\end{align*}
Using Lemma~\ref{lemma:rk:chasles}~(3), and noticing that
$(C_i \diff F) \diff (C \diff F) = C_i \diff (F \cup C)$, the term
$\rk_{\Sober (X \diff F) \cup \{\emptyset\}} (C_i \diff F)$ is smaller
than or equal to
$\beta + \rk_{\Sober (X \diff (F \cup C)) \cup \{\emptyset\}} (C_i
\diff (F \cup C))$.  The latter is smaller than or equal to
$\beta \oplus \rk_{\Sober (X \diff (F \cup C)) \cup \{\emptyset\}}
(C_i \diff (F \cup C))$, since the inequality
$\beta+\gamma \leq \beta \oplus\gamma$ holds for all ordinals
$\gamma$.
Using the associativity and the commutativity of $\oplus$, we obtain:
\begin{align*}
  \varphi (P)
  & \leq (||F|| \otimes \omega) \oplus
    (\beta \otimes m) \oplus 
    \bigoplus_{i=1}^k \rk_{\Sober (X \diff (F \cup C)) \cup \{\emptyset\}} (C_i \diff (F
\cup C)).
\end{align*}
We note that, since $P$ is reduced, $C = C_m$ is not included in $F$,
so $\beta = \rk_{\Sober (X \diff F) \cup \{\emptyset\}} (C \diff F)$
is non-zero.  Therefore Lemma~\ref{lemma:mset:+} applies, to the
effect that
$(||F|| \otimes \omega) \oplus (\beta \otimes m) < (||F|| + \beta)
\otimes \omega$.  It is easy to see that $\beta \leq ||C \diff F||$,
so
$||F|| + \beta \leq ||F|| + ||C \diff F|| = ||F|| + ||(F \cup C) \diff
F|| \leq ||F \cup C||$, using Lemma~\ref{lemma:stature:C}~(2).  It
follows that:
\begin{align*}
  \varphi (P)
  & < (||F \cup C|| \otimes \omega) \oplus 
    \bigoplus_{i=1}^k \rk_{\Sober (X \diff (F \cup C)) \cup \{\emptyset\}} (C_i \diff (F
    \cup C)) = \varphi (P').
\end{align*}
This concludes the proof.    \qed

\begin{proposition}
  \label{prop:dim:mset:ubound}
  For every Noetherian space $X$,
  $\sob {X^\circledast} \leq (||X|| \otimes \omega) +1$, or
  equivalently, $\rsob {X^\circledast} \leq ||X|| \otimes \omega$.
\end{proposition}
\proof We use Lemma~\ref{lemma:dim:mset:phi}, and we observe that the
largest value taken by $\varphi$ is obtained as
$\varphi (X \mid) = ||X|| \otimes \omega$.  Therefore $\varphi$ takes
its values in (the set of ordinals strictly below)
$(||X|| \otimes \omega)+1$.  \qed

Given an ordinal $\alpha$, written in Cantor normal form as
$\omega^{\alpha_1} + \cdots + \omega^{\alpha_m}$ with
$\alpha \geq \alpha_1 \geq \cdots \geq \alpha_m$, $m \in \nat$, the
ordinal $\omega \times \alpha$ is equal to
$\omega^{1+\alpha_1} + \cdots + \omega^{1+\alpha_m}$.  This is in
general smaller than or equal to 
$\alpha \otimes \omega = \omega^{\alpha_1+1} + \cdots +
\omega^{\alpha_m+1}$, and equal to it if and only if every $\alpha_i$
is finite, if and only if $\alpha < \omega^\omega$.

\begin{proposition}
  \label{prop:dim:mset:lbound}
  For every Noetherian space $X$,
  $\sob {X^\circledast} \geq (\omega \times ||X||)+1$; equivalently,
  $\rsob {X^\circledast} \geq \omega \times ||X||$.
\end{proposition}
\proof This is clear if $||X||=0$, namely
if $X$ is empty, since then $X^\circledast$ has exactly one element,
the empty multiset $\mopen\mclose$.  Henceforth, we assume that $X$ is
non-empty, and we prove the claim by induction on $||X||$.

For every increment $(F, F') \in \Inc (\Hoare_0 X)$, $F' \diff F$ is
the equivalence class of some point $x$ with respect to the
specialization preordering of $X$, and $F' = F \cup C$, where
$C \eqdef \dc x$, by Lemma~\ref{lemma:+1}.  We will simply write such
increments as $(F, F \cup C)$ in the sequel, without restating the
associated requirements on $F$ and $C$.  Given such an increment, for
every $n \in \nat$, we form the m-product $F \mid C^n$, meaning
$F \mid C, \cdots, C$ with $n$ copies of $C$.  We claim that the map
$f \colon \Inc (\Hoare_0 X) \times_\lex \nat \to \Sober
(X^\circledast)$ defined by $f ((F, F \cup C), n) \eqdef (F \mid C^n)$
is strictly monotonic.

Given two increments $(F, F \cup C) < (F', F' \cup C')$, we have
$f ((F, F \cup C), n) \subsetneq f ((F', F' \cup C'), n')$ for all
$n, n' \in \nat$.  Inclusion follows from the fact that both $F$ and
$C$ are included in $F'$, so
$(F \mid C^n) \subseteq (F \cup C)^\circledast \subseteq
{F'}^\circledast \subseteq (F' \mid {C'}^{n'})$.  The inclusion is
strict: given $x \in C \diff F$, the multiset that contains $n+1$
copies of $x$ is in $F' \mid {C'}^{n'}$ since in $(F \cup
C)^\circledast \subseteq {F'}^\circledast$,
but not in $F \mid C^n$.

Next, we show that $f ((F, F \cup C), n) = F \mid C^n$ is strictly
included in $F \mid C^{n'}$ for all $n < n'$.  The inclusion is
obvious.  Given $x \in C \diff F$, the multiset containing exactly
$n'$ copies of $x$ is in $F \mid C^{n'}$ but not in $F \mid C^n$.

Since $f$ is strictly monotonic,
$\rk_{\Hoare_0 (X^\circledast)} (F \mid C^n) \geq \rk_{\Inc (\Hoare_0
  X) \times_\lex \nat} ((F, \allowbreak F \cup C), n) = \omega \times
\rk_{\Inc (\Hoare_0 X)} (F, F \cup C) + n$.  When $(F, F \cup C)$
varies and $n=0$, the supremum of those quantities is
$\omega \times ||X||$, using Lemma~\ref{prop:Inc} and the fact that
$|\Hoare_0 X| = ||X||+1$.  Therefore
$\sob {X^\circledast} \geq (\omega \times ||X||)+1$, proving the
claim.  \qed

As we said, the lower bound and the upper bound only match provided
that $||X|| < \omega^\omega$.  We state the following for the record.
\begin{theorem}
  \label{thm:dim:mset:small}
  Let $X$ be a Noetherian space such that $||X|| < \omega^\omega$,
  namely $||X||$ is of the form
  $\omega^{n_1} + \cdots + \omega^{n_m}$, where
  $n_1 \geq \cdots \geq n_m$ are natural numbers.  Then
  $\sob {X^\circledast} = \omega^{n_1+1} + \cdots + \omega^{n_m+1}+1$.
\end{theorem}

\subsection{The bounds on $\sob {X^\circledast}$ are tight.}
\label{sec:bound-dim-x}

The upper and lower bounds of Proposition~\ref{prop:dim:mset:ubound}
and Proposition~\ref{prop:dim:mset:lbound} do not match, unless
$||X|| < \omega^\omega$.  This cannot be improved upon, as we now see,
by realizing that both the lower bound and the upper bound can be
attained.
\begin{proposition}
  \label{prop:dim:mset:ubound:tight}
  For every Noetherian space $Z$, let $X \eqdef Z^*$ and
  $\alpha \eqdef ||Z||$.  If $\alpha \geq 2$, then the lower bound on
  $\sob {X^\circledast}$ given in Lemma~\ref{prop:dim:mset:lbound} is
  attained:
  $\sob {X^\circledast} = \omega^{\omega^{\alpha'}} + 1 = (\omega
  \times ||X||)+1$.
\end{proposition}
\proof For every proper irreducible closed subset $P$ of
$X^\circledast$, $P$ is an m-product $F \mid C_1, \cdots, C_m$, which
we write in reduced form.  Since $P \neq X^\circledast$, $F$ is a
proper closed subset of $X = Z^*$, so $||F|| < ||Z^*||$.  By
Theorem~\ref{thm:stature:X*}, $||Z^*|| = \omega^{\omega^{\alpha'}}$.
Hence $||F|| < \omega^{\omega^{\alpha'}}$.  We write $||F||$ in Cantor
normal form as $\omega^{\beta_1} + \cdots + \omega^{\beta_n}$, with
$\omega^{\alpha'} > \beta_1 \geq \cdots \geq \beta_n$.  Since
$\omega^{\alpha'}$ is additively indecomposable and
$\omega^{\alpha'} > 1$ (we recall that $\alpha \geq 2$, so
$\alpha' \geq 1$), we also have
$\omega^{\alpha'} > \beta_1+1 \geq \cdots \geq \beta_n+1$.  Therefore
$||F|| \otimes \omega < \omega^{\omega^{\alpha'}}$.

By Theorem~\ref{thm:dim:X*}, $\sob X = \omega^{\adjsw\alpha}+1$.  By
Corollary~\ref{corl:stature:subspace}, $\sob (X \diff F) \leq \sob X$,
so, for every $i \in \{1, \cdots, m\}$,
$\rk_{\Sober (X \diff F)} (C_i \diff F) \leq \sob (X \diff F)-1 \leq
\omega^{\adjsw\alpha}$.  Since $\alpha \geq 2$, $\adjsw\alpha \geq 1$,
so $\omega^{\adjsw\alpha}$ is infinite, and therefore
$1+\omega^{\adjsw\alpha} = \omega^{\adjsw\alpha}$.  It follows that
$\rk_{\Sober (X \diff F) \cup \{\emptyset\}} (C_i \diff F) =
1+\rk_{\Sober (X \diff F)} (C_i \diff F) \leq \omega^{\adjsw\alpha}$.
We note that $\adjsw\alpha < \omega^{\alpha'}$.  (If $\alpha \geq 2$
is finite, then $\adjsw\alpha = \alpha$ and $\alpha' = \alpha-1$, and
$\alpha < \omega^{\alpha-1}$ because the right-hand side is infinite.
If $\alpha$ is infinite, then $\adjsw\alpha = \alpha'$, and this is
never a critical ordinal.  We have $\alpha' \leq \omega^{\alpha'}$,
and the inequality is strict, since $\alpha'$ is not critical.)  It
follows that $\omega^{\adjsw\alpha} < \omega^{\omega^{\alpha'}}$, so
$\rk_{\Sober (X \diff F) \cup \{\emptyset\}} (C_i \diff F) <
\omega^{\omega^{\alpha'}}$.

We use the map $\varphi$ of Lemma~\ref{lemma:dim:mset:phi}:
$\varphi (P)$ is the natural sum of
$||F|| \otimes \omega < \omega^{\omega^{\alpha'}}$, and of finitely
many terms
$\rk_{\Sober (X \diff F) \cup \{\emptyset\}} (C_i \diff F)$, which are
all strictly less than $\omega^{\omega^{\alpha'}}$.  Since
$\omega^{\omega^{\alpha'}}$ is $\oplus$-indecomposable,
$\varphi (P) < \omega^{\omega^{\alpha'}}$.  Adding one and taking
suprema over all proper irreducible closed subsets $P$ of
$X^\circledast$, we obtain that the rank of $X^\circledast$ in
$\Sober (X^\circledast)$ is less than or equal to
$\omega^{\omega^{\alpha'}}$.  Therefore
$\sob {X^\circledast} \leq \omega^{\omega^{\alpha'}}+1$.

By Proposition~\ref{prop:dim:mset:lbound}, we have $\sob
{X^\circledast} \geq (\omega \times ||X||)+1$.
We note that $\omega \times ||X|| = \omega^{1+\omega^{\alpha'}} =
\omega^{\omega^{\alpha'}}$, because $\omega^{\alpha'}$ is infinite,
and the result follows.
\qed

We now give a class of examples where the upper bound, instead of the
lower bound, is attained.
\begin{proposition}
  \label{prop:dim:mset:lbound:tight}
  For every ordinal $\alpha$, with its Alexandroff topology, the upper
  bound on $\sob {\alpha^\circledast}$ given in
  Lemma~\ref{prop:dim:mset:ubound} is reached:
  $\sob {\alpha^\circledast} = (||\alpha|| \otimes \omega)+1 = (\alpha
  \otimes \omega)+1$.
\end{proposition}
\proof The key is that every ordinal $\beta < \alpha$ is a closed
subset, and is irreducible provided that $\beta \neq 0$.  This will
allow us to form m-products of the form $\beta \mid \gamma^n$ for
various ordinals $\beta$ and $\gamma$, $\gamma \neq 0$.

We recall that $||\alpha|| = \alpha$, by
Lemma~\ref{lemma:ordinal}~(1).  This means that
$(\alpha \otimes \omega)+1$ is indeed the upper bound of
Proposition~\ref{prop:dim:mset:ubound}.

Let us write $\alpha$ in Cantor normal form as
$\omega^{\alpha_1} + \cdots + \omega^{\alpha_m}$ with
$\alpha \geq \alpha_1 \geq \cdots \geq \alpha_m$.  For each
$i \in \{1, \cdots, m+1\}$, let
$\overline\alpha_i \eqdef \omega^{\alpha_1} + \cdots +
\omega^{\alpha_{i-1}}$.  We note that $\overline\alpha_1=0$ and that
$\overline\alpha_{m+1} = \alpha$.

For the time being, let us fix $i \in \{1, \cdots, m\}$.  For every
$n \in \nat$ and every ordinal $\beta \leq \omega^{\alpha_i}$, we form
the m-product
$P_{i, n, \beta} \eqdef \overline\alpha_i \mid
\underbrace{\overline\alpha_{i+1}, \cdots, \overline\alpha_{i+1}}_n,
\overline\alpha_i + \beta$.  Let also
$P_{i, n} \eqdef \overline\alpha_i \mid
\underbrace{\overline\alpha_{i+1}, \cdots, \overline\alpha_{i+1}}_n$;
this is simply $P_{i, n, 0}$, written in reduced form.

For each fixed $n$, we claim that
$\rk_{\Sober (\alpha^\circledast)} (P_{i, n, \beta}) \geq \rk_{\Sober
  (\alpha^\circledast)} (P_{i, n}) + \beta$.  It suffices to observe
that the map $\beta \mapsto P_{i, n, \beta}$ is strictly monotonic,
and to induct on $\beta$.  It is clear that $\beta \leq \beta'$
implies $P_{i, n, \beta} \subseteq P_{i, n, \beta'}$.  If
$\beta < \beta'$, we reason by contradiction, and we assume that
$P_{i, n, \beta} = P'_{i, n, \beta'}$.  Since $\beta' \neq 0$,
$P_{i, n, \beta'}$ is already in reduced form.  If $\beta=0$, then the
reduced form of $P_{i, n, \beta}$ is the expression we gave for
$P_{i, n}$, which does not have the same number of terms to the right
of the vertical bar.  This would contradict the equality
$P_{i, n, \beta} = P'_{i, n, \beta'}$, by
Lemma~\ref{lemma:mset:reduced}.  Hence $\beta \neq 0$.
Lemma~\ref{lemma:mset:reduced} then implies that the lists
$\underbrace{\overline\alpha_{i+1}, \cdots, \overline\alpha_{i+1}}_n,
\overline\alpha_i + \beta$ and
$\underbrace{\overline\alpha_{i+1} \cdots, \overline\alpha_{i+1}}_n,
\overline\alpha_i + \beta'$ are equal up to permutation.  It is easy
to see that this implies
$\overline\alpha_i + \beta = \overline\alpha_i + \beta'$.  Since
ordinal addition is left-cancellative, $\beta=\beta'$, which is
impossible.

We observe that $P_{i, n+1} = P_{i, n, \omega^{\alpha_i}}$, because
$\overline\alpha_i + \omega^{\alpha_i} = \overline\alpha_{i+1}$.  It
follows that
$\rk_{\Sober (\alpha^\circledast)} (P_{i, n+1}) \geq \rk_{\Sober
  (\alpha^\circledast)} (P_{i, n}) + \omega^{\alpha_i}$.  By induction
on $n \in \nat$, and observing that
$P_{i, 0} = {\overline\alpha_i}^\circledast$, we have that
$\rk_{\Sober (\alpha^\circledast)} (P_{i, n}) \geq \rk_{\Sober
  (\alpha^\circledast)} ({\overline\alpha_i}^\circledast) +
\omega^{\alpha_i} \times n$.  Since $P_{i, n}$ is included in
${\overline\alpha_{i+1}}^\circledast$ for every $n \in \nat$, it
follows that
$\rk_{\Sober (\alpha^\circledast)}
({\overline\alpha_{i+1}}^\circledast) \geq \rk_{\Sober
  (\alpha^\circledast)} ({\overline\alpha_i}^\circledast) + \sup_{n
  \in \nat} (\omega^{\alpha_i} \times n) = \rk_{\Sober
  (\alpha^\circledast)} ({\overline\alpha_i}^\circledast) +
\omega^{\alpha_i+1}$.

We now induct on $i \in \{1, \cdots, m+1\}$, and we obtain that
$\rk_{\Sober (\alpha^\circledast)} ({\overline\alpha_i}^\circledast)
\geq \omega^{\alpha_1+1} + \cdots + \omega^{\alpha_{i-1}+1}$.  When
$i=m+1$, this implies that
$\rk_{\Sober (\alpha^\circledast)} (\alpha^\circledast) \geq
\omega^{\alpha_1+1} + \cdots + \omega^{\alpha_m+1} = \alpha \otimes
\omega$.  Hence
$\sob {\alpha^\circledast} \geq (\alpha \otimes \omega)+1$.  The
reverse inequality is by Proposition~\ref{prop:dim:mset:ubound}.  \qed

\subsection{The stature of $X^\circledast$}
\label{sec:stature-xcircledast}

Working as for Lemma~\ref{lemma:stature:X*:ubound:infinite}, we obtain
the following, non-optimal upper bound.

\begin{lemma}
  \label{lemma:stature:mset:ubound:1}
  For every non-empty Noetherian space $X$,
  $||X^\circledast|| \leq \omega^{\alpha \otimes \omega}$, where
  $\alpha \eqdef ||X||$.
\end{lemma}
\proof $X^\circledast$ is irreducible closed, so the space
$X^\circledast$ has exactly one component.  By
Proposition~\ref{prop:dim:stature}~(3),
$||X^\circledast|| \leq \omega^{\sob {X^\circledast}-1}$, and
$\sob {X^\circledast}-1 \leq ||X|| \otimes \omega$ by
Proposition~\ref{prop:dim:mset:ubound}.
\qed

In order to obtain a better upper bound, we observe the following.
\begin{lemma}
  \label{lemma:uplus}
  Let $X$ be a Noetherian space, and $F$ be a subset of $X$.
  The map $\uplus \colon F^\circledast \times (X \diff F)^\circledast
  \to X^\circledast$, which maps $(m_1, m_2)$ to $m_1 \uplus m_2$, is
  continuous and bijective.
\end{lemma}
\proof Bijectivity is obvious.  The inverse function maps every
$m \in X^\circledast$ to the pair $(m \cap F, m \diff F)$, where
$m \cap F$ denotes the multiset of those elements of $m$ that are in
$F$ (i.e., recalling that a multiset is a function with values in
$\nat$, this is the function that maps every $x \in F$ to $m (x)$, and
all other elements to $0$), and $m \diff F$ is the multiset of those
elements of $m$ that are not in $F$.

We recall that a base of the topology on $X^\circledast$ is given by
the sets $\langle U_1, \cdots, U_n \rangle$, where $U_1$, \ldots,
$U_n$ range over the open subsets of $X$ and $n \in \nat$.  Let
$(m_1, m_2)$ be such that
$m_1 \uplus m_2 \in \langle U_1, \cdots, U_n \rangle$.  We can write
$m_1 \uplus m_2$ as $m \cup \mopen x_1, \cdots, x_n\mclose$ where each
$x_i$ is in $U_i$.  We look at the indices $i$ such that $x_i \in F$.
Up to permutation, we assume that those are the indices $1$, \ldots,
$k$, where $0 \leq k\leq n$.  The open set
$\langle U_1 \cap F, \cdots, U_k \cap F\rangle \times \langle U_{k+1}
\diff F, \cdots, U_n \diff F \rangle$ is then an open neighborhood of
$(m_1, m_2)$ whose image by $\uplus$ is included in
$\langle U_1, \cdots, U_n \rangle$.  \qed

\begin{remark}
  \label{rem:uplus:inv}
  The inverse of the bijection $\uplus$ of Lemma~\ref{lemma:uplus} is
  not continuous in general, even when $F$ is closed.  The problem is
  that the map $m \mapsto m \cap F$ is not even monotonic, and every
  continuous map is necessarily monotonic with respect to the
  underlying specialization preorderings.  For example, let
  $X \eqdef \{1, 2\}$ in its Alexandroff topology, with $1 < 2$, and
  let $F \eqdef \{1\}$.  Let $m \eqdef \mopen 1\mclose$ and
  $m' \eqdef \mopen 2 \mclose$.  Then $m \leq^\circledast m'$, but
  $m \cap F = \mopen 1 \mclose \not\leq^\circledast m' \cap F =
  \mempty$.
\end{remark}


\begin{proposition}
  \label{prop:stature:mset:ubound}
  For every non-empty Noetherian space $X$,
  $||X^\circledast|| \leq \omega^{\adjmul\alpha}$, where
  $\alpha \eqdef ||X||$.
\end{proposition}
\proof We prove this by induction on $\alpha$.  We first note that
there is a continuous map
$f \colon X^\circledast \times X_\bot \to X^\circledast$, which maps
$(m, x)$ to $m \uplus \mopen x\mclose$ if $x \neq \bot$, and to $m$ if
$x=\bot$.  In order to see that $f$ is continuous, we note that for
all open subsets $U_1$, \ldots, $U_n$ of $X$,
$f^{-1} (\langle U_1, \cdots, U_n \rangle)$ is equal to the union of
$\langle U_1, \cdots, U_n \rangle \times X_\bot$ and of the sets
$\langle U_1, \cdots, U_{i-1}, U_{i+1}, \cdots, U_n \rangle \times
U_i$, $1\leq i\leq n$, which are all open.

We first deal with the case where $\alpha = \omega^\beta$ for some
ordinal $\beta$.  If $\beta=0$, namely if $\alpha=1$, then there is
just one non-empty closed subset of $X$, and that is $X$ itself.  As
such, it is irreducible, and the m-products are just
$\emptyset \mid X^n$, $n \in \nat$, plus $X \mid$.  They form a chain,
so any finite union of m-products is equal to just one of them.
In other words, the sets $\emptyset \mid X^n$, $n \in \nat$, and $X
\mid$, exhaust all the non-empty closed subsets of $X$.
Then $||X^\circledast|| = \omega$, which is indeed equal to
$\omega^{\adjmul\alpha}$, since $\adjmul\alpha = \omega^{\adjsw\beta}
= 1$.

Henceforth, we assume that $\alpha = \omega^\beta$ and that $\beta
\geq 1$.

For every proper closed subset $F$ of $X$,
$||F|| < \alpha = \omega^\beta$.  Writing $||F||$ in Cantor normal
form as $\omega^{\beta_1} + \cdots + \omega^{\beta_n}$, this entails
that $\beta_i < \beta$ for every $i$.  By
Lemma~\ref{lemma:alpha':inc}, $\adjsw{\beta_i} < \adjsw\beta$ for
every $i$, so $\adjmul{||F||} < \omega^{\adjsw\beta}$.  Using the
induction hypothesis, we obtain that
$||F^\circledast|| \leq \omega^{\adjmul{||F||}} =
\omega^{\omega^{\adjsw{\beta_1}} + \cdots + \omega^{\adjsw{\beta_n}}}
< \omega^{\omega^{\adjsw\beta}}$.

For every $n \in \nat$, $F \mid X^n$ is a closed subset of
$X^\circledast$, although not necessarily an m-product.  Indeed,
letting $C_1$, \ldots, $C_m$ be the components of $X$, $F \mid X^n$ is
equal to the finite union of the m-products
$F \mid C_{i_1}, \cdots, C_{i_n}$, where
$1 \leq i_1, \cdots, i_n \leq m$.  Then $f$ restricts to a continuous
map from $(F \mid X^n) \times X_\bot$ to $F \mid X^{n+1}$, and that
restriction is surjective.  By Lemma~\ref{lemma:stature:mono},
Theorem~\ref{thm:stature:prod}, and Proposition~\ref{prop:lifting},
$||F \mid X^{n+1}|| \leq ||(F \mid X^n) \times X_\bot|| = ||F \mid
X^n|| \otimes ||X_\bot|| = ||F \mid X^n|| \otimes (1+\alpha)$.  Since
$\beta \geq 1$, $\alpha = \omega^\beta$ is infinite, so
$1+\alpha=\omega^\beta$, and we obtain that
$||F \mid X^{n+1}|| \leq ||F \mid X^n|| \otimes \omega^\beta$.  By
induction on $n$,
$||F \mid X^n|| \leq ||F^\circledast|| \otimes \omega^{\beta \otimes
  n}$.  We recall that
$||F^\circledast|| \leq \omega^{\omega^{\adjsw{\beta_1}} + \cdots +
  \omega^{\adjsw{\beta_n}}}$.  For each $i$,
$\adjsw{\beta_i} < \adjsw\beta$, so
$\omega^{\adjsw{\beta_i}} < \omega^{\adjsw\beta}$.  Also,
$\beta \otimes n$ is the natural sum of $n$ copies of $\beta$, and
$\beta < \omega^{\adjsw\beta}$ by Lemma~\ref{lemma:alpha':omega}.
Since $\omega^{\adjsw\beta}$ is additively indecomposable and
$\oplus$-indecomposable,
$\omega^{\adjsw{\beta_1}} + \cdots + \omega^{\adjsw{\beta_n}} \oplus
(\beta \otimes n) < \omega^{\adjsw\beta}$.  It follows that
$||F \mid X^n|| < \omega^{\omega^{\adjsw\beta}} =
\omega^{\adjmul\alpha}$, for every proper closed subset $F$ of $X$ and
for every $n \in \nat$.

For every finite family of proper closed subsets $\bA_1$, \ldots,
$\bA_k$ of $X^\circledast$, there is a surjective continuous map
$g \colon \bA_1 + \cdots + \bA_k \to \bigcup_{i=1}^k \bA_i$, which
maps every $m \in \bA_i$ to itself.  By
Lemma~\ref{lemma:stature:mono},
$||\bigcup_{i=1}^k \bA_i|| \leq ||\bA_1 + \cdots + \bA_k|| =
\bigoplus_{i=1}^k ||\bA_i||$.

Every proper closed subset $\bA$ of $X^\circledast$ is a finite union
of m-products
$\bA_i \eqdef F_i \mid C_{i1}, \cdots, C_{in_i} \subseteq F_i \mid
X^{n_i}$.  Each of them is proper, so $F_i$ is a proper closed subset
of $X$.  It follows that
$||\bA|| \leq \bigoplus_{i=1}^k ||F_i \mid X^{n_i}||$.  Since each
term $||F_i \mid X^{n_i}||$ is strictly less than
$\omega^{\adjmul\alpha}$, and since the latter is
$\oplus$-indecomposable, $||\bA|| < \omega^{\adjmul\alpha}$.  By
Proposition~\ref{prop:limit},
$||X^\circledast|| \leq \omega^{\adjmul\alpha}$.

Finally, we deal with the general case, where $\alpha$ is not
necessarily of the form $\omega^\beta$.  We write $\alpha$ in Cantor
normal form as $\omega^{\alpha_1} + \cdots + \omega^{\alpha_m}$, with
$\alpha_1 \geq \cdots \geq \alpha_m$, $m \in \nat$.  By
Lemma~\ref{lemma:decomp:C:*}, there are closed subsets
$X=F_m \supseteq \cdots \supseteq F_1 \supseteq F_0 = \emptyset$ such
that $||F_i \diff F_{i-1}|| = \omega^{\alpha_i}$ for every $i$,
$1\leq i\leq m$.  Using Lemma~\ref{lemma:uplus}, the multiset union
map from
$(F_1 \diff F_0)^\circledast \times (F_2 \diff F_1)^\circledast \times
\cdots \times (F_m \diff F_{m-1})^\circledast$ to $X^\circledast$ is
continuous and surjective.  By Lemma~\ref{lemma:stature:mono} and
Theorem~\ref{thm:stature:prod},
$||X^\circledast|| \leq \bigotimes_{i=1}^m ||(F_i \diff
F_{i-1})^\circledast||$.  Since
$||F_i \diff F_{i-1}|| = \omega^{\alpha_i}$, we have already shown
that
$||(F_i \diff F_{i-1})^\circledast|| \leq
\omega^{\adjmul{\omega^{\alpha_i}}}$, so
$||X^\circledast|| \leq \omega^{\bigoplus_{i=1}^m
  \adjmul{\omega^{\alpha_i}}}$.  Now the $\adjmul\_$ operation
commutes with finite natural sums, as one can check by oneself, or by
referring to \cite[Notation~1]{vdMRW:bigVeblen}; so
$||X^\circledast|| \leq \omega^{\adjmul\alpha}$.  \qed

For every m-product $P \eqdef F \mid C_1, \cdots, C_m$, for every
$n \in \nat$ and every irreducible closed subset $C$ of $X$, we write
$P\mpl C^n$ for the m-product $F \mid C_1, \cdots, C_m, \allowbreak C^n$.
For every closed subset $\bA$ of $X^\circledast$, which one can write
as a finite union of m-products $P_1$, \ldots, $P_k$, we write
$\bA \mpl  C^n$ for the union of the m-products $P_1\mpl C^n$, \ldots,
$P_k\mpl C^n$.  Clearly, $\bA\mpl C^n$ is closed.  This is the set of
multisets that one can write as $m \uplus m'$, where $m \in \bA$ and
$m'$ is a multiset containing at most $n$ elements, all taken from
$C$.

\begin{lemma}
  \label{lemma:stature:mset:inc:FC}
  Let $F$ be a closed subset of a Noetherian space $X$, and $C$ be an
  irreducible closed subset of $X$ that is not included in $F$.  The
  map $(n, \bA) \mapsto (F \mid C^n) \cup (\bA\mpl C^{n+1})$ is a strictly
  monotonic map from
  $\nat \times_\lex (\Hoare_0 (F^\circledast) \diff
  \{F^\circledast\})$ to
  $\Hoare_0 ((F \cup C)^\circledast) \diff \{(F \cup
  C)^\circledast\}$.
\end{lemma}
\proof Let $f$ denote this map.  Let us also fix a point $x$ in
$C \diff F$.

Let $n \in \nat$ and $\bA$ be a proper closed subset of
$F^\circledast$.  It is clear that $f (n, \bA)$ is a closed subset
of $X^\circledast$ included in $(F \cup C)^\circledast$, hence a closed
subset of $(F \cup C)^\circledast$.  There is a multiset $m$ in
$F^\circledast \diff \bA$, and this allows us to form the multiset
$x^{n+1} \uplus m$ obtained by adding $n+1$ copies of $x$ to $m$.  It is
clear that $x^{n+1} \uplus m$ is not in $F \mid C^n$, since any element of
$F \mid C^n$ can contain at most $n$ elements of $C \diff F$.  We
claim that $x^{n+1} \uplus m$ is not in $\bA\mpl C^{n+1}$ either.  If it were,
$x^{n+1} \uplus m$ would be equal to the union $m_1 \uplus m_2$ of a multiset
$m_1$ in $\bA$, hence consisting of elements of $F$, plus a multiset $m_2$
of at most $n+1$ elements of $C$.  The $n+1$ copies of $x$ in
$x^{n+1} \uplus m$ cannot be in $m_1$, hence would comprise the totality of
$m_2$; this entails that $m_1 = m$, which is impossible since $m$ is
not in $\bA$.  It follows that $x^{n+1} \uplus m$ is not in
$f (n, \bA) = (F \mid C^n) \cup (\bA\mpl C^{n+1})$, showing that
$f (n, \bA)$ is a proper closed subset of $(F \cup C)^\circledast$.

If $n < n'$, then $(F \mid C^n) \subseteq (F \mid C^{n'})$, and
$(\bA\mpl C^{n+1}) \subseteq (F \mid C^{n'})$, since $\bA \subseteq
F^\circledast$ and $n+1 \leq n'$.  Therefore $f (n, \bA) \subseteq f
(n', \bA)$.  The inclusion is strict: given any $m \in F^\circledast
\diff \bA$, $x^{n+1} \uplus m$ is in $f (n', \bA)$ but not in $f (n, \bA)$.

It remains to verify that $\bA \subsetneq \bB$ implies
$f (n, \bA) \subsetneq f (n, \bB)$.  The inclusion is clear.  We pick
$m \in \bB \diff \bA$, and we note that $x^{n+1} \uplus m$ is in
$\bB\mpl C^{n+1} \subseteq f (n, \bB)$ but not in $f (n, \bA)$.  \qed

\begin{corollary}
  \label{corl:stature:mset:lbound:F}
  Let $F$ be a closed subset of a Noetherian space $X$, and $C$ be an
  irreducible closed subset of $X$ that is not included in $F$.  Then
  $||(F \cup C)^\circledast|| \geq ||F^\circledast|| \times \omega$.
\end{corollary}
\proof Let us call $f$ the map of
Lemma~\ref{lemma:stature:mset:inc:FC}.  For every $n \in \nat$ and
every proper closed subset $\bA$ of $F^\circledast$,
$\rk_{\Hoare_0 ((F \cup C)^\circledast)} (f (n, \bA)) \geq |\Hoare_0
(F^\circledast) \diff \{F^\circledast\}| \times n + \rk_{\Hoare_0
  (F^\circledast)} (\bA)$.  We note that
$|\Hoare_0 (F^\circledast) \diff \{F^\circledast\}| =
||F^\circledast||$, and that the least ordinal strictly larger than
all those ranks is $||F^\circledast|| \times \omega$.
\qed

A very similar argument shows the following.
\begin{lemma}
  \label{lemma:stature:mset:lbound:eps}
  Let $F$ be a closed subset of a Noetherian space $X$, such that
  $||X \diff F||$ is a critical ordinal $\epsilon$.  Then
  $||X^\circledast|| \geq ||F^\circledast|| \times \epsilon \times
  \omega$.
\end{lemma}
\proof
For every closed subset $C$ of $X$, for every $n \in \nat$, let $C^n$
denote the set of multisets of at most $n$ elements, all in $C$.  For
every irreducible closed subset $\bA$ of $F^\circledast$, let us write
$\bA \mpl  C^n$ for the set of multisets of the form $m \uplus m'$ where
$m \in \bA$ and $m' \in C^n$.  Writing $\bA$ as a finite union of
m-products $F_i \mid C_{i1}, \cdots, C_{in_i}$, $\bA \mpl  C^n$ is the
union of the m-products
$F_i \mid C_{i1}, \cdots, C_{in_i}, C, \cdots, C$, with $n$ copies of
$C$; so $\bA \mpl  C^n$ is closed.

For every irreducible closed subset $C$ of $X \diff F$, we write
$cl (C)$ for its closure in $X$.  Then $cl (C)$ is irreducible closed
in $X$ and $cl (C) \cap (X \diff F) = C$
\cite[Lemma~8.4.10]{JGL:topo}.  For every $n \in \nat$, we form the
set $F \mid X^n, cl (C)$ of all multisets $m \uplus m'$ where
$m \in F^\circledast$, and $m'$ contains at most $n$ elements, or
exactly $n+1$ of which one is in $cl (C)$.  Write $X$ as a finite
union of components $C_1$, \ldots, $C_m$, $F \mid X^n, cl (C)$ is
equal to the union of the m-products
$F \mid C_{i_1}, \cdots, C_{i_n}, cl (C)$, where $i_1$, \ldots, $i_n$
range over $\{1, \cdots, m\}$.  Hence $F \mid X^n, cl (C)$ is closed.

Let $Y \eqdef \Step (\Sober (X \diff F))$, and $Z$ be the poset
$\Hoare_0 (F^\circledast) \diff \{F^\circledast\}$ of all proper
closed subsets of $F^\circledast$.  For every $n \in \nat$, for every
$(C, C^+) \in Y$, for every $\bA \in Z$, we let
$f (n, (C, C^+), \bA) \eqdef (\bA \mpl  cl (C^+)^{n+1}) \cup (F \mid X^n, cl
(C))$.  We claim that $f$ is strictly monotonic from
$\nat \times_\lex Y \times_\lex Z$ to $\Hoare_0 (X^\circledast)$.  We
will write $x^n$ for the multiset containing exactly $n$ copies of
$x$.

If $n < n'$, then for all $(C, C^+)$ and $(C', {C'}^+)$ in $Y$, for
all $\bA, \bA' \in Z$, $\bA \mpl  cl (C^+)^{n+1}$ is included in
$F^\circledast \mpl  X^{n+1} = F\mid X^{n+1}$, hence in
$F \mid X^{n'}, cl (C')$; also, $F \mid X^n, cl (C)$ is included in
$F \mid X^{n'}$, hence in $F \mid X^{n'}, cl (C')$.  It follows that
$f (n, (C, C^+), \bA) \subseteq f (n', (C', {C'}^+), \bA')$.  With the
aim of showing that this inclusion is strict, we observe that $C$ is a
proper subset of $C^+$, so we can pick an element $x$ in
$C^+ \diff C$.  In particular, $x$ is in $C^+$, hence not in $F$.
Then $x^{n+2}$ is in $\bA' \mpl  cl (C^+)^{n'+1}$, and therefore in
$f (n', (C', {C'}^+), \bA')$.  We claim that it is not in
$f (n, (C, C^+), \bA)$.  It is not in $\bA \mpl  cl (C^+)^{n+1}$ since all
the multisets in that set are in $F^\circledast \mpl  cl (C^+)^{n+1}$, and
can therefore only contain at most $n+1$ elements of $C^+$.  If
$x^{n+2}$ were in $F \mid X^n, cl (C)$, then one of the $n+2$ copies
of $x$ would have to be in $F$, which is impossible.

Let us now fix $n \in \nat$, let $(C, C^+)$ and $(C', {C'}^+)$ be two
elements of $Y$ such that $(C, C^+) < (C', {C'}^+)$, namely such that
$C^+ \subseteq C'$, and let $\bA, \bA' \in Z$.  It is clear that
$F \mid X^n, cl (C)$ is included in $F \mid X^n, cl (C')$, since
$C \subsetneq C^+ \subseteq C'$.  Since $\bA \subseteq F^\circledast$
and $C^+ \subseteq C'$, $\bA \mpl  cl (C^+)^{n+1}$ is also included in
$F \mid X^n, cl (C')$, so
$f (n, (C, C^+), \bA) \subseteq f (n, (C', {C'}^+), \bA')$.  In order
to show that the inclusion is strict, let us pick an element $x$ of
$C^+ \diff C$, and let us note that $x$ is in $C'$ (and hence, not in
$F$).  Since $\bA$ is in $Z$, it is strictly included in
$F^\circledast$, so there is a multiset $m'$ in
$F^\circledast \diff \bA$.  We see that $x^{n+1} \uplus m'$ is in
$F \mid X^n, cl (C')$, hence in $f (n, (C', {C'}^+), \bA')$, and we
claim that it is not in $f (n, (C, C^+), \bA)$.  We first assume that
$x^{n+1} \uplus m'$ is in $\bA \mpl  cl (C^+)^{n+1}$, namely that it
splits as $m_1 \uplus m_2$ where $m_1 \in \bA$ and $m_2$ consists of
at most $n+1$ elements, all from $cl (C^+)$.  Since
$\bA \subseteq F^\circledast$ and $x$ is not in $F$, the $n+1$ copies
of $x$ cannot be part of $m_1$; so $m_2 = x^{n+1}$, and therefore
$m_1 = m'$.  This entails that $m'$ is in $\bA$, which is impossible.
Second, we assume that $x^{n+1} \uplus m'$ is in $F \mid X^n, cl (C)$.
Among the $n+1$ copies of $x$, none is in $F$, so that at least one is
in $cl (C)$.  Remembering that $x$ is not in $F$, $x$ would be in
$cl (C) \cap (X \diff F)$, which, as we have said earlier, is equal to
$C$.  This is impossible, since we have chosen $x$ in $C^+ \diff C$.

Finally, let $n \in \nat$, $(C, C^+) \in Y$, and let $\bA, \bA' \in Z$
be such that $\bA \subsetneq \bA'$.  It is clear that
$f (n, (C, C^+), \bA) \subseteq f (n, (C, C^+), \bA')$.  In order to
show that the inclusion is strict, let $m' \in \bA' \diff \bA$.  Let
also $x \in C^+ \diff C$ (hence $x \in X \diff F$, as before).  Then
$x^{n+1} \uplus m'$ is in $\bA' \mpl  cl ({C'}^+)^{n+1}$, hence in
$f (n, (C, C^+), \bA')$.  If it were in $F \mid X^n, cl (C)$, then
among the $n+1$ copies of $x$, since none is $F$, one would be in
$cl (C)$, hence in $cl (C) \cap (X \diff F) = C$, which is impossible.
If $x^{n+1} \uplus m'$ were in $\bA \mpl  cl (C^+)^{n+1}$, it would split
as $m_1 \uplus m_2$ with $m_1 \in \bA$ and $m_2 \in cl (C^+)^{n+1}$.
Since $\bA \subseteq F$, the $n+1$ copies of $x$ cannot be in $m_1$,
and are therefore in $m_2$.  It follows that $m_2 = x^{n+1}$, so that
$m_1 = m'$; but $m'$ is not in $\bA$, while $m_1$ is.  In any case, we
conclude that $x^{n+1} \uplus m'$ is not in $f (n, (C, C^+), \bA)$.

This finishes to prove that $f$ is strictly monotonic.

We claim that $\rsob {(X \diff F)} \geq \epsilon$. Otherwise,
$\rsob (X \diff F) < \epsilon$, so
$\omega^{\rsob (X \diff F)} < \omega^\epsilon$, and since
$\omega^\epsilon$ is $\oplus$-indecomposable,
$\omega^{\rsob (X \diff F)} \otimes n < \omega^\epsilon$ for every
$n \in \nat$.  By Proposition~\ref{prop:dim:stature}~(3), and letting
$n$ be the number of components of $X \diff F$,
$\epsilon = ||X \diff F|| \leq \omega^{\rsob {(X\diff F)}} \otimes n <
\omega^\epsilon = \epsilon$, which is impossible.

By Proposition~\ref{prop:Inc},
$|Y| = |\Sober (X \diff F)|-1 = \rsob {(X \diff F)}$, so
$|Y| \geq \epsilon$.  Since $f$ is strictly monotonic, for all
$n \in \nat$, $(C, C^+) \in Y$, and $\bA \in Z$, we have
$\rk_{\Hoare_0 (X^\circledast)} (f (n, (C, C^+), \bA)) \geq |Z| \times
(|Y| \times n + \rk_Y (C, C^+)) + \rk_Z (\bA)$.  Since
$|Y| \geq \epsilon$ and $|Z| = ||F^\circledast||$, we obtain that
$\rk_{\Hoare_0 (X^\circledast)} (f (n, (C, C^+), \bA)) \geq
||F^\circledast|| \times (\epsilon \times n + \rk_Y (C, C^+)) + \rk_Z
(\bA)$.  When $n$, $C$ and $\bA$ vary, the smallest ordinal strictly
larger than the latter is at least
$||F^\circledast|| \times \epsilon \times \omega$, so
$||X^\circledast|| \geq ||F^\circledast|| \times \epsilon \times
\omega$.  \qed

\begin{theorem}
  \label{thm:stature:mset}
  For every Noetherian space $X$,
  $||X^\circledast|| = \omega^{\adjmul\alpha}$, where
  $\alpha \eqdef ||X||$.
\end{theorem}
\proof Considering Proposition~\ref{prop:stature:mset:ubound}, it
suffices to show that $||X^\circledast|| \geq \omega^{\adjmul\alpha}$.
We do this by induction on $\alpha$.  When $\alpha=0$, $X$ is empty,
and there is exactly one element of $X^\circledast$, so
$||X^\circledast|| = 1  = \omega^0 = \omega^{\adjmul\alpha}$.

Otherwise, let us write $\alpha$ in Cantor normal form as
$\omega^{\alpha_1} + \cdots + \omega^{\alpha_m}$ with
$\alpha \geq \alpha_1 \geq \cdots \geq \alpha_m$, $m \geq 1$.
Let $\beta \eqdef \omega^{\alpha_1} + \cdots + \omega^{\alpha_{m-1}}$.

If $\alpha_m=0$, then $||X|| = \beta+1$.  There is a proper closed
subset $F$ of $X$ such that $\rk_{\Hoare_0 X} (F) = \beta$, namely
such that $||F|| = \beta$, by Lemma~\ref{lemma:rk:C}.  By
Lemma~\ref{lemma:+1}, we can write $X$ as $F \cup C$, where
$C \eqdef \dc x$, and $x$ is a point of $X \diff F$.  By
Corollary~\ref{corl:stature:mset:lbound:F}, $||X^\circledast||$ is
then larger than or equal to $||F^\circledast|| \times \omega$, hence
to $\omega^{\adjmul\beta+1}$, using the induction hypothesis.  We now
note that $\adjmul\beta+1 = \adjmul\alpha$.

If $\alpha_m \geq 1$, then $\omega^{\alpha_m}$ is a limit ordinal.
For every $\gamma < \omega^{\alpha_m}$, $\alpha' \eqdef \beta+\gamma$
is strictly smaller than $\alpha$, so there is a proper closed subset
$F$ of $X$ such that $\rk_{\Hoare_0 X} (F) = \alpha'$.  By induction
hypothesis, $||F^\circledast|| \geq \omega^{\adjmul{\alpha'}}$.  Since
$\gamma < \omega^{\alpha_m}$, the Cantor normal form of $\alpha'$ is
of the form
$\omega^{\alpha_1} + \cdots + \omega^{\alpha_{m-1}}+ \omega^{\gamma_1}
+ \cdots + \omega^{\gamma_k}$, with
$\alpha_m > \gamma_1 \geq \cdots \geq \gamma_k$.  Using
Lemma~\ref{lemma:alpha':inc}, the Cantor normal form of
$\adjmul{\alpha'}$ is equal to
$\omega^{\adjsw{\alpha_1}} + \cdots + \omega^{\adjsw{\alpha_{m-1}}}+
\omega^{\adjsw{\gamma_1}} + \cdots + \omega^{\adjsw{\gamma_k}}$, which
is larger than or equal to $\adjmul\beta + \gamma$.  In particular,
$||X^\circledast|| \geq ||F^\circledast|| \geq
\omega^{\adjmul\beta+\gamma}$.  By letting $\gamma$ vary among the
ordinals strictly smaller than $\omega^{\alpha_m}$ and taking suprema,
$||X^\circledast|| \geq \omega^{\adjmul\beta + \omega^{\alpha_m}}$.

When 
$\alpha_m$ is not of the form $\epsilon+n$ with $\epsilon$ a critical
ordinal and $n \in \nat$, $\alpha_m = \adjsw{\alpha_m}$, and therefore
we have obtained the desired lower bound
$||X^\circledast|| \geq \omega^{\adjmul\alpha}$.

If $\alpha_m=\epsilon$ for some critical ordinal $\epsilon$, then we
use Lemma~\ref{lemma:stature:mset:lbound:eps} instead.  There is a
proper closed subset $F$ of $X$ such that
$\rk_{\Hoare_0 X} (F) = \beta$, by Lemma~\ref{corl:decomp:C}, and
$||X \diff F|| = \omega^{\epsilon} = \epsilon$.  By induction
hypothesis, $||F^\circledast|| \geq \omega^{\adjmul\beta}$.  By
Lemma~\ref{lemma:stature:mset:lbound:eps},
$||X^\circledast|| \geq ||F^\circledast|| \times \epsilon \times
\omega \geq \omega^{\adjmul\beta} \times \epsilon \times \omega$.  Now
$\epsilon \times \omega = \omega^\epsilon \times \omega =
\omega^{\epsilon+1} = \omega^{\adjsw{\alpha_m}}$, so
$||X^\circledast|| \geq \omega^{\adjmul\beta} \times
\omega^{\adjsw{\alpha_m}} = \omega^{\adjmul\alpha}$.

If $\alpha_m=\epsilon+n$ for some critical ordinal $\epsilon$ and some
$n \in \nat$ such that $n \geq 1$,
$\omega^{\alpha_m} = \epsilon \times \omega^n$ is the supremum of
ordinals of the form $\epsilon \times \omega^{n-1} \times k$,
$k \in \nat$.  Given
$\gamma \eqdef \epsilon \times \omega^{n-1} \times k =
\omega^{\epsilon+n-1} \times k$, $\alpha' \eqdef \beta+\gamma$ is
strictly smaller than $\alpha$, so there is a proper closed subset $F$
of $X$ such that $\rk_{\Hoare_0 X} (F) = \alpha'$.  By induction
hypothesis,
$||F^\circledast|| \geq \omega^{\adjmul{\alpha'}} =
\omega^{\adjmul\beta + \omega^{\epsilon+n} \times k}$.  By taking
suprema over $k$,
$||X^\circledast|| \geq \omega^{\adjmul\beta + \omega^{\epsilon+n+1}}
= \omega^{\adjmul\beta + \omega^{\adjsw{\alpha_m}}} =
\omega^{\adjmul\alpha}$.  \qed

\begin{table}
  \centering\footnotesize
  \begin{tabular}{cclcl}
    \hline
    $X$ & $\sob X$ && $||X||$ & \\
    \hline \hline
    Finite $T_0$ & $\leq \card X$ && $\card X$
                              & Lem.~\ref{lemma:finite}
    \\
    \hline
    Ordinal $\alpha$ $\rlap{(Alex.)}$
        & $\alpha$ $/$ $\alpha+1$ 
                   & Lem.~\ref{lemma:ordinal}
                              & $\alpha$ & Lem.~\ref{lemma:ordinal} \\
    \hline
    Ordinal $\alpha$ $\rlap{(Scott)}$
        & $\alpha$ $/$ $\alpha+1$ 
                   & Lem.~\ref{lemma:ordinal}
                              & $\alpha$ $/$ $\alpha-1$ 
                              & Lem.~\ref{lemma:ordinal} \\
    \hline
    Cofinite $\rlap{topology}$ & 
                     $1$ $/$ $2$
                   & Thm.~\ref{thm:dim:cofinite}
                              & $\min (\card X, \omega)$
                              & Thm.~\ref{thm:stature:cofinite} \\
    \hline
    \hline
    $X+Y$ & $\max (\sob X, \sob Y)$
                   & Prop.~\ref{prop:dim:sum}
                              & $||X|| \oplus ||Y||$
                              & Prop.~\ref{prop:stature:sum} \\
    \hline
    $X +_\lex Y$ & $\sob X + \sob Y$
                   & Prop.~\ref{prop:dim:+lex}
                              & $||X||+||Y||$
                              & Prop.~\ref{prop:stature:+lex} \\
    \hline
    $X_\bot$ & $1+\sob X$ & Prop.~\ref{prop:lifting}
                              & $1+||X||$ & Prop.~\ref{prop:lifting} \\
    \hline
    $X \times Y$ & $(\sob X \oplus \sob Y)-1$ 
                   & Prop.~\ref{prop:dim:prod}
                              & $||X|| \otimes ||Y||$
                              & Thm.~\ref{thm:stature:prod} \\
    \hline
    $\HVz X$, $\Hoare_\fin X$,
        & $||X||+1$
                   & Thm.~\ref{thm:dim:H}
                              & $\geq 1+||X||$,
                              & Prop.~\ref{prop:stature:H} \\
    $\pow X$, $\pow_\fin X$ &&&  $\leq \omega^{||X||}$ & \\
    \hline
    $X^*$ & $\omega^{\adjsw{||X||}}+1$ 
                   & Thm.~\ref{thm:dim:X*}
                              & $\omega^{\omega^{||X||'}}$ 
                 & Thm.~\ref{thm:stature:X*} \\
    & \multicolumn{2}{l}{($\adjsw\alpha \eqdef \alpha+1$ if
      $\alpha=\epsilon+n$,
      $\epsilon$ critical, $n \in \nat$,}
                   & \multicolumn{2}{l}{($\alpha' \eqdef \alpha-1$ if $\alpha$ finite,}
    \\
        &  $\alpha$ otherwise)
                   &&  $\adjsw\alpha$
        otherwise) & \\
    \hline
    $\bigtriangleright_{n=1}^{+\infty} X_n$
        & $\bigoplus_{n=1}^{+\infty} \rsob {X_n} + 1$ $/$
                   & Thm.~\ref{thm:dim:prefix}
                              & $\bigotimes_{n=1}^{+\infty} ||X_n||$ $/$
                              & Thm.~\ref{thm:stature:prefix} \\
        & $\bigoplus_{n=1}^k \rsob {X_n} + \omega+1$ &
                   & $\bigotimes_{m=1}^k ||X_m|| \times \omega$
                              & \\
    \hline
    $X^\triangleright$ & $\omega^{\alpha_1+1}+1$ 
                   & Cor.~\ref{corl:dim:prefix:homog}
                              & $\omega^{\omega^{\beta_1+1}}$ $/$
                                $\omega$ 
                              & Cor.~\ref{corl:stature:prefix:homog} \\
        & \multicolumn{2}{l}{where
          $\sob X-1 =_{\text{CNF}} \omega^{\alpha_1} + \cdots$}
                   & \multicolumn{2}{l}{where
                     $||X|| =_{\text{CNF}} \omega^{\alpha_1} +
                     \cdots$,} \\
    &&& \multicolumn{2}{l}{\qquad\qquad
                      $\alpha_1 =_{\text{CNF}} \omega^{\beta_1}+ \cdots$}
    \\
    \hline
    $X^\circledast$
        & $\geq (\omega \times ||X||)+1$,
          & Prop.~\ref{prop:dim:mset:ubound},
                              & $\omega^{\adjmul\alpha}$
                              & Thm.~\ref{thm:stature:mset} \\
    & $\leq  (||X|| \otimes \omega) +1$ &
                                          Prop.~\ref{prop:dim:mset:lbound}
                              & \multicolumn{2}{l}{($\adjmul\alpha \eqdef \omega^{\adjsw{\alpha_1}} + \cdots +
        \omega^{\adjsw{\alpha_m}}$} \\
    &&& \multicolumn{2}{l}{\strut if $\alpha =_{\text{CNF}} \omega^{\alpha_1} +
        \cdots + \omega^{\alpha_m}$)} \\
    \hline
  \end{tabular}
  \caption[Statures and sobrification ranks]{Statures and sobrification ranks of
    Noetherian constructions (all spaces assumed non-empty)}
  \label{tab:summary}
\end{table}

\section{Conclusion and Open Problems}
\label{sec:bf-conclusion}

We have developed a theory of statures of Noetherian spaces that
generalizes and extends the theory of maximal order types of wpos.  In
the process, we have also studied the related notion of sobrification rank.  We
have also given an extensive list of explicit formulae for sobrification ranks
and statures of several families of Noetherian spaces, arising or not
from wqos (see Table~\ref{tab:summary}, where only the cases of
non-empty spaces $X$, $Y$, $X_n$ are shown, in order to avoid a
proliferation of cases).  Among the questions that remain, let us cite
the following.
\begin{enumerate}
\item We have $1+||X|| \leq ||\HVz X|| \leq \omega^{||X||}$
  (Proposition~\ref{prop:stature:H}), and both the lower bounds and
  upper bounds are attained.  What is the exact set of ordinals
  between those bounds that one can obtain as $||\HVz X||$? Can we
  reach the upper bound $\omega^{||X||}$ for all infinite values that
  $||X||$ may take?
\item Similarly, the sobrification rank of $X^\circledast$ lies between
  $(\omega \times ||X||)+1$ (Proposition~\ref{prop:dim:mset:lbound})
  and $(||X|| \otimes \omega) +1$
  (Proposition~\ref{prop:dim:mset:ubound}), and those bounds are
  attained.  What other values can $\sob {X^\circledast}$ evaluate to?
\item There are many other examples of Noetherian spaces we have not
  considered.  An outstanding one is the space of finite trees with
  function symbols taken from $X$ \cite[Section~10]{FGL:partI}, for
  which we expect the sobrification rank and stature to be equal to Schmidt's
  formula $f^+ {\binom {||X||} \omega}$, as in the wqo case
  \cite[Chapter~II]{Schmidt:maxot}.
\item For a Noetherian space $X$, the spaces of infinite words
  $X^\omega$, and of finite-or-infinite words $X^{\leq \omega}$ are
  Noetherian as well \cite{GL:infwords}.  (This would fail with wqos
  instead of Noetherian spaces, and is similar to a well-known result
  on bqos due to Nash-Williams \cite{NW:bqo}.)  \blue{This is more
    generally true of the spaces $X^{<\alpha}$ of transfinite words
    over $X$ of ordinal length $< \alpha$, where $\alpha$ is an
    arbitrarily large ordinal \cite{GLHL:words:II}.}  What are the
  sobrification ranks and statures of those spaces?  \blue{An upper
    bound is given \cite[Proposition~33]{GLHL:words:II}.}
\item What about Noetherian spaces obtained as spectra of Noetherian
  rings?  There is an abundant literature on Krull dimension of
  Noetherian rings and variations, see
  \cite{Bass:krull,Gullliksen:length,Rhodes:krull} for example.  It is
  not yet clear to us what the precise relationship to our notion of
  sobrification rank is, and whether one can derive corresponding
  results on statures in general.  For the Noetherian ring of
  polynomials on $m$ variables over a field $K$, the reduced
  sobrification rank of its spectrum is exactly its Krull dimension,
  namely $m$.  The ordinal rank of the family of all ideals in that
  ring is $\omega^m+1$ \cite[Lemma~3.18]{AP:monomial};
  when $K$ is algebraically closed, this is not far from the stature
  of the spectrum, which is the ordinal rank of the subfamily of
  \emph{radical} ideals, minus $1$, by Hilbert's Nullstellensatz.
\item An application of the theory of maximal order types in computer
  science consists in evaluating the precise complexity of
  reachability and related questions on well-structured transition
  systems, see \cite{FFSS:ack,SS:higman} for example.  The present
  theory should find similar complexity-theoretic applications in
  relation with the topological well-structured transition systems of
  \cite{JGL-icalp10}, and this remains to be developed.
\end{enumerate}


\section*{Acknowledgments}
\label{sec:acknowledgments}

A first version of this paper was submitted to the Pacific Journal of
Mathematics.  While it was rejected there, mostly for reasons of
length, the referee did an extraordinary job of checking every result,
every proof, proper terminology, and even more minor details, in
almost no time.  We owe her or him the discovery of some flaws in the
previous proofs of Theorem~\ref{thm:dim:prefix},
Lemma~\ref{lemma:stature:mset:lbound:eps}, and
Lemma~\ref{lemma:stature:X*:ubound:finite}, and we are grateful to her
or him for that.  We also thank the reviewer of this version, who made
us clarify quite a few points, who argued for a better name for what
we had variously called ``almost open'', then ``initial'' maps, and
who suggested the adjective ``full''.  Finally, we thank Maurice
Pouzet for pointing us to the origins of the discovery of the
coincidence $o (P) = ||P||$ for wpos $P$.


\bibliographystyle{plain}
\DeclareRobustCommand{\VAN}[3]{#3}
\ifarxiv

\else
\bibliography{stature}
\fi

\end{document}